\documentclass[preprint,11pt]{elsarticle}

\usepackage{amssymb}
\usepackage[english]{babel}
\usepackage{graphicx}
\usepackage{subcaption}
\usepackage[T1]{fontenc}
\usepackage[utf8x]{inputenc}
\usepackage{color}

\usepackage{algorithm}
\usepackage{algpseudocode}
\usepackage{tikz}
 \usepackage{pgfplots}
 \tikzset{math3d/.style= {x= {(-0.353cm,-0.353cm)}, z={(0cm,1cm)},y={(1cm,0cm)}}}
\usepackage[colorlinks=true]{hyperref}
\usepackage{amsmath}
\usepackage{amssymb, amsthm,mathrsfs,cases}
\usepackage{multirow}

\usepackage[most]{tcolorbox}
\usepackage{array,tabularx}
\usepackage{colortbl}
\newcolumntype{Y}{>{\raggedleft\arraybackslash}X}
\tcbset{enhanced,fonttitle=\bfseries\large,fontupper=\normalsize\sffamily,
colback=white,colframe=gray!50!black,colbacktitle=white,
coltitle=black,center title}

\usepackage[text={6.5in,10.1in}]{geometry}
\newtheorem{theorem}{Theorem}
\newtheorem{corollary}{Corollary}
\newtheorem{lemma}{Lemma}

\newtheorem{remark}{Remark}

\newcommand{\bm}[1]{\mbox{\boldmath $#1$}}

\journal{Journal of Scientific Computing}

\begin{document}

\begin{frontmatter}

\title{\small{\bf Stabilization of the Gradient Method for Solving Linear Algebraic Systems}\footnote{This research was supported by the NSERC Canadian grant number RGPIN-2024-05291. I confirm that this work is original and has not been published elsewhere, nor is it currently under consideration for publication elsewhere.} \\ {\it A Method Related to the Normal Equation}}

\author{Ibrahima Dione}

\address{Universit{\'e} de Moncton, D{\'e}partement de math{\'e}matiques
         et de statistique, \\
         18 Antonine-Maillet Ave, Moncton, NB, E1A 3E9, Canada,  \\
         ibrahima.dione@umoncton.ca.}

\begin{abstract} Although it is relatively easy to apply, the gradient method often displays a disappointingly slow rate of convergence. Its convergence is specially based on the structure of the matrix of the algebraic linear system, and on the choice of the stepsize defining the new iteration. We propose here a simple and robust stabilization of the gradient method, which no longer assumes a structure on the matrix (neither symmetric, nor positive definite) to converge, and which no longer requires an approach on the choice of the stepsize. We establish the global convergence of the proposed stabilized algorithm under the only assumption of nonsingular matrix. Several numerical examples illustrating its performances are presented, where we have tested small and large scale linear systems, with and not structured matrices, and with well and ill conditioned matrices. 
\end{abstract}

\end{frontmatter}

\section{Introduction}
This paper is devoted to the solving of the linear system  
\begin{align}\label{Eq1}
    \mathbf{A}\bm{x} = \bm{b},
\end{align}
by means of a stabilized gradient iterative method, where $\mathbf{A} \in \mathbb{R}^{n\times n}$ is a $n\times n$ matrix (which is ill-conditioned or not) and $\bm{b} \in \mathbb{R}^n$ is a given $n$-vector. The only requirement on this matrix is that it is nonsingular (whatever the structure of the eigenvalues), the exact solution being thus given as follows $\bm{x}^{\star} = \mathbf{A}^{-1}\bm{b} \in \mathbb{R}^n$. Throughout this work, $\|\cdot\|_{n}$ denotes the euclidean vector norm in $\mathbb{R}^n$, $\left<\cdot, \cdot\right>_{n}$ its associated scalar product, and $\|\cdot\|$ is the induced matrix norm to $\|\cdot\|_{n}$. 

The linear system \eqref{Eq1} is the standard formulation of many problems coming from different fields in science and engineering whose economy, chemistry, physics, models of partial differential equations, ect. The solving of this linear system is still challenging, and also represents one of the most time consuming part of the numerical simulation process. We may distinguish two main methods for its solving: Direct and iterative methods. The former methods, which are widely used in many industrial codes, are based on the factorization of the matrix $\mathbf{A}$ into easily invertible matrices and require a predictable amount of ressources in time and storage \cite{AllaireKaber08,BjÄorck96,QuarteroniSaccoSaleri07}. Unfortunately, these methods are inefficient when the problem size becomes large since the factorization becomes impractical due to lack of time or storage capacity. The use of the latter methods, the {\it iterative methods}, becomes thus mandatory. Iterative methods require fewer storage, and often, fewer operations than direct methods during the solving process \cite{AllaireKaber08,BjÄorck96,QuarteroniSaccoSaleri07,Hackbusch10,SibonyMardon88}. However, these methods may fail to converge in a reasonable amount of time, and sometimes may not have the same reliability than direct methods. 

The gradient method which is traced back to A. L. Cauchy \cite{Cauchy47}, is one of the most straightforward iterative method to build an approximation solution of the linear system \eqref{Eq1}. A simple way to establish this iterative method is to reformulate equation \eqref{Eq1} under the linear fixed-point equation $\bm{x} = \bm{x} - \alpha(\mathbf{A}\bm{x} - \bm{b})$, from which the following iterative scheme is built
\begin{align}\label{Eq2}
      \bm{x}^{[k+1]} = \bm{x}^{[k]} - \alpha_{k}(\mathbf{A}\bm{x}^{[k]} - \bm{b}), \;\; k = 0, 1, 2, \cdots,
\end{align}  
where $\alpha_{k}$ is a parameter (a stepsize) to be defined. The gradient method \eqref{Eq2} can be written under the iterative form (see \cite{AllaireKaber08} - Definition $9.1.1$)
\begin{align}\label{LIM1}
   \mathbf{M}_{\alpha_{k}}\bm{x}^{[k+1]} = \mathbf{N}_{\alpha_{k}}\bm{x}^{[k]} + \bm{b},
\end{align}
where the regular decomposition of the matrix $\mathbf{A} = \mathbf{M}_{\alpha} - \mathbf{N}_{\alpha}$, with $\mathbf{M}_{\alpha} = \frac{1}{\alpha}\mathbf{I}$ and $\mathbf{N}_{\alpha} = \left(\frac{1}{\alpha}\mathbf{I} - \mathbf{A}\right)$, is made. In the case where the matrix $\mathbf{A}$ is symmetric, the gradient method is based on the natural idea of choosing the direction of the greatest rate of change of the quadratic functional 
\begin{align*}
    J(\bm{x}) = \frac{1}{2}\left<\mathbf{A}\bm{x}, \bm{x}\right>_{n} - \left<\bm{b}, \bm{x}\right>_{n}, 
\end{align*}
namely, the direction opposite to the gradient $\nabla J(\bm{x}) = \mathbf{A}\bm{x} - \bm{b}$. A process of finding the appropriate stepsize $\alpha_{k}$ is also needed. Many iterative algorithms, typically called {\it line search algorithms}, have been proposed for the choice of this parameter. The most known are the {\it constant stepsize},  the {\it exact line search stepsize} and the {\it backtracking stepsize} (see \cite{Beck14} - Chap. 4). The main advantage of the constant stepsize strategy is of course its simplicity, but its drawback is how to fixe its value since a large stepsize can cause a non converging gradient method, while a small stepsize can slow the convergence.\\ 
As for the exact line search stepsize, even in the simplest case where the marix $\mathbf{A}$ is symmetric and positive definite, the gradient method may not be very efficient \cite{Akaike59}. With the exact line search, the gradient method converges linearly and behaves increasingly badly when the condition number of the matrix deteriorates. Attempts to increase the performance of this method have also been investigated \cite{ForsytheMotzkin51,Humphrey66,Schinzinger66}. The exact line search method is mainly applied in the special cases of quadratic problems, where the stepsize $\alpha_{k}$ is analytically computed. Otherwise, its approximation is used, generating the so called {\it inexact line search} procedures widely investigated in the litterature \cite{Armijo66,Goldstein65,Wolfe68,Powell76,DennisSchnabel83,Fletcher87,PotraShi95,Lemaréchal81,MoréThuente90}. The {\it backtracking stepsize} is a version of this approximation. \\
Later on, Barzilai and Borwein \cite{BarzilaiBorwein88} proposed a surprising superlinear convergence result of the steepest descent method, based on a different strategy for choosing the step lenght. For the simplest case of symmetric and positive definite matrix $\mathbf{A}$ with $n = 2$, Barzilai and Borwein show that the gradient method \eqref{Eq2}, with the stepsize $\alpha_{k}$ defined from informations of the previous iteration, converges superlinearly. Their works have also been the origin of many others researches on the steepest descent method. Raydan \cite{Raydan93,Raydan97} established convergence of the Barzilai and Borwein's approach in the quadratic case by providing numerical results, whereas Fletcher \cite{Fletcher90} analysed its links with the spectrum of the matrix $\mathbf{A}$. The possibility to apply the Barzilai and Borwein's approach to constrained optimization problems was considered by Friedlander et al. \cite{FriedlanderMartinezRaydan95} (in the case of quadratic functionals), and by Birgin et al. \cite{BirginMartinezRaydan00}. Others more rescent investigations on the Barzilai and Borwein's approach have also been made \cite{DaiYuanYuan02,DaiYuan03}, and some explained why it is competitive to the conjugate gradient based method \cite{Fletcher05}. More recently, under a quadratic functional, Yuan \cite{Yuan06} proposed a new computational process of the stepsize $\alpha_{k}$, enabling a fast convergence of the steepest descent method. Moreover, its approach also possesses the monotoniticity property that permits a generalization to non-linear functions.   \\      
Despite these performances on the gradient method, it is supplanted by the use of the conjugate gradient method or its variants. In the solving of the linear system \eqref{Eq1}, the conjugate gradient method \cite{HestenesStiefel52} and some of its preconditionings \cite{GolubLoan13} remain still the methods of choices. The main feature of these methods is that, they are able to deliver a reasonably accurate approximation solution in few iterations, specially when the matrix has certain eigenvalue structure.

We propose a simple and robust stabilization of the gradient method \eqref{Eq2}, that may be seen as a preconditioned iterative method, defined as follows  
\begin{align*}
    \left(\mathbf{I} + \gamma\mathbf{A}^{t}\mathbf{A}\right)\bm{x}^{[k+1]}_{\gamma} & = \left(\mathbf{I} - \alpha_{k}\mathbf{A}\right)\bm{x}^{[k]}_{\gamma} + \alpha_{k}\bm{b} + \gamma\mathbf{A}^{t}\bm{b},  \; \forall \; k = 0, 1, 2, \cdots, 
\end{align*}
where $\gamma > 0$ is a stabilization parameter which can be chosen as large as possible, and independently of the linear system data and the stepsize parameter $\alpha_{k}$. This stabilized iterative method is obtained from \eqref{Eq2}, through the idea that we need to incorporate terms that gradually control the residual of the iterations. One of the nice properties of this iterative scheme is that we no longer need to assume that the matrix $\mathbf{A}$ is positive definite, nor that it is symmetric, but only that it is nonsingular in order to the linear system \eqref{Eq1} to have a unique solution. We also show that the choice of the stepsize parameter $\alpha_{k}$ is no longer necessary, since this stabilized method converges whatever $\alpha_{k} \in \mathbb{R}$ (positive or negative) and through the established error estimate 
\begin{align*}
     \|\bm{x}^{[k]}_{\gamma} - \bm{x}^{\star}\|_{n} \leq \left(\frac{\kappa\left(\mathbf{V}\right)}{1 + \gamma\sigma^{2}_{n}}\right)^{k}\left(\prod_{i=0}^{k-1}\left\|\mathbf{I} - \alpha_{i}\mathbf{A}\right\|\right)\|\bm{x}^{[0]}_{\gamma} - \bm{x}^{\star}\|_{n}, \, \forall \, \gamma > 0, \, \forall \, k \in \mathbb{N}^{\star},
\end{align*}
where $\kappa\left(\mathbf{V}\right) =\|\mathbf{V}\|\|\mathbf{V}^{-1}\|$ is the condition number of the orthogonal matrix $\mathbf{V} \in \mathbb{R}^{n\times n}$ coming from the Singular Value Decomposition (SVD) of $\mathbf{A} = \mathbf{U}\mathbf{\Sigma}\mathbf{V}^{t}$, and $\sigma_{n}$ is the smallest singular value of $\mathbf{A}$. Whereas, taking a constant stepsize $\alpha_{i} = \alpha$ this estimate simplifies to
\begin{align*}
    \|\bm{x}^{[k]}_{\gamma} - \bm{x}^{\star}\|_{n} \leq \left(\left(\frac{\kappa\left(\mathbf{V}\right)}{1 + \gamma\sigma^{2}_{n}}\right)\left\|\mathbf{I} - \alpha\mathbf{A}\right\|\right)^{k}\|\bm{x}^{[0]}_{\gamma} - \bm{x}^{\star}\|_{n}, \, \forall \, \gamma > 0, \, \forall \, k \in \mathbb{N}^{\star}.
\end{align*}

This work is organized as follows: we introduce in section \ref{sec2}, a stabilized gradient method built from the gradient method. We then present in section \ref{sec3}, the convergence analysis of the stabilized method and propose a pseudo-algorithm that underlines the differents steps of its solving. Several numerical examples are presented in section \ref{sec5} to illustrate the theoretical results. We begin by analysing well-conditioned problems (subsections \ref{NEx1} - \ref{NEx2}), specially examples in the litterature whose the convergence is clearly proved to be slow or impossible. We specially investigated the cases of rank-deficient and severely ill-conditioned matrix through both small and large scale problems (subsection \ref{NEx3}). Through these examples, we also investigate the accuracy and stability of this approach with respect to the stepsize parameter $\alpha_{k}$ and to the stabilization parameter $\gamma$. We also note that the stabilized method may behave as a direct method, when the stabilization parameter $\gamma$ is chosen greater than the desired accuracy $\varepsilon$ (namely $\gamma > \frac{1}{\varepsilon}$).

\section{Stabilization of the Gradient Method}\label{sec2}
\subsection{The Stabilized Gradient Method}
When using the gradient method in the solving of the linear system \eqref{Eq1}, it is important to properly choose the stepsize parameter $\alpha_{k}$. Although clear recommandations have been made on this choice (see \cite{AllaireKaber08} - Chap. $9$ - Theorem $9.1.1$, \cite{SibonyMardon88} - Chap. $1$ - Theorem $4.1$), we may also face in some applications non convergence of this iterative method. The aim of this work is, from the formulation \eqref{Eq2}, to define a stable iterative method that not only solves convergence difficulties but also behaves much like a direct method since we may be able to control its iterations. This method will be a very effective technique for the solving of differential and integral equations.

The iterative method \eqref{Eq2} can also be written under the following minimization problem: Given $\bm{x}^{[k]}$, compute $\bm{x}^{[k+1]} \in \mathbb{R}^{n}$ from the following minimization problem 
\begin{align}\label{Eq4}
   \mathcal{F}\left(\bm{x}^{[k+1]}\right) := \min_{\bm{x} \in \mathbb{R}^{n}} \mathcal{F}(\bm{x}),
\end{align}
where the functional $\mathcal{F}(\cdot)$ is defined under the following form 
\begin{align}\label{Eq41}
    \mathcal{F}(\bm{x}) = \frac{1}{2}\left\|\bm{x} - \bm{x}^{[k]}\right\|^{2}_{n} + \alpha_{k}\left<\mathbf{A}\bm{x}^{[k]} - \bm{b}, \,\bm{x} - \bm{x}^{[k]}\right>_{n},
\end{align}
The approach we propose here is formulated from the problem \eqref{Eq4}, with the additional worry to obtain stability. For that, we add the equation \eqref{Eq1} as a constraint of the minimization problem \eqref{Eq4} by considering the following affine constrained problem: Given $\bm{x}^{[k]}$, compute  $\bm{x}^{[k+1]} \in \bm{K}$ by
\begin{align}\label{Eq5}
   \mathcal{F}\left(\bm{x}^{[k+1]}\right) := \min_{\bm{x} \in \bm{K}} \mathcal{F}(\bm{x}), 
\end{align}
where the constraint space is defined as $\bm{K} = \left\{\bm{x} \in \mathbb{R}^n\mid \mathbf{A}\bm{x} = \bm{b}\right\}$. Instead of working in this constrained set, we relaxe the constraint by applying the regularization method to the minimization problem \eqref{Eq5}. We thus interest to the regularized problem: Given $\bm{x}^{[k]}_{\gamma}$, compute $\bm{x}^{[k+1]}_{\gamma} \in \mathbb{R}^n$ by
\begin{align}\label{Eq52}
    \mathcal{F}_{\gamma}(\bm{x}^{[k+1]}_{\gamma}) := \min_{\bm{x} \in \mathbb{R}^n} \mathcal{F}_{\gamma}(\bm{x}), \; \forall \; \gamma > 0,
\end{align}
where the regularized functional $\mathcal{F}_{\gamma}(\cdot)$ is defined here as follows  
\begin{align}\label{Def34}
    \mathcal{F}_{\gamma}(\bm{x}) := \mathcal{F}(\bm{x}) + \frac{\gamma}{2}\left\|\mathbf{A}\bm{x} - \bm{b}\right\|^{2}_{n}, \; \forall \; \gamma > 0,
\end{align}
and where $\gamma$ is a regularisation parameter which controls the weight given to the residual $\left\|\mathbf{A}\bm{x} - \bm{b}\right\|^{2}_{n}$, relative to the minimisation of the functional $\mathcal{F}(\cdot)$. Taking the first order necessary optimality condition of the problem \eqref{Eq52} by differentiating the functional $\mathcal{F}_{\gamma}(\cdot)$ with respect to $\bm{x}$ and setting it equals to zero at $\bm{x}^{[k+1]}_{\gamma}$, we obtain the iterative method: Given $\bm{x}^{[k]}_{\gamma}$, compute $\bm{x}^{[k+1]}_{\gamma} \in \mathbb{R}^n$ by
\begin{align}\label{EqP1}
    \left(\mathbf{I} + \gamma\mathbf{A}^{t}\mathbf{A}\right)\bm{x}^{[k+1]}_{\gamma} = \left(\mathbf{I} - \alpha_{k}\mathbf{A}\right)\bm{x}^{[k]}_{\gamma} + \alpha_{k}\bm{b} + \gamma\mathbf{A}^{t}\bm{b},  \, \forall \, \gamma > 0, \; \forall \; k = 0, 1, 2, \cdots.
\end{align}

Comparing the gradient method \eqref{Eq2} to the stabilized gradient method \eqref{EqP1}, we may see that the latter is obtained from the former by simply adding the terms $\gamma\mathbf{A}^{t}\mathbf{A}$ and $\gamma\mathbf{A}^{t}\bm{b}$ in the left and right hand sides, respectively. These added terms play a stabilization role to the gradient method by controlling the residual, as the iterations go along. 
\begin{remark} \text{}  
\begin{itemize}
    \item The only difference between the gradient method \eqref{Eq2} and the stabilized iterative method \eqref{EqP1} is the residual that is taken account as a stabilizer of the iterative problem. This residual plays the role of a stabilizer since its goal is not to ensure the uniqueness of solution, but to bring stability to the whole system by permitting the convergence of the method in few iterations.
    \item What we mean from the formulation \eqref{EqP1} is, when solving the linear system \eqref{Eq1} by a gradient iterative method, we also need to control its residual for stability purpose. This thus defines a preconditioned gradient method under the form
    \begin{align}\label{EqP1Gen}
        \mathbf{M}_{\gamma}\bm{x}^{[k+1]}_{\gamma} = \mathbf{N}_{k}\bm{x}^{[k]}_{\gamma} + \bm{c}_{(\gamma,k)},  \; \forall \; \gamma > 0, \; \forall \; k = 0, 1, 2, \cdots,
    \end{align}
    where the matrix $\mathbf{M}_{\gamma} = \left(\mathbf{I} + \gamma\mathbf{A}^{t}\mathbf{A}\right)$ is symmetric and still nonsingular for any $\gamma > 0$ whatever the properties of the matrix $\mathbf{A}$ (even if $\mathbf{A}$ is not symmetric), the matrix $\mathbf{N}_{k} = \left(\mathbf{I} - \alpha_{k}\mathbf{A}\right)$, and the vector $\bm{c}_{(\gamma,k)} = \alpha_{k}\bm{b} + \gamma\mathbf{A}^{t}\bm{b}$.
    \item The regularisation parameter $\gamma$ in \eqref{EqP1} controls the amount of residual we would like to get rid of. The larger this parameter is, the less is the residual and the fewer is the performed iterations. 
\end{itemize}
\end{remark}

\subsection{Singular Value Expansion of the Iterate $\bm{x}^{[k]}_{\gamma}$}
Similarly to the Truncated SVD and the Tikhonov solutions (\cite{HansenBook10}, Chap.$4$), the iterate $\bm{x}^{[k]}_{\gamma}$ from the stabilized gradient method \eqref{EqP1Gen} also admits a simple filtered SVD expression. Taking the initial vector $\bm{x}^{[0]}_{\gamma} = \bm{0}$, we can express from \eqref{EqP1Gen} the $k^{\text{th}}$ iteration in the following
\begin{align}\label{Iter3}
    \bm{x}^{[k]}_{\gamma} = \sum_{\ell=0}^{k-1}\left(\mathbf{M}_{\gamma}^{-1}\mathbf{N}_{k}\right)^{\ell}\mathbf{M}_{\gamma}^{-1}\bm{c}_{(\gamma,k)},  \; \forall \; \gamma > 0, \; \forall \; k = 1, 2, 3, \cdots.
\end{align}
In order to have some insight into $\bm{x}^{[k]}_{\gamma}$, we use SVD decomposition of the matrix $\mathbf{A} \in \mathcal{M}_{n,n}(\mathbb{R})$ \cite{BjÄorck96}   
\begin{align}\label{Decom1}
    \mathbf{A} = \mathbf{V}\mathbf{\Sigma}\mathbf{V}^{t},  \, \text{ and } \, \mathbf{A}^{t}\mathbf{A} = \mathbf{V}\mathbf{\Sigma}^{2}\mathbf{V}^{t},
\end{align}
where $\mathbf{\Sigma} \in \mathcal{M}_{n,n}(\mathbb{R})$ is a diagonal matrix with the singular values $\sigma_{i}$, $i = 1, \cdots, n$, satisfying 
\begin{align*}
     \sigma_{1} \geq \sigma_{2} \geq \cdots \geq \sigma_{n} > 0, \; \text{ and }\, \mathbf{\Sigma} = \text{diag}(\sigma_{1}, \sigma_{2}, \cdots, \sigma_{n}).
\end{align*}
The matrix $\mathbf{V} \in \mathcal{M}_{n,n}(\mathbb{R})$ is orthogonal, that is $\mathbf{V}\mathbf{V}^{t} = \mathbf{V}^{t}\mathbf{V} = \mathbf{I}$ (see \cite{HansenBook10} - Chap 3.2). We thus use the decomposition in \eqref{Decom1} to rewrite the matrix $\mathbf{M}_{\gamma}^{-1}\mathbf{N}_{k}$ and the vector $\mathbf{M}_{\gamma}^{-1}\bm{c}_{(\gamma,k)}$ by
\begin{align}
    \mathbf{M}_{\gamma}^{-1}\mathbf{N}_{k} & = \mathbf{V}\left(\mathbf{I} + \gamma\mathbf{\Sigma}^{2}\right)^{-1}\left(\mathbf{I} - \alpha_{k}\mathbf{\Sigma}\right)\mathbf{V}^{t},  \label{InterF1}\\ 
    \mathbf{M}_{\gamma}^{-1}\bm{c}_{(\gamma,k)} & = \mathbf{V}\left(\mathbf{I} + \gamma\mathbf{\Sigma}^{2}\right)^{-1}\left(\alpha_{k}\mathbf{I} + \gamma\mathbf{\Sigma}\right)\mathbf{V}^{t}\bm{b}.\label{InterF12}
\end{align}
The $k^{\text{th}}$ iteration is thus express from \eqref{Iter3}, using \eqref{InterF1}-\eqref{InterF12} and fixing the stepsize $\alpha_{k} = \alpha$, as follows
\begin{align}\label{BGay5} 
    \bm{x}^{[k]}_{\gamma} & = \mathbf{V}\varphi^{[k]}_{\gamma}\left(\mathbf{I} + \gamma\mathbf{\Sigma}^{2}\right)^{-1}\left(\alpha\mathbf{I} + \gamma\mathbf{\Sigma}\right)\mathbf{V}^{t}\bm{b},
\end{align}
in which $\varphi^{[k]}_{\gamma} = \text{diag}\left(\varphi_{(\gamma,1)}^{[k]}, \cdots, \varphi_{(\gamma,n)}^{[k]}\right)$ is a diagonal matrix with entries defined by
\begin{align*}
    \varphi_{(\gamma,i)}^{[k]} & = \sum_{\ell=0}^{k-1}\left(\frac{1 - \alpha\sigma_{i}}{1 + \gamma\sigma^{2}_{i}}\right)^{\ell} = \left(\frac{1 - \left(\frac{1 - \alpha\sigma_{i}}{1 + \gamma\sigma^{2}_{i}}\right)^{k}}{1 - \left(\frac{1 - \alpha\sigma_{i}}{1 + \gamma\sigma^{2}_{i}}\right)}\right),  \; i = 1, \cdots, n.
\end{align*}
Multiplying diagonal matrices $\varphi^{[k]}_{\gamma}$ and $\left(\mathbf{I} + \gamma\mathbf{\Sigma}^{2}\right)^{-1}\left(\alpha\mathbf{I} + \gamma\mathbf{\Sigma}\right)$, we obtain the filtered expansion 
\begin{align}\label{SVDSol}
    \bm{x}^{[k]}_{\gamma} & = \sum^{n}_{i = 1}\left(\phi_{(\gamma,i)}^{[k]}\frac{\bm{v}^{t}_{i}\bm{b}}{\sigma_{i}}\right)\bm{v}_{i}, 
\end{align}
where $\bm{v}_{i}$ , $i = 1, \cdots, n$, are columns of $\mathbf{V}$, and $\phi_{(\gamma,i)}^{[k]}$, $i = 1, \cdots, n$, are the filter factors defined by
\begin{align}\label{FilFac} 
     \phi_{(\gamma,i)}^{[k]} = \sigma_{i}\left(\frac{\alpha + \gamma\sigma_{i}}{1 + \gamma\sigma^{2}_{i}}\right)\varphi_{(\gamma,i)}^{[k]} 
     = \left(\frac{\alpha\sigma_{i} + \gamma\sigma^{2}_{i}}{1 + \gamma\sigma^{2}_{i}}\right)\left(\frac{1 - \left(\frac{1 - \alpha\sigma_{i}}{1 + \gamma\sigma^{2}_{i}}\right)^{k}}{1 - \left(\frac{1 - \alpha\sigma_{i}}{1 + \gamma\sigma^{2}_{i}}\right)}\right), \, \forall \; k = 1, 2, 3, \cdots.
\end{align}

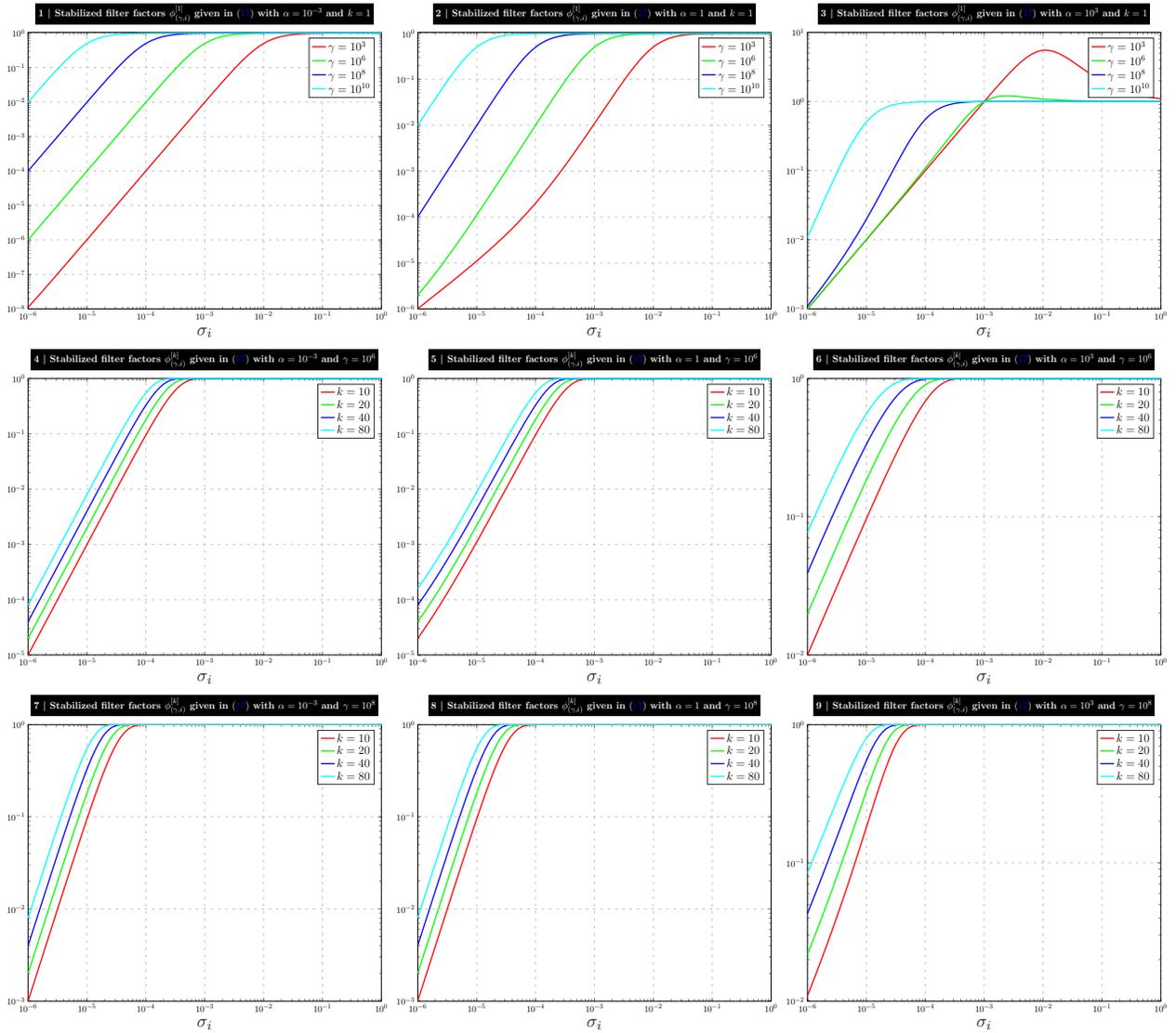
\begin{figure}[h]  
    \centering
    \scalebox{0.33}{
%
%
\definecolor{mycolor1}{rgb}{0.00000,1.00000,1.00000}%
\begin{tikzpicture}

\begin{axis}[%
width=6.028in,
height=4.719in,
at={(1.011in,0.676in)},
scale only axis,
xmode=log,
xmin=1e-06,
xmax=1,
xlabel style={font=\color{white!15!black}, scale=2},
xlabel={$\sigma_{i}$},
ymode=log,
ymin=1e-08,
ymax=1.05,
axis background/.style={fill=white},
title style={font=\bfseries, scale=1.1},
title={\normalsize\colorbox{black}{\color{white}1 | Stabilized filter factors $\phi_{(\gamma,i)}^{[1]}$ given in \eqref{FilFac} with $\alpha = 10^{-3}$ and $k = 1$}},
xmajorgrids,
ymajorgrids,
legend style={legend cell align=left, align=left, draw=white!15!black, font=\Large},
grid style = loosely dashed,
]
\addplot [color=red, line width=1.5pt]
  table[row sep=crcr]{%
1e-06	1.099999989e-08\\
1.51991108295293e-06	2.46212075150027e-08\\
2.65608778294668e-06	7.3204105725728e-08\\
5.33669923120631e-06	2.90140203441816e-07\\
1.4174741629268e-05	2.02340367870486e-06\\
8.69749002617783e-05	7.57275791420741e-05\\
0.000811130830789687	0.00653713340479613\\
0.0014174741629268	0.0196979693978369\\
0.00187381742286038	0.033922700117042\\
0.00247707635599171	0.0578141284950967\\
0.00327454916287773	0.0968455643276337\\
0.00376493580679246	0.124152837313844\\
0.00432876128108305	0.157814512699639\\
0.0049770235643321	0.198534180247906\\
0.00572236765935022	0.246683058413145\\
0.00657933224657568	0.302107557624013\\
0.00756463327554629	0.363968292159851\\
0.00869749002617783	0.430679088585257\\
0.01	0.500005\\
0.0114975699539774	0.569330814826775\\
0.0132194114846603	0.636041330606184\\
0.0151991108295293	0.697901625768909\\
0.0174752840000768	0.753325563149199\\
0.0200923300256505	0.801473797620341\\
0.0231012970008316	0.842192778571781\\
0.0265608778294668	0.875853757727674\\
0.0305385550883341	0.903160350538937\\
0.0351119173421513	0.92497525345628\\
0.0403701725859656	0.942190539249156\\
0.0464158883361278	0.955645040926091\\
0.0613590727341318	0.974127936635586\\
0.0811130830789687	0.985029657905855\\
0.123284673944207	0.993464478253231\\
0.215443469003188	0.997850660081333\\
0.756463327554629	0.999825409864237\\
1	0.999900109989001\\
};
\addlegendentry{$\gamma = 10^{3}$}

\addplot [color=green, line width=1.5pt]
  table[row sep=crcr]{%
1e-06	1.000998999001e-06\\
8.69749002617783e-05	0.00750792554702321\\
0.000151991108295293	0.0225798257314897\\
0.000200923300256505	0.0388038551787009\\
0.000265608778294668	0.0658992284262069\\
0.000351119173421513	0.109754034683379\\
0.000403701725859655	0.140136697235704\\
0.000464158883361277	0.177255412247813\\
0.00053366992312063	0.221671330489659\\
0.000613590727341318	0.273516854386047\\
0.000705480231071865	0.332311062859463\\
0.000811130830789687	0.396839888364195\\
0.00093260334688322	0.465169353749309\\
0.00107226722201032	0.534831643821322\\
0.00123284673944207	0.603161090120123\\
0.0014174741629268	0.667689879223893\\
0.00162975083462064	0.726484037141143\\
0.00187381742286038	0.778329500251986\\
0.00215443469003188	0.82274535152096\\
0.00247707635599171	0.859863997021175\\
0.0028480358684358	0.890246590481696\\
0.00327454916287773	0.914695497180588\\
0.00376493580679246	0.934101267784655\\
0.00432876128108305	0.949336972120269\\
0.00572236765935022	0.970366581449371\\
0.00756463327554629	0.982824986433812\\
0.0114975699539774	0.992492247096791\\
0.0200923300256505	0.997529094036983\\
0.0613590727341318	0.999734478044175\\
1	0.999999001000999\\
};
\addlegendentry{$\gamma = 10^{6}$}

\addplot [color=blue, line width=1.5pt]
  table[row sep=crcr]{%
1e-06	9.99910008999101e-05\\
8.11130830789687e-06	0.0065363356340723\\
1.4174741629268e-05	0.0196965937386147\\
1.87381742286038e-05	0.0339209079637322\\
2.47707635599171e-05	0.0578118179616917\\
3.27454916287773e-05	0.0968426364686812\\
3.76493580679247e-05	0.124149572768293\\
4.32876128108306e-05	0.157810903520286\\
4.97702356433211e-05	0.198530231203124\\
5.72236765935022e-05	0.246678790739784\\
6.57933224657568e-05	0.302103011844516\\
7.56463327554629e-05	0.363963528890663\\
8.69749002617783e-05	0.430674186396301\\
0.0001	0.50000005\\
0.000114975699539774	0.56932591263782\\
0.000132194114846603	0.636036567336997\\
0.000151991108295293	0.697897079989412\\
0.000174752840000768	0.753321295475839\\
0.000200923300256505	0.801469848575559\\
0.000231012970008316	0.842189169392428\\
0.000265608778294668	0.875850493182123\\
0.000305385550883341	0.903157422679984\\
0.000351119173421513	0.924972645444635\\
0.000403701725859655	0.942188228715752\\
0.000464158883361277	0.955643002644402\\
0.000613590727341318	0.974126364928187\\
0.000811130830789687	0.985028455660641\\
0.00123284673944207	0.993463680482507\\
0.00215443469003188	0.99785020155191\\
0.00756463327554629	0.999825279014929\\
1	0.999999990010001\\
};
\addlegendentry{$\gamma = 10^{8}$}

\addplot [color=mycolor1, line width=1.5pt]
  table[row sep=crcr]{%
1e-06	0.00990099108910892\\
1.51991108295293e-06	0.0225796786578835\\
2.00923300256505e-06	0.0388036639832277\\
2.65608778294668e-06	0.0658989828018304\\
3.51119173421513e-06	0.109753725226666\\
4.03701725859656e-06	0.140136353578549\\
4.64158883361278e-06	0.17725503418227\\
5.33669923120631e-06	0.221670919272545\\
6.13590727341318e-06	0.273516413080088\\
7.05480231071865e-06	0.332310596528201\\
8.11130830789687e-06	0.396839404014457\\
9.3260334688322e-06	0.465168859951847\\
1.07226722201032e-05	0.53483115002386\\
1.23284673944207e-05	0.603160605770385\\
1.4174741629268e-05	0.667689412892632\\
1.62975083462064e-05	0.726483595835184\\
1.87381742286038e-05	0.778329089034872\\
2.15443469003188e-05	0.822744973455417\\
2.47707635599171e-05	0.859863653364021\\
2.8480358684358e-05	0.890246281024984\\
3.27454916287773e-05	0.914695220639233\\
3.76493580679247e-05	0.934101022160278\\
4.32876128108306e-05	0.949336755004259\\
5.72236765935022e-05	0.970366413570827\\
7.56463327554629e-05	0.982824857809388\\
0.000114975699539774	0.992492161638102\\
0.000200923300256505	0.997529044886201\\
0.000613590727341318	0.999734461913926\\
0.572236765935022	0.999999999694788\\
1	0.999999999900101\\
};
\addlegendentry{$\gamma = 10^{10}$}

\end{axis}
\end{tikzpicture}%
%
%
\definecolor{mycolor1}{rgb}{0.00000,1.00000,1.00000}%
\begin{tikzpicture}

\begin{axis}[%
width=6.028in,
height=4.719in,
at={(1.011in,0.676in)},
scale only axis,
xmode=log,
xmin=1e-06,
xmax=1,
xlabel style={font=\color{white!15!black}, scale=2},
xlabel={$\sigma_{i}$},
ymode=log,
ymin=1e-06,
ymax=1.05,
axis background/.style={fill=white},
title style={font=\bfseries, scale=1.1},
title={\normalsize\colorbox{black}{\color{white}2 | Stabilized filter factors $\phi_{(\gamma,i)}^{[1]}$ given in \eqref{FilFac} with $\alpha = 1$ and $k = 1$}},
xmajorgrids,
ymajorgrids,
legend style={legend cell align=left, align=left, draw=white!15!black, font=\Large},
grid style = loosely dashed,
]
\addplot [color=red, line width=1.5pt]
  table[row sep=crcr]{%
1e-06	1.0099999899e-06\\
2.31012970008316e-06	2.36349656626252e-06\\
4.03701725859656e-06	4.19999165756463e-06\\
6.13590727341318e-06	6.51239840221622e-06\\
9.3260334688322e-06	1.01957736036861e-05\\
1.23284673944207e-05	1.38483574291016e-05\\
1.62975083462064e-05	1.89535457868717e-05\\
2.15443469003188e-05	2.61858141901488e-05\\
2.8480358684358e-05	3.65913701883699e-05\\
3.76493580679247e-05	5.1823365114382e-05\\
4.97702356433211e-05	7.4539152811508e-05\\
6.57933224657568e-05	0.000109076213627685\\
8.69749002617783e-05	0.000162608932247843\\
0.000114975699539774	0.000247137144310338\\
0.000151991108295293	0.000382915619829009\\
0.000200923300256505	0.000604381036448669\\
0.000265608778294668	0.000970404408240374\\
0.000351119173421513	0.00158201553017545\\
0.000464158883361277	0.0026129641128646\\
0.000613590727341318	0.00436210349449443\\
0.00093260334688322	0.00954705793191885\\
0.00247707635599171	0.0601456667494039\\
0.00327454916287773	0.0998000401798892\\
0.00376493580679246	0.127447060552114\\
0.00432876128108305	0.161456502773931\\
0.0049770235643321	0.202519125437715\\
0.00572236765935022	0.250989528803625\\
0.00657933224657568	0.306694662388614\\
0.00756463327554629	0.368774863794612\\
0.00869749002617783	0.435625842895379\\
0.01	0.504999999999999\\
0.0114975699539774	0.574277569136897\\
0.0132194114846603	0.640847902240945\\
0.0151991108295293	0.70248873053351\\
0.0174752840000768	0.75763203353968\\
0.0200923300256505	0.805458742810148\\
0.0231012970008316	0.845834768646073\\
0.0265608778294668	0.879147980965944\\
0.0305385550883341	0.906114826391192\\
0.0351119173421513	0.927606974297603\\
0.0403701725859656	0.944522077503463\\
0.0464158883361278	0.957701852447525\\
0.0613590727341318	0.97571393228275\\
0.0811130830789687	0.986242832621598\\
0.123284673944207	0.99426950143819\\
0.247707635599171	0.998775945794602\\
1	1\\
};
\addlegendentry{$\gamma = 10^{3}$}

\addplot [color=green, line width=1.5pt]
  table[row sep=crcr]{%
1e-06	1.999998000002e-06\\
1.32194114846603e-06	3.06946418449788e-06\\
1.74752840000768e-06	4.80136924615317e-06\\
2.31012970008316e-06	7.64678812268117e-06\\
3.05385550883341e-06	1.23797735234834e-05\\
4.03701725859656e-06	2.03341942081032e-05\\
5.33669923120631e-06	3.38160948210546e-05\\
7.05480231071865e-06	5.68222098992633e-05\\
1.07226722201032e-05	0.00012568392116312\\
1.62975083462064e-05	0.000281831429739137\\
2.8480358684358e-05	0.000838930706912773\\
6.57933224657568e-05	0.00437561361674367\\
0.000132194114846603	0.0173050671518052\\
0.000200923300256505	0.0389967887923754\\
0.000265608778294668	0.066147085751417\\
0.000351119173421513	0.110066304638079\\
0.000403701725859655	0.14048347854673\\
0.000464158883361277	0.177636914749824\\
0.00053366992312063	0.222086285941727\\
0.000613590727341318	0.273962172217687\\
0.000705480231071865	0.332781633495754\\
0.000811130830789687	0.397328641281326\\
0.00093260334688322	0.465667640279871\\
0.00107226722201032	0.535329930351884\\
0.00123284673944207	0.603649843037254\\
0.0014174741629268	0.668160449860184\\
0.00162975083462064	0.726929354972783\\
0.00187381742286038	0.778744455704054\\
0.00215443469003188	0.823126854022971\\
0.00247707635599171	0.860210778332201\\
0.0028480358684358	0.890558860436395\\
0.00327454916287773	0.914974552548877\\
0.00376493580679246	0.934349125109864\\
0.00432876128108305	0.949556061913972\\
0.00572236765935022	0.9705359861614\\
0.00756463327554629	0.982954780172107\\
0.0114975699539774	0.992578482681576\\
0.0200923300256505	0.997578691645327\\
0.0613590727341318	0.999750754931729\\
1	1\\
};
\addlegendentry{$\gamma = 10^{6}$}

\addplot [color=blue, line width=1.5pt]
  table[row sep=crcr]{%
1e-06	0.000100989901009899\\
1.4174741629268e-05	0.019710475390766\\
2.15443469003188e-05	0.0443776066482197\\
2.8480358684358e-05	0.0750537244508827\\
3.27454916287773e-05	0.0968721812272037\\
3.76493580679247e-05	0.124182515000676\\
4.32876128108306e-05	0.157847323421029\\
4.97702356433211e-05	0.198570080655023\\
5.72236765935022e-05	0.246721855443689\\
6.57933224657568e-05	0.302148882892163\\
7.56463327554629e-05	0.364011594607011\\
8.69749002617783e-05	0.430723653939402\\
0.0001	0.50005\\
0.000114975699539774	0.56937538018092\\
0.000132194114846603	0.636084633053344\\
0.000151991108295293	0.697942951037059\\
0.000174752840000768	0.753364360179744\\
0.000200923300256505	0.801509698027457\\
0.000231012970008316	0.842225589293172\\
0.000265608778294668	0.875883435414505\\
0.000305385550883341	0.903186967438507\\
0.000351119173421513	0.924998962653048\\
0.000403701725859655	0.942211544098294\\
0.000464158883361277	0.955663570759617\\
0.000613590727341318	0.974142224884661\\
0.000811130830789687	0.985040587407798\\
0.00123284673944207	0.993471730714357\\
0.00215443469003188	0.997854828530631\\
0.00756463327554629	0.999826599403394\\
1	1\\
};
\addlegendentry{$\gamma = 10^{8}$}

\addplot [color=mycolor1, line width=1.5pt]
  table[row sep=crcr]{%
1e-06	0.00990198019801979\\
1.51991108295293e-06	0.0225811627642729\\
2.00923300256505e-06	0.0388055933193645\\
2.65608778294668e-06	0.0659014613750826\\
3.51119173421513e-06	0.109756847926213\\
4.03701725859656e-06	0.140139821391659\\
4.64158883361278e-06	0.17725884920729\\
5.33669923120631e-06	0.221675068827066\\
6.13590727341318e-06	0.273520866258404\\
7.05480231071865e-06	0.332315302234565\\
8.11130830789687e-06	0.396844291543628\\
9.3260334688322e-06	0.465173842817153\\
1.07226722201032e-05	0.534836132889165\\
1.23284673944207e-05	0.603165493299557\\
1.4174741629268e-05	0.667694118598995\\
1.62975083462064e-05	0.7264880490135\\
1.87381742286038e-05	0.778333238589392\\
2.15443469003188e-05	0.822748788480438\\
2.47707635599171e-05	0.859867121177132\\
2.8480358684358e-05	0.890249403724531\\
3.27454916287773e-05	0.914698011192916\\
3.76493580679247e-05	0.93410350073353\\
4.32876128108306e-05	0.949338945902195\\
5.72236765935022e-05	0.970368107617946\\
7.56463327554629e-05	0.982826155746771\\
0.000114975699539774	0.99249302399395\\
0.000200923300256505	0.997529540862283\\
0.000613590727341318	0.999734624682802\\
0.572236765935022	0.999999999869367\\
1	1\\
};
\addlegendentry{$\gamma = 10^{10}$}

\end{axis}
\end{tikzpicture}%
%
%
\definecolor{mycolor1}{rgb}{0.00000,1.00000,1.00000}%
\begin{tikzpicture}

\begin{axis}[%
width=6.028in,
height=4.719in,
at={(1.011in,0.676in)},
scale only axis,
xmode=log,
xmin=1e-06,
xmax=1,
xlabel style={font=\color{white!15!black}, scale=2},
xlabel={$\sigma_{i}$},
ymode=log,
ymin=0.001,
ymax=10,
axis background/.style={fill=white},
title style={font=\bfseries, scale=1.1},
title={\normalsize\colorbox{black}{\color{white}3 | Stabilized filter factors $\phi_{(\gamma,i)}^{[1]}$ given in \eqref{FilFac} with $\alpha = 10^{3}$ and $k = 1$}},
xmajorgrids,
ymajorgrids,
legend style={legend cell align=left, align=left, draw=white!15!black, font=\Large},
grid style = loosely dashed,
]
\addplot [color=red, line width=1.5pt]
  table[row sep=crcr]{%
1e-06	0.0010000099899999\\
0.00123284673944207	1.22936066133057\\
0.00162975083462064	1.61345688134168\\
0.00215443469003188	2.10322740977061\\
0.00247707635599171	2.39168392105667\\
0.0028480358684358	2.70938257741981\\
0.00327454916287773	3.05427589243529\\
0.00376493580679246	3.42167029882176\\
0.00432876128108305	3.80344657706651\\
0.0049770235643321	4.18746431524584\\
0.00572236765935022	4.55745991928397\\
0.00657933224657568	4.89379942698993\\
0.00756463327554629	5.17534649855504\\
0.00869749002617783	5.38238015301713\\
0.01	5.5\\
0.0114975699539774	5.52103187925865\\
0.0132194114846603	5.44741953700137\\
0.0151991108295293	5.28959349513484\\
0.0174752840000768	5.06410242402003\\
0.0200923300256505	4.79040393261827\\
0.0231012970008316	4.48782484293865\\
0.0265608778294668	4.17337121923559\\
0.0305385550883341	3.8605906786466\\
0.0403701725859656	3.27606033181073\\
0.0533669923120631	2.77633497857842\\
0.0705480231071865	2.36985819005817\\
0.0811130830789687	2.19941754836407\\
0.093260334688322	2.04871201233297\\
0.107226722201032	1.91593948234477\\
0.123284673944207	1.79929268639705\\
0.14174741629268	1.69703405260259\\
0.162975083462064	1.60753844827661\\
0.187381742286038	1.52931438090968\\
0.215443469003188	1.46101123008466\\
0.247707635599171	1.40141776408898\\
0.28480358684358	1.34945550160634\\
0.327454916287773	1.30416927824954\\
0.376493580679247	1.26471654577372\\
0.432876128108306	1.23035636582116\\
0.497702356433211	1.20043868108916\\
0.657933224657568	1.15172504487204\\
0.869749002617783	1.11482832579604\\
1	1.0998900109989\\
};
\addlegendentry{$\gamma = 10^{3}$}

\addplot [color=green, line width=1.5pt]
  table[row sep=crcr]{%
1e-06	0.001000998999001\\
6.13590727341318e-06	0.00617332420978745\\
1.62975083462064e-05	0.0165587189833818\\
3.27454916287773e-05	0.0337815360170074\\
5.72236765935022e-05	0.0603007679272428\\
8.69749002617783e-05	0.0938297459191088\\
0.000132194114846603	0.147098805445448\\
0.000305385550883341	0.364639485146338\\
0.000351119173421513	0.422336259338373\\
0.000403701725859655	0.487264789572944\\
0.000464158883361277	0.559139416761047\\
0.00053366992312063	0.637041738009918\\
0.000613590727341318	0.719280003857335\\
0.000705480231071865	0.803352269786411\\
0.000811130830789687	0.886081558411785\\
0.00093260334688322	0.963954170841958\\
0.00107226722201032	1.03361646091397\\
0.00123284673944207	1.09240276016771\\
0.0014174741629268	1.13873108615084\\
0.00162975083462064	1.17224718661243\\
0.00187381742286038	1.19369990777225\\
0.00215443469003188	1.2046293560342\\
0.00247707635599171	1.20699208935842\\
0.0028480358684358	1.20282881513669\\
0.00376493580679246	1.18220645032003\\
0.00657933224657568	1.12597952096464\\
0.01	1.08910891089109\\
0.0151991108295293	1.06119964254162\\
0.0231012970008316	1.04133633863644\\
0.0351119173421513	1.02764680267954\\
0.0613590727341318	1.01602764248537\\
0.123284673944207	1.00804498567909\\
0.327454916287773	1.00304450108224\\
1	1.000998999001\\
};
\addlegendentry{$\gamma = 10^{6}$}

\addplot [color=blue, line width=1.5pt]
  table[row sep=crcr]{%
1e-06	0.0010998900109989\\
1.32194114846603e-06	0.00149643248264058\\
1.74752840000768e-06	0.00205228721203021\\
2.31012970008316e-06	0.00284228278236983\\
3.05385550883341e-06	0.00398274453483366\\
4.03701725859656e-06	0.00565754770013101\\
5.33669923120631e-06	0.00816149088087347\\
7.05480231071865e-06	0.011972239755669\\
9.3260334688322e-06	0.0178681157366043\\
1.23284673944207e-05	0.0271154475317231\\
1.62975083462064e-05	0.0417494832515845\\
2.47707635599171e-05	0.0811505158873074\\
3.27454916287773e-05	0.126416939749758\\
3.76493580679247e-05	0.157124747383372\\
4.32876128108306e-05	0.194267224163955\\
4.97702356433211e-05	0.238419532553104\\
5.72236765935022e-05	0.289786559348493\\
6.57933224657568e-05	0.348019930538176\\
7.56463327554629e-05	0.412077310954615\\
8.69749002617783e-05	0.48019119704062\\
0.0001	0.55\\
0.000114975699539774	0.618842923282138\\
0.000132194114846603	0.684150349400949\\
0.000151991108295293	0.743813998683071\\
0.000174752840000768	0.796429064084547\\
0.000200923300256505	0.841359149925538\\
0.000231012970008316	0.878645490036097\\
0.000265608778294668	0.908825667797202\\
0.000305385550883341	0.932731725961061\\
0.000351119173421513	0.951316171066291\\
0.000403701725859655	0.965526926641367\\
0.000464158883361277	0.976231685973962\\
0.000613590727341318	0.990002181356295\\
0.000811130830789687	0.997172334565224\\
0.00123284673944207	1.00152196256395\\
0.0028480358684358	1.00227553960319\\
0.215443469003188	1.00004620043491\\
1	1.0000099899999\\
};
\addlegendentry{$\gamma = 10^{8}$}

\addplot [color=mycolor1, line width=1.5pt]
  table[row sep=crcr]{%
1e-06	0.0108910891089109\\
2.31012970008316e-06	0.0528563382168821\\
3.05385550883341e-06	0.088098129184037\\
3.51119173421513e-06	0.112879547473217\\
4.03701725859656e-06	0.143607634501921\\
4.64158883361278e-06	0.181073874227403\\
5.33669923120631e-06	0.225824623347748\\
6.13590727341318e-06	0.277974044574801\\
7.05480231071865e-06	0.337021008597471\\
8.11130830789687e-06	0.401731820714933\\
9.3260334688322e-06	0.470156708122773\\
1.07226722201032e-05	0.539818998194786\\
1.23284673944207e-05	0.608053022470862\\
1.4174741629268e-05	0.672399824961902\\
1.62975083462064e-05	0.730941227329897\\
1.87381742286038e-05	0.782482793110073\\
2.15443469003188e-05	0.82656381350055\\
2.47707635599171e-05	0.863334934287392\\
2.8480358684358e-05	0.893372103271535\\
3.27454916287773e-05	0.917488564875808\\
3.76493580679247e-05	0.936582073985632\\
4.32876128108306e-05	0.951529843839217\\
4.97702356433211e-05	0.96312760535219\\
6.57933224657568e-05	0.978905914809077\\
8.69749002617783e-05	0.988087817552196\\
0.000132194114846603	0.9950623507128\\
0.000265608778294668	0.998960492892219\\
0.0014174741629268	1.0000207767534\\
1	1.0000000999\\
};
\addlegendentry{$\gamma = 10^{10}$}

\end{axis}
\end{tikzpicture}%
}     \\
        \scalebox{0.33}{
%
%
\definecolor{mycolor1}{rgb}{0.00000,1.00000,1.00000}%
\begin{tikzpicture}

\begin{axis}[%
width=6.028in,
height=4.719in,
at={(1.011in,0.676in)},
scale only axis,
xmode=log,
xmin=1e-06,
xmax=1,
xlabel style={font=\color{white!15!black}, scale=2},
xlabel={$\sigma_{i}$},
ymode=log,
ymin=1e-05,
ymax=1,
axis background/.style={fill=white},
title style={font=\bfseries, scale=1.1},
title={\normalsize\colorbox{black}{\color{white}4 | Stabilized filter factors $\phi_{(\gamma,i)}^{[k]}$ given in \eqref{FilFac} with $\alpha = 10^{-3}$ and $\gamma = 10^{6}$}},
xmajorgrids,
ymajorgrids,
legend style={legend cell align=left, align=left, draw=white!15!black, font=\Large},
grid style = loosely dashed,
]
\addplot [color=red, line width=1.5pt]
  table[row sep=crcr]{%
1e-06	1.00099449001892e-05\\
3.76493580679247e-05	0.0140652287051556\\
6.57933224657568e-05	0.0422752401857486\\
8.69749002617783e-05	0.0725927772694276\\
0.000114975699539774	0.123070892655362\\
0.000132194114846603	0.159068264763226\\
0.000151991108295293	0.204183411461427\\
0.000174752840000768	0.259786849940163\\
0.000200923300256505	0.326837053651596\\
0.000231012970008316	0.405431437561687\\
0.000265608778294668	0.494248350253688\\
0.000305385550883341	0.590019916917618\\
0.000351119173421513	0.687319972954973\\
0.000403701725859655	0.779049929079573\\
0.000464158883361277	0.857881589800044\\
0.00053366992312063	0.91841125017348\\
0.000613590727341318	0.959049856184499\\
0.000705480231071865	0.982390711643763\\
0.000811130830789687	0.993627261041316\\
0.00093260334688322	0.99808502615159\\
0.0014174741629268	0.99998357763567\\
1	1\\
};
\addlegendentry{$k = 10$}

\addplot [color=green, line width=1.5pt]
  table[row sep=crcr]{%
1e-06	2.00197896013072e-05\\
2.8480358684358e-05	0.0160858291167905\\
4.97702356433211e-05	0.0482770111114709\\
6.57933224657568e-05	0.0827632844387345\\
8.69749002617783e-05	0.139915843227166\\
0.0001	0.180457168750088\\
0.000114975699539774	0.230995340691737\\
0.000132194114846603	0.292833816671669\\
0.000151991108295293	0.366675957406828\\
0.000174752840000768	0.452084492478494\\
0.000200923300256505	0.546851647663536\\
0.000231012970008316	0.646488224560038\\
0.000265608778294668	0.744215268778883\\
0.000305385550883341	0.831916331475763\\
0.000351119173421513	0.902231200687122\\
0.000403701725859655	0.951181066160258\\
0.000464158883361277	0.979802357482236\\
0.00053366992312063	0.993343275901745\\
0.000613590727341318	0.99832308572149\\
0.00093260334688322	0.99999633287516\\
1	1\\
};
\addlegendentry{$k = 20$}

\addplot [color=blue, line width=1.5pt]
  table[row sep=crcr]{%
1e-06	4.00391784107702e-05\\
1.87381742286038e-05	0.0139449076512482\\
3.27454916287773e-05	0.0419631344674624\\
4.32876128108306e-05	0.0721491602896376\\
5.72236765935022e-05	0.122580678368792\\
6.57933224657568e-05	0.158676807626382\\
7.56463327554629e-05	0.204071013569331\\
8.69749002617783e-05	0.260255243268364\\
0.0001	0.328349547746878\\
0.000114975699539774	0.408631833962182\\
0.000132194114846603	0.499915989156841\\
0.000151991108295293	0.598900657073443\\
0.000174752840000768	0.699788596617451\\
0.000200923300256505	0.794656570774748\\
0.000231012970008316	0.875029424625285\\
0.000265608778294668	0.93457417127414\\
0.000305385550883341	0.971747880375434\\
0.000351119173421513	0.990441261880919\\
0.000403701725859655	0.997616711698751\\
0.00053366992312063	0.99995568802428\\
1	1\\
};
\addlegendentry{$k = 40$}

\addplot [color=mycolor1, line width=1.5pt]
  table[row sep=crcr]{%
1e-06	8.0076753685694e-05\\
1.4174741629268e-05	0.015944895912812\\
2.47707635599171e-05	0.0478895064407731\\
3.27454916287773e-05	0.0821653642805906\\
4.32876128108306e-05	0.139092819248776\\
4.97702356433211e-05	0.179568664700702\\
5.72236765935022e-05	0.230135334028231\\
6.57933224657568e-05	0.292175285974264\\
7.56463327554629e-05	0.366497048559447\\
8.69749002617783e-05	0.452777694888051\\
0.0001	0.548885669988178\\
0.000114975699539774	0.65028369219707\\
0.000132194114846603	0.749915982099018\\
0.000151991108295293	0.839119317103886\\
0.000174752840000768	0.90987311327908\\
0.000200923300256505	0.957834076074016\\
0.000231012970008316	0.984382355290512\\
0.000265608778294668	0.995719460935535\\
0.000351119173421513	0.999908630525571\\
0.0464158883361278	1\\
1	1\\
};
\addlegendentry{$k = 80$}

\end{axis}
\end{tikzpicture}%
%
%
\definecolor{mycolor1}{rgb}{0.00000,1.00000,1.00000}%
\begin{tikzpicture}

\begin{axis}[%
width=6.028in,
height=4.719in,
at={(1.011in,0.676in)},
scale only axis,
xmode=log,
xmin=1e-06,
xmax=1,
xlabel style={font=\color{white!15!black}, scale=2},
xlabel={$\sigma_{i}$},
ymode=log,
ymin=1e-05,
ymax=1,
axis background/.style={fill=white},
title style={font=\bfseries, scale=1.1},
title={\normalsize\colorbox{black}{\color{white}5 | Stabilized filter factors $\phi_{(\gamma,i)}^{[k]}$ given in \eqref{FilFac} with $\alpha = 1$ and $\gamma = 10^{6}$}},
xmajorgrids,
ymajorgrids,
legend style={legend cell align=left, align=left, draw=white!15!black, font=\Large},
grid style = loosely dashed,
]
\addplot [color=red, line width=1.5pt]
  table[row sep=crcr]{%
1e-06	1.99998000013374e-05\\
1.32194114846603e-06	3.06942178760077e-05\\
1.74752840000768e-06	4.80126550832373e-05\\
2.31012970008316e-06	7.64652499789255e-05\\
3.05385550883341e-06	0.000123790838816841\\
4.03701725859656e-06	0.0002033233365145\\
5.33669923120631e-06	0.000338109494078479\\
7.05480231071865e-06	0.000568076826647146\\
1.07226722201032e-05	0.00125612860966016\\
1.62975083462064e-05	0.00281474267937058\\
3.27454916287773e-05	0.0109836224789717\\
6.57933224657568e-05	0.0429045428940378\\
8.69749002617783e-05	0.0733982672417109\\
0.000114975699539774	0.124077619431644\\
0.000132194114846603	0.160178155767319\\
0.000151991108295293	0.205390947217496\\
0.000174752840000768	0.261078085403024\\
0.000200923300256505	0.328187022780221\\
0.000231012970008316	0.406802170708156\\
0.000265608778294668	0.495588726819402\\
0.000305385550883341	0.591268969479884\\
0.000351119173421513	0.688415025382268\\
0.000403701725859655	0.779939401264527\\
0.000464158883361277	0.858539212296726\\
0.00053366992312063	0.918845187526285\\
0.000613590727341318	0.959300180088154\\
0.000705480231071865	0.982514424687809\\
0.000811130830789687	0.993678712747936\\
0.00093260334688322	0.998102792804942\\
0.0014174741629268	0.999983808709698\\
1	1\\
};
\addlegendentry{$k = 10$}

\addplot [color=green, line width=1.5pt]
  table[row sep=crcr]{%
1e-06	3.99992000106353e-05\\
1.32194114846603e-06	6.13874936169788e-05\\
1.74752840000768e-06	9.60230049513449e-05\\
2.31012970008316e-06	0.000152924653023359\\
3.05385550883341e-06	0.000247566353461863\\
4.03701725859656e-06	0.000406605332649769\\
5.33669923120631e-06	0.000676104670127027\\
7.05480231071865e-06	0.00113583094201326\\
1.07226722201032e-05	0.00251067936023637\\
1.62975083462064e-05	0.00562156258239008\\
6.57933224657568e-05	0.0839682859871292\\
8.69749002617783e-05	0.141409228849336\\
0.0001	0.182093062402769\\
0.000114975699539774	0.232759983219465\\
0.000132194114846603	0.294699269949619\\
0.000151991108295293	0.368596453236092\\
0.000174752840000768	0.45399440412834\\
0.000200923300256505	0.548667323639098\\
0.000231012970008316	0.648116335323444\\
0.000265608778294668	0.745569267488328\\
0.000305385550883341	0.832938944689965\\
0.000351119173421513	0.902914803592468\\
0.000403701725859655	0.951573332884185\\
0.000464158883361277	0.979988845542369\\
0.00053366992312063	0.993413896412356\\
0.000613590727341318	0.998343524659143\\
0.00093260334688322	0.99999640060486\\
1	1\\
};
\addlegendentry{$k = 20$}

\addplot [color=blue, line width=1.5pt]
  table[row sep=crcr]{%
1e-06	7.99968000853089e-05\\
1.32194114846603e-06	0.000122771218809496\\
1.74752840000768e-06	0.00019203678948526\\
2.31012970008316e-06	0.000305825920097096\\
3.05385550883341e-06	0.000495071417824484\\
4.03701725859656e-06	0.000813045337403054\\
5.33669923120631e-06	0.00135175222272922\\
7.05480231071865e-06	0.00227037177209768\\
1.07226722201032e-05	0.0050150552096229\\
1.87381742286038e-05	0.0146829740248753\\
3.76493580679247e-05	0.0565055744550416\\
4.97702356433211e-05	0.0960230324985727\\
6.57933224657568e-05	0.160885898922642\\
7.56463327554629e-05	0.206473427651151\\
8.69749002617783e-05	0.262821887694908\\
0.0001	0.331028241430319\\
0.000114975699539774	0.411342756650604\\
0.000132194114846603	0.502550880190399\\
0.000151991108295293	0.601329561133956\\
0.000174752840000768	0.701877889276833\\
0.000200923300256505	0.796298815248905\\
0.000231012970008316	0.876177886533797\\
0.000265608778294668	0.935265002353574\\
0.000305385550883341	0.972090603798696\\
0.000351119173421513	0.990574464638511\\
0.000403701725859655	0.997654857912055\\
0.00053366992312063	0.999956623239534\\
1	1\\
};
\addlegendentry{$k = 40$}

\addplot [color=mycolor1, line width=1.5pt]
  table[row sep=crcr]{%
1e-06	0.000159987200682678\\
1.32194114846603e-06	0.000245527364846851\\
1.74752840000768e-06	0.000384036700842081\\
2.31012970008316e-06	0.000611558310700821\\
3.05385550883341e-06	0.000989897739940128\\
4.03701725859656e-06	0.00162542963208555\\
5.33669923120631e-06	0.00270167721138672\\
7.05480231071865e-06	0.00453558895621167\\
1.07226722201032e-05	0.0100049596404902\\
3.76493580679247e-05	0.109818268965589\\
4.97702356433211e-05	0.182825642226923\\
5.72236765935022e-05	0.233648230298797\\
6.57933224657568e-05	0.295887525373137\\
7.56463327554629e-05	0.370315578976286\\
8.69749002617783e-05	0.456568430738302\\
0.0001	0.552476786236189\\
0.000114975699539774	0.653482649852292\\
0.000132194114846603	0.752544373200653\\
0.000151991108295293	0.841061881174356\\
0.000174752840000768	0.911123207097963\\
0.000200923300256505	0.958505827331001\\
0.000231012970008316	0.984668084216763\\
0.000265608778294668	0.995809380079717\\
0.000351119173421513	0.999911159283149\\
0.0533669923120631	1\\
1	1\\
};
\addlegendentry{$k = 80$}

\end{axis}
\end{tikzpicture}%
%
%
\definecolor{mycolor1}{rgb}{0.00000,1.00000,1.00000}%
\begin{tikzpicture}

\begin{axis}[%
width=6.028in,
height=4.719in,
at={(1.011in,0.676in)},
scale only axis,
xmode=log,
xmin=1e-06,
xmax=1,
xlabel style={font=\color{white!15!black}, scale=2},
xlabel={$\sigma_{i}$},
ymode=log,
ymin=0.00996502018460165,
ymax=1,
axis background/.style={fill=white},
title style={font=\bfseries, scale=1.1},
title={\normalsize\colorbox{black}{\color{white}6 | Stabilized filter factors $\phi_{(\gamma,i)}^{[k]}$ given in \eqref{FilFac} with $\alpha = 10^{3}$ and $\gamma = 10^{6}$}},
xmajorgrids,
ymajorgrids,
legend style={legend cell align=left, align=left, draw=white!15!black, font=\Large},
grid style = loosely dashed,
]
\addplot [color=red, line width=1.5pt]
  table[row sep=crcr]{%
1e-06	0.00996502018460165\\
2.31012970008316e-06	0.0229147643397131\\
4.03701725859656e-06	0.0398011243696143\\
6.13590727341318e-06	0.0600462242134489\\
9.3260334688322e-06	0.0902338069469232\\
1.23284673944207e-05	0.118006620695228\\
1.62975083462064e-05	0.153777938997631\\
2.15443469003188e-05	0.199435885562967\\
2.47707635599171e-05	0.22660180896834\\
2.8480358684358e-05	0.256990128954455\\
3.27454916287773e-05	0.290825152528274\\
3.76493580679247e-05	0.328290098388497\\
4.32876128108306e-05	0.369502018926683\\
4.97702356433211e-05	0.414479918631163\\
5.72236765935022e-05	0.463105791201829\\
6.57933224657568e-05	0.515079356073398\\
7.56463327554629e-05	0.569869147833421\\
8.69749002617783e-05	0.626665590889887\\
0.0001	0.684345956794772\\
0.000114975699539774	0.741466432251717\\
0.000132194114846603	0.796301783392353\\
0.000151991108295293	0.846955506645802\\
0.000174752840000768	0.891557894956266\\
0.000200923300256505	0.928549421747034\\
0.000231012970008316	0.957007460930797\\
0.000265608778294668	0.976921199555485\\
0.000305385550883341	0.98927981657192\\
0.000351119173421513	0.995862375431783\\
0.000403701725859655	0.998744168355332\\
0.00053366992312063	0.99996032043827\\
0.093260334688322	1\\
1	1\\
};
\addlegendentry{$k = 10$}

\addplot [color=green, line width=1.5pt]
  table[row sep=crcr]{%
1e-06	0.0198307387419237\\
1.74752840000768e-06	0.0344353499670161\\
2.65608778294668e-06	0.0519362426510446\\
4.03701725859656e-06	0.0780181192381429\\
5.33669923120631e-06	0.102003766150821\\
7.05480231071865e-06	0.13289177392626\\
9.3260334688322e-06	0.172325473977712\\
1.23284673944207e-05	0.222087678862549\\
1.4174741629268e-05	0.251393139428774\\
1.62975083462064e-05	0.283908223472903\\
1.87381742286038e-05	0.319780892119801\\
2.15443469003188e-05	0.359097098675649\\
2.47707635599171e-05	0.401855238108957\\
2.8480358684358e-05	0.447936331528883\\
3.27454916287773e-05	0.497071035713454\\
3.76493580679247e-05	0.548805808077065\\
4.32876128108306e-05	0.602472295862472\\
4.97702356433211e-05	0.65716623431383\\
5.72236765935022e-05	0.711744608558986\\
6.57933224657568e-05	0.76485196909381\\
7.56463327554629e-05	0.814987450014452\\
8.69749002617783e-05	0.860621418974403\\
0.0001	0.900362525008192\\
0.000114975699539774	0.933160394347343\\
0.000132194114846603	0.958507036550863\\
0.000151991108295293	0.976577383053958\\
0.000174752840000768	0.988240309853684\\
0.000200923300256505	0.994894814867317\\
0.000231012970008316	0.998151641584384\\
0.000305385550883341	0.999885077667267\\
0.00215443469003188	0.999999999999984\\
1	1\\
};
\addlegendentry{$k = 20$}

\addplot [color=blue, line width=1.5pt]
  table[row sep=crcr]{%
1e-06	0.0392682192847971\\
1.51991108295293e-06	0.0591156968033652\\
2.00923300256505e-06	0.0774481042795211\\
2.65608778294668e-06	0.101175112001381\\
3.51119173421513e-06	0.131673908190924\\
4.64158883361278e-06	0.170521047064563\\
5.33669923120631e-06	0.193602763992692\\
6.13590727341318e-06	0.219404601040973\\
7.05480231071865e-06	0.248123324275253\\
8.11130830789687e-06	0.27993011145689\\
9.3260334688322e-06	0.31495487897378\\
1.07226722201032e-05	0.35326699622917\\
1.23284673944207e-05	0.394852420622544\\
1.4174741629268e-05	0.439587768305692\\
1.62975083462064e-05	0.487212567590266\\
1.87381742286038e-05	0.537301965274666\\
2.15443469003188e-05	0.589243471074029\\
2.47707635599171e-05	0.642222843822308\\
2.8480358684358e-05	0.695225705954213\\
3.27454916287773e-05	0.747062456881663\\
3.76493580679247e-05	0.796423801175011\\
4.32876128108306e-05	0.841971724443146\\
4.97702356433211e-05	0.882465009105439\\
5.72236765935022e-05	0.916908829305187\\
6.57933224657568e-05	0.944705403560939\\
7.56463327554629e-05	0.965770356347846\\
8.69749002617783e-05	0.98057361115129\\
0.0001	0.990072373577257\\
0.000114975699539774	0.995532467116198\\
0.000132194114846603	0.998278333984207\\
0.000174752840000768	0.999861709687663\\
0.00093260334688322	1\\
1	1\\
};
\addlegendentry{$k = 40$}

\addplot [color=mycolor1, line width=1.5pt]
  table[row sep=crcr]{%
1e-06	0.0769944455237954\\
1.32194114846603e-06	0.100543922624537\\
1.74752840000768e-06	0.130788566631007\\
2.31012970008316e-06	0.169270549474871\\
2.65608778294668e-06	0.19211382071427\\
3.05385550883341e-06	0.217630639907206\\
3.51119173421513e-06	0.246009798283575\\
4.03701725859656e-06	0.277413997070418\\
4.64158883361278e-06	0.311964666637131\\
5.33669923120631e-06	0.349723497759774\\
6.13590727341318e-06	0.390670823123999\\
7.05480231071865e-06	0.434681464501105\\
8.11130830789687e-06	0.481499355613514\\
9.3260334688322e-06	0.53071318215817\\
1.07226722201032e-05	0.58173642183356\\
1.23284673944207e-05	0.633796407173608\\
1.4174741629268e-05	0.685938130567406\\
1.62975083462064e-05	0.737049049162632\\
1.87381742286038e-05	0.785910528661314\\
2.15443469003188e-05	0.831279073944687\\
2.47707635599171e-05	0.871995506517405\\
2.8480358684358e-05	0.907112629688892\\
3.27454916287773e-05	0.936022599281262\\
3.76493580679247e-05	0.958556731271971\\
4.32876128108306e-05	0.975027064124527\\
4.97702356433211e-05	0.986185525915414\\
5.72236765935022e-05	0.993095857352564\\
6.57933224657568e-05	0.996942507604639\\
8.69749002617783e-05	0.999622615416298\\
0.000200923300256505	0.999999999320724\\
1	1\\
};
\addlegendentry{$k = 80$}

\end{axis}
\end{tikzpicture}%
}     \\
        \scalebox{0.33}{
%
%
\definecolor{mycolor1}{rgb}{0.00000,1.00000,1.00000}%
\begin{tikzpicture}

\begin{axis}[%
width=6.028in,
height=4.719in,
at={(1.011in,0.676in)},
scale only axis,
xmode=log,
xmin=1e-06,
xmax=1,
xlabel style={font=\color{white!15!black}, scale=2},
xlabel={$\sigma_{i}$},
ymode=log,
ymin=0.000999460209933994,
ymax=1,
axis background/.style={fill=white},
title style={font=\bfseries, scale=1.1},
title={\normalsize\colorbox{black}{\color{white}7 | Stabilized filter factors $\phi_{(\gamma,i)}^{[k]}$ given in \eqref{FilFac} with $\alpha = 10^{-3}$ and $\gamma = 10^{8}$}},
xmajorgrids,
ymajorgrids,
legend style={legend cell align=left, align=left, draw=white!15!black, font=\Large},
grid style = loosely dashed,
]
\addplot [color=red, line width=1.5pt]
  table[row sep=crcr]{%
1e-06	0.000999460209933994\\
3.05385550883341e-06	0.00927840543821531\\
5.33669923120631e-06	0.0280393242136423\\
7.05480231071865e-06	0.0484346007280873\\
9.3260334688322e-06	0.0829551888164417\\
1.07226722201032e-05	0.108027395459552\\
1.23284673944207e-05	0.140021370005453\\
1.4174741629268e-05	0.180394023008473\\
1.62975083462064e-05	0.230598841007256\\
1.87381742286038e-05	0.291847436316353\\
2.15443469003188e-05	0.364732474984047\\
2.47707635599171e-05	0.448714556333333\\
2.8480358684358e-05	0.541553515891687\\
3.27454916287773e-05	0.638894360301752\\
3.76493580679247e-05	0.734356336168727\\
4.32876128108306e-05	0.820486900880407\\
4.97702356433211e-05	0.890636738037128\\
5.72236765935022e-05	0.941142481366464\\
6.57933224657568e-05	0.972589734583126\\
7.56463327554629e-05	0.989165217605272\\
8.69749002617783e-05	0.996422260181104\\
0.000114975699539774	0.999780473263699\\
0.000613590727341318	1\\
1	1\\
};
\addlegendentry{$k = 10$}

\addplot [color=green, line width=1.5pt]
  table[row sep=crcr]{%
1e-06	0.00199792149915673\\
2.65608778294668e-06	0.0140056780418191\\
4.03701725859656e-06	0.0320439207185979\\
5.33669923120631e-06	0.055292444724927\\
7.05480231071865e-06	0.0945232909084853\\
8.11130830789687e-06	0.122918903384675\\
9.3260334688322e-06	0.159028814281312\\
1.07226722201032e-05	0.204384872749329\\
1.23284673944207e-05	0.260436755952702\\
1.4174741629268e-05	0.328246042479764\\
1.62975083462064e-05	0.408021856540623\\
1.87381742286038e-05	0.498519946548279\\
2.15443469003188e-05	0.596435171660106\\
2.47707635599171e-05	0.696084359601248\\
2.8480358684358e-05	0.789826821208726\\
3.27454916287773e-05	0.869602716978118\\
3.76493580679247e-05	0.929433443866297\\
4.32876128108306e-05	0.967775047244478\\
4.97702356433211e-05	0.98803967693284\\
5.72236765935022e-05	0.996535792500302\\
6.57933224657568e-05	0.999248677349777\\
0.000114975699539774	0.999999951808012\\
1	1\\
};
\addlegendentry{$k = 20$}

\addplot [color=blue, line width=1.5pt]
  table[row sep=crcr]{%
1e-06	0.00399185130799662\\
2.00923300256505e-06	0.0160152605892088\\
3.05385550883341e-06	0.0366002765641607\\
4.03701725859656e-06	0.0630610285821761\\
5.33669923120631e-06	0.107527635006195\\
6.13590727341318e-06	0.139563011334733\\
7.05480231071865e-06	0.180111929292801\\
8.11130830789687e-06	0.23072874996006\\
9.3260334688322e-06	0.292767464790904\\
1.07226722201032e-05	0.366996569289899\\
1.23284673944207e-05	0.453046208054237\\
1.4174741629268e-05	0.548746620555902\\
1.62975083462064e-05	0.649561877666391\\
1.87381742286038e-05	0.748517755990059\\
2.15443469003188e-05	0.837135429326993\\
2.47707635599171e-05	0.907635283521017\\
2.8480358684358e-05	0.955827234916771\\
3.27454916287773e-05	0.982996548580511\\
3.76493580679247e-05	0.995020361155428\\
4.32876128108306e-05	0.998961552419905\\
6.57933224657568e-05	0.999999435514274\\
1	1\\
};
\addlegendentry{$k = 40$}

\addplot [color=mycolor1, line width=1.5pt]
  table[row sep=crcr]{%
1e-06	0.00796776773912809\\
1.74752840000768e-06	0.0241313233344737\\
2.65608778294668e-06	0.0548567089446212\\
3.51119173421513e-06	0.0938652683180131\\
4.03701725859656e-06	0.12214536383851\\
4.64158883361278e-06	0.158163740003327\\
5.33669923120631e-06	0.203493077722365\\
6.13590727341318e-06	0.259648188536647\\
7.05480231071865e-06	0.327783551512027\\
8.11130830789687e-06	0.408221743861987\\
9.3260334688322e-06	0.499822141141715\\
1.07226722201032e-05	0.599306656709241\\
1.23284673944207e-05	0.700841549476151\\
1.4174741629268e-05	0.796370387540284\\
1.62975083462064e-05	0.877193122415294\\
1.87381742286038e-05	0.936756680947725\\
2.15443469003188e-05	0.973475131619499\\
2.47707635599171e-05	0.991468759149758\\
2.8480358684358e-05	0.998048766824902\\
3.76493580679247e-05	0.999975203196978\\
1	1\\
};
\addlegendentry{$k = 80$}

\end{axis}
\end{tikzpicture}%
%
%
\definecolor{mycolor1}{rgb}{0.00000,1.00000,1.00000}%
\begin{tikzpicture}

\begin{axis}[%
width=6.028in,
height=4.719in,
at={(1.011in,0.676in)},
scale only axis,
xmode=log,
xmin=1e-06,
xmax=1,
xlabel style={font=\color{white!15!black}, scale=2},
xlabel={$\sigma_{i}$},
ymode=log,
ymin=0.001,
ymax=1,
axis background/.style={fill=white},
title style={font=\bfseries, scale=1.1},
title={\normalsize\colorbox{black}{\color{white}8 | Stabilized filter factors $\phi_{(\gamma,i)}^{[k]}$ given in \eqref{FilFac} with $\alpha = 1$ and $\gamma = 10^{8}$}},
xmajorgrids,
ymajorgrids,
legend style={legend cell align=left, align=left, draw=white!15!black, font=\Large},
grid style = loosely dashed,
]
\addplot [color=red, line width=1.5pt]
  table[row sep=crcr]{%
1e-06	0.00100944018047136\\
4.64158883361278e-06	0.0213366708901693\\
7.05480231071865e-06	0.0485016625283951\\
9.3260334688322e-06	0.0830406236174168\\
1.07226722201032e-05	0.108122938510329\\
1.23284673944207e-05	0.140127280299613\\
1.4174741629268e-05	0.18051007646729\\
1.62975083462064e-05	0.230724099656778\\
1.87381742286038e-05	0.29197998731902\\
2.15443469003188e-05	0.364869189120445\\
2.47707635599171e-05	0.448850962202457\\
2.8480358684358e-05	0.541683935832261\\
3.27454916287773e-05	0.639012470489089\\
3.76493580679247e-05	0.734456232384665\\
4.32876128108306e-05	0.820564515007259\\
4.97702356433211e-05	0.890691101798097\\
5.72236765935022e-05	0.941176119469788\\
6.57933224657568e-05	0.972607745346745\\
7.56463327554629e-05	0.98917340274143\\
8.69749002617783e-05	0.996425367590102\\
0.000114975699539774	0.999780725283436\\
0.000613590727341318	1\\
1	1\\
};
\addlegendentry{$k = 10$}

\addplot [color=green, line width=1.5pt]
  table[row sep=crcr]{%
1e-06	0.00201786139146494\\
3.51119173421513e-06	0.0244091310622361\\
5.33669923120631e-06	0.0553931711930581\\
7.05480231071865e-06	0.0946509137887719\\
8.11130830789687e-06	0.123061035662332\\
9.3260334688322e-06	0.159185502064064\\
1.07226722201032e-05	0.20455530718855\\
1.23284673944207e-05	0.260618905915059\\
1.4174741629268e-05	0.328436265228354\\
1.62975083462064e-05	0.408214589151125\\
1.87381742286038e-05	0.498707661643224\\
2.15443469003188e-05	0.596608853071478\\
2.47707635599171e-05	0.696234738134843\\
2.8480358684358e-05	0.789946385325793\\
3.27454916287773e-05	0.86968800353761\\
3.76493580679247e-05	0.929486507480653\\
4.32876128108306e-05	0.967802906725421\\
4.97702356433211e-05	0.988051564773886\\
5.72236765935022e-05	0.996539751079367\\
6.57933224657568e-05	0.999249664385012\\
0.000114975699539774	0.999999951918599\\
1	1\\
};
\addlegendentry{$k = 20$}

\addplot [color=blue, line width=1.5pt]
  table[row sep=crcr]{%
1e-06	0.0040316510183346\\
2.31012970008316e-06	0.0212054428566986\\
3.51119173421513e-06	0.0482224564452588\\
4.64158883361278e-06	0.0826538095353632\\
5.33669923120631e-06	0.107717938971293\\
6.13590727341318e-06	0.139773957401092\\
7.05480231071865e-06	0.180343032096495\\
8.11130830789687e-06	0.230978052826379\\
9.3260334688322e-06	0.29303098006074\\
1.07226722201032e-05	0.367267740678098\\
1.23284673944207e-05	0.453315597709756\\
1.4174741629268e-05	0.549002150139558\\
1.62975083462064e-05	0.649790027506429\\
1.87381742286038e-05	0.748705991504796\\
2.15443469003188e-05	0.837275582579694\\
2.47707635599171e-05	0.907726665683993\\
2.8480358684358e-05	0.9558774789623\\
3.27454916287773e-05	0.983018783577985\\
3.76493580679247e-05	0.995027847372725\\
4.32876128108306e-05	0.998963347184669\\
6.57933224657568e-05	0.999999436996464\\
1	1\\
};
\addlegendentry{$k = 40$}

\addplot [color=mycolor1, line width=1.5pt]
  table[row sep=crcr]{%
1e-06	0.00804704782673553\\
2.00923300256505e-06	0.0319294959747747\\
2.65608778294668e-06	0.0550573177716756\\
3.51119173421513e-06	0.0941195075849028\\
4.03701725859656e-06	0.122428548357326\\
4.64158883361278e-06	0.158475966840018\\
5.33669923120631e-06	0.203832723566362\\
6.13590727341318e-06	0.260011155634623\\
7.05480231071865e-06	0.328162454967233\\
8.11130830789687e-06	0.408605244765292\\
9.3260334688322e-06	0.500194804846123\\
1.07226722201032e-05	0.5996498880134\\
1.23284673944207e-05	0.701136164292559\\
1.4174741629268e-05	0.796600939421261\\
1.62975083462064e-05	0.877352975166054\\
1.87381742286038e-05	0.936851321294413\\
2.15443469003188e-05	0.973520763975224\\
2.47707635599171e-05	0.991485631774208\\
2.8480358684358e-05	0.998053203137278\\
3.76493580679247e-05	0.999975277698251\\
1	1\\
};
\addlegendentry{$k = 80$}

\end{axis}
\end{tikzpicture}%
%
%
\definecolor{mycolor1}{rgb}{0.00000,1.00000,1.00000}%
\begin{tikzpicture}

\begin{axis}[%
width=6.028in,
height=4.719in,
at={(1.011in,0.676in)},
scale only axis,
xmode=log,
xmin=1e-06,
xmax=1,
xlabel style={font=\color{white!15!black}, scale=2},
xlabel={$\sigma_{i}$},
ymode=log,
ymin=0.01,
ymax=1,
axis background/.style={fill=white},
title style={font=\bfseries, scale=1.1},
title={\normalsize\colorbox{black}{\color{white}9 | Stabilized filter factors $\phi_{(\gamma,i)}^{[k]}$ given in \eqref{FilFac} with $\alpha = 10^{3}$ and $\gamma = 10^{8}$}},
xmajorgrids,
ymajorgrids,
legend style={legend cell align=left, align=left, draw=white!15!black, font=\Large},
grid style = loosely dashed,
]
\addplot [color=red, line width=1.5pt]
  table[row sep=crcr]{%
1e-06	0.0109446203635166\\
1.32194114846603e-06	0.014863956934546\\
1.51991108295293e-06	0.0173679439738397\\
1.74752840000768e-06	0.0203343709572257\\
2.00923300256505e-06	0.0238594968581399\\
2.31012970008316e-06	0.0280620338414962\\
2.65608778294668e-06	0.0330887141829384\\
3.05385550883341e-06	0.0391211723670256\\
3.51119173421513e-06	0.0463843567888646\\
4.03701725859656e-06	0.0551566405504001\\
4.64158883361278e-06	0.0657816763189914\\
5.33669923120631e-06	0.0786817754875455\\
6.13590727341318e-06	0.0943720920068643\\
8.11130830789687e-06	0.136724728727525\\
1.4174741629268e-05	0.28943372517131\\
1.62975083462064e-05	0.347184253789481\\
1.87381742286038e-05	0.413894292739998\\
2.15443469003188e-05	0.489061985735209\\
2.47707635599171e-05	0.571013368651488\\
2.8480358684358e-05	0.656596575192757\\
3.27454916287773e-05	0.741153653267332\\
3.76493580679247e-05	0.819018981284653\\
4.32876128108306e-05	0.884678494435437\\
4.97702356433211e-05	0.934361602444825\\
5.72236765935022e-05	0.967349572661665\\
6.57933224657568e-05	0.986121478619047\\
7.56463327554629e-05	0.99506596520806\\
8.69749002617783e-05	0.998559755252602\\
0.000114975699539774	0.99993527979642\\
0.00376493580679246	1\\
1	1\\
};
\addlegendentry{$k = 10$}

\addplot [color=green, line width=1.5pt]
  table[row sep=crcr]{%
1e-06	0.0217694560121315\\
1.32194114846603e-06	0.0295069766533399\\
1.74752840000768e-06	0.0402552552722251\\
2.00923300256505e-06	0.0471497181259562\\
2.31012970008316e-06	0.0553365899396712\\
2.65608778294668e-06	0.0650825653595967\\
3.05385550883341e-06	0.0767118786066808\\
3.51119173421513e-06	0.0906172050230124\\
4.64158883361278e-06	0.127236123698646\\
6.13590727341318e-06	0.179838092263977\\
9.3260334688322e-06	0.302737663135967\\
1.07226722201032e-05	0.358695561104102\\
1.23284673944207e-05	0.422933114019888\\
1.4174741629268e-05	0.495095569076078\\
1.62975083462064e-05	0.573831601499604\\
1.87381742286038e-05	0.656480099917253\\
2.15443469003188e-05	0.738942345579153\\
2.47707635599171e-05	0.815970470124257\\
2.8480358684358e-05	0.882074087830656\\
3.27454916287773e-05	0.932998568783152\\
3.76493580679247e-05	0.967245870864755\\
4.32876128108306e-05	0.986700950354322\\
4.97702356433211e-05	0.995691600766389\\
5.72236765935022e-05	0.998933949594624\\
7.56463327554629e-05	0.999975655300672\\
0.327454916287773	1\\
1	1\\
};
\addlegendentry{$k = 20$}

\addplot [color=blue, line width=1.5pt]
  table[row sep=crcr]{%
1e-06	0.043065002809199\\
1.32194114846603e-06	0.0581432916354591\\
1.74752840000768e-06	0.078890024967418\\
2.31012970008316e-06	0.107611041693191\\
3.05385550883341e-06	0.147539044893996\\
6.13590727341318e-06	0.327334445098807\\
7.05480231071865e-06	0.382316672948253\\
8.11130830789687e-06	0.444611091323354\\
9.3260334688322e-06	0.513825233590907\\
1.07226722201032e-05	0.588728616652416\\
1.23284673944207e-05	0.666993809105216\\
1.4174741629268e-05	0.745071515633392\\
1.62975083462064e-05	0.818380496119607\\
1.87381742286038e-05	0.881994078247139\\
2.15443469003188e-05	0.931848901068286\\
2.47707635599171e-05	0.966133132133713\\
2.8480358684358e-05	0.986093479239029\\
3.27454916287773e-05	0.995510808214895\\
3.76493580679247e-05	0.998927167024591\\
4.97702356433211e-05	0.999981437696043\\
1	1\\
};
\addlegendentry{$k = 40$}

\addplot [color=mycolor1, line width=1.5pt]
  table[row sep=crcr]{%
1e-06	0.0842754111514417\\
1.51991108295293e-06	0.130784578701674\\
3.51119173421513e-06	0.316108929332654\\
4.03701725859656e-06	0.364846755567994\\
4.64158883361278e-06	0.419787541313086\\
5.33669923120631e-06	0.480868602817124\\
6.13590727341318e-06	0.54752105124947\\
7.05480231071865e-06	0.618467307482285\\
8.11130830789687e-06	0.691543160118964\\
9.3260334688322e-06	0.763634096507064\\
1.07226722201032e-05	0.830855849239362\\
1.23284673944207e-05	0.889106876825747\\
1.4174741629268e-05	0.935011467858545\\
1.62975083462064e-05	0.96701435581024\\
1.87381742286038e-05	0.986074602431258\\
2.15443469003188e-05	0.995355427714402\\
2.47707635599171e-05	0.998853035260929\\
3.27454916287773e-05	0.999979847157116\\
0.572236765935022	1\\
1	1\\
};
\addlegendentry{$k = 80$}

\end{axis}
\end{tikzpicture}%
}     
    \caption{\footnotesize The stabilized filter factors $\phi_{(\gamma,i)}^{[k]}$ represented as functions of the singular values $\sigma_{i} \in [10^{-6}, \,10^{0}]$ (in {\tiny\colorbox{black}{\color{white}1}}, {\tiny\colorbox{black}{\color{white}2}}, {\tiny\colorbox{black}{\color{white}3}}), for each of the stepsize $\alpha = 10^{i}, i = -3, 0, 3$ and with respect to the parameter $\gamma = 10^{i}, i = 3, 6, 8, 10$. The same filter factors are represented in {\tiny\colorbox{black}{\color{white}4}}, {\tiny\colorbox{black}{\color{white}5}}, {\tiny\colorbox{black}{\color{white}6}}, {\tiny\colorbox{black}{\color{white}7}}, {\tiny\colorbox{black}{\color{white}8}}, {\tiny\colorbox{black}{\color{white}9}}, for each $\alpha = 10^{i}, i = -3, 0, 3$ and with respect to the iterations $k = 10, 20, 40, 80$. The figure is displayed on a log-log scale.}
    \label{FiltFacs} 
\end{figure}

The behavior of filter factors in \eqref{FilFac} is illustrated in Figure \ref{FiltFacs} (see {\tiny\colorbox{black}{\color{white}1}}, {\tiny\colorbox{black}{\color{white}2}}, {\tiny\colorbox{black}{\color{white}3}}) as a function of the singular values $\sigma_{i} \in [10^{-6}, \, 10^{0}]$ (for simplicity), for each of the stepsize $\alpha = 10^{i}, i = -3, 0, 3$ and with respect to the stabilization parameter $\gamma = 10^{i}, i = 3, 6, 8, 10$. We see that as $\gamma$ increases, the filter factors are shifted toward the unit value. More SVD components are thus effectively included in the iteration $\bm{x}^{[k]}_{\gamma}$ in \eqref{SVDSol} as $\gamma$ increases, and the limit asymptotic behavior of the filter factors when $\gamma \to \infty$ becomes
\begin{align}\label{FilFacBis1} 
     \phi_{(\gamma,i)}^{[k]} \longrightarrow \phi_{(\infty,i)} = 1, \; i = 1, \cdots, n,
\end{align}
provides the the naïve solution $\bm{x}^{\star}_{\text{naïve}}$ of the linear system \eqref{Eq1}:
\begin{remark}\label{RMK1} \text{ } 
\begin{itemize}
    \item The SVD expansion reached by taking $\gamma \to \infty$ in \eqref{SVDSol} is the naïve solution $\bm{x}^{\star}_{\text{naive}}$ of \eqref{Eq1}:
\begin{align}\label{SVD1}
    \bm{x}^{[k]}_{\infty} & = \sum^{n}_{i = 1}\left(\frac{\bm{v}^{t}_{i}\bm{b}}{\sigma_{i}}\right)\bm{v}_{i}, \, \forall \; k = 1, 2, 3, \cdots \nonumber\\
    & \equiv \bm{x}^{\star}_{\text{naive}}.
\end{align}
    \item In a context of noisy coefficients encountered among the values $\bm{v}^{t}_{i}\bm{b}$, this solution may not be accurate when the stabilization parameter $\gamma$ is very large, specially when $\mathbf{A}$ is close to a singular matrix (i.e. it has small singular values $\sigma_{i}$, compared to the largest singular value). 
\end{itemize}
\end{remark}

We also represent in Figure \ref{FiltFacs} (in {\tiny\colorbox{black}{\color{white}4}}, {\tiny\colorbox{black}{\color{white}5}}, {\tiny\colorbox{black}{\color{white}6}} with $\gamma = 10^{6}$, and in {\tiny\colorbox{black}{\color{white}7}}, {\tiny\colorbox{black}{\color{white}8}}, {\tiny\colorbox{black}{\color{white}9}} with $\gamma = 10^{8}$) the filter factors in \eqref{FilFac} as functions of the singular values $\sigma_{i}$, for each of the stepsize $\alpha = 10^{i}, i = -3, 0, 3$ and with respect to different fixed iterations $k = 10, 20, 40, 80$ (with an increasing iteration number). We also note as the number of iteration $k$ increases, the filter factors are more and more shifted toward larger values (more near than one). The SVD components are thus more and more included in the iteration $\bm{x}^{[k]}_{\gamma}$ provided in \eqref{SVDSol}.

\section{Convergence Analysis of the Stabilized Gradient Method}\label{sec3}
\begin{tcolorbox}[breakable, colback=gray!5!]
\begin{lemma}\label{Lem1} Let $\mathbf{A} \in \mathbb{R}^{m\times n}$ be a matrix of $\text{rank}(\mathbf{A}) = r \leq \min(m, n)$, whose the Singular Value Decomposition (SVD) is $\mathbf{A} = \mathbf{U}\mathbf{\Sigma}\mathbf{V}^{t}$, where $\mathbf{U} \in \mathbb{R}^{m\times m}$ and $\mathbf{V} \in \mathbb{R}^{n\times n}$ are orthogonal matrices, and $\mathbf{\Sigma} \in \mathbb{R}^{m\times n}$ is a matrix defined as follows 
\begin{align*}
    \mathbf{\Sigma} = \begin{bmatrix}
\mathbf{\Sigma}_{1} & \mathbf{0} \\
\mathbf{0}          & \mathbf{0} \\
\end{bmatrix},
\end{align*}
where $\mathbf{\Sigma}_{1} = \text{diag}(\sigma_{1}, \sigma_{2}, \cdots, \sigma_{r})$ is a diagonal matrix with entries $\sigma_1 \geq \sigma_2 \geq \cdots \geq \sigma_r > 0$. We obtain the following estimate  
\begin{align}\label{Lim1}
    \left\|(\mathbf{I} + \gamma\mathbf{A}^{t}\mathbf{A})^{-1}\right\| \leq \left(\frac{\kappa\left(\mathbf{V}\right)}{1 + \gamma\sigma^{2}_{r}}\right), \; \forall \; \gamma > 0,
\end{align}
where $\kappa\left(\mathbf{V}\right) =\|\mathbf{V}\|\|\mathbf{V}^{-1}\|$ is the condition number of $\mathbf{V}$.
\end{lemma}
\end{tcolorbox}
\begin{proof}
From the SVD of $\mathbf{A} = \mathbf{U}\mathbf{\Sigma}\mathbf{V}^{t}$ and the orthogonality of the matrices $\mathbf{U}$ and $\mathbf{V}$, we obtain   
\begin{align*}
    \mathbf{I} + \gamma\mathbf{A}^{t}\mathbf{A} & =  \mathbf{V}\left(\mathbf{I} + \gamma\mathbf{\Sigma}^{2}\right)\mathbf{V}^{-1}, \; \forall \; \gamma > 0.
\end{align*}
So, the following inverse is obtained
\begin{align*}
    \left(\mathbf{I} + \gamma\mathbf{A}^{t}\mathbf{A}\right)^{-1} & = \mathbf{V}\left(\mathbf{I} + \gamma\mathbf{\Sigma}^{2}\right)^{-1}\mathbf{V}^{-1}, \; \forall \; \gamma > 0.
\end{align*}
We thus establish the following estimate by using the sub-multiplicity of the matrix norm
\begin{align}\label{Est12}
    \left\|\left(\mathbf{I} + \gamma\mathbf{A}^{t}\mathbf{A}\right)^{-1}\right\| = \left\|\mathbf{V}\left(\mathbf{I} + \gamma\mathbf{\Sigma}^{2}\right)^{-1}\mathbf{V}^{-1}\right\| \leq \left\|\mathbf{V}\right\|\left\|\mathbf{V}^{-1}\right\|\left\|\left(\mathbf{I} + \gamma\mathbf{\Sigma}^{2}\right)^{-1}\right\|, \; \forall \; \gamma > 0.
\end{align}
Since the induced matrix norm $\left\|\cdot\right\|$ (from the euclidean vector norm in $\mathbb{R}^{n}$) of a diagonal matrix is equal to the absolute value of the largest diagonal element, then the norm of the matrix $\left(\mathbf{I} + \gamma\mathbf{\Sigma}^{2}\right)^{-1} = \text{diag}\left(\frac{1}{1 + \gamma\sigma^2_{1}}, \frac{1}{1 + \gamma\sigma^2_{2}}, \cdots, \frac{1}{1 + \gamma\sigma^2_{r}}\right)$ is $\frac{1}{1 + \gamma\sigma^{2}_{r}}$. Hence, we establish from \eqref{Est12} the estimate \eqref{Lim1}, where $\kappa\left(\mathbf{V}\right) =\|\mathbf{V}\|\|\mathbf{V}^{-1}\|$.
\end{proof}

\begin{tcolorbox}[breakable, colback=gray!5!]
\begin{theorem}\label{CorMSN1} Let $\bm{x}^{\star}$ be the exact solution of \eqref{Eq1}. For any $\alpha_{k} \in \mathbb{R}$, the sequence $\left\{\bm{x}^{[k]}_{\gamma}\right\}_{k \in \mathbb{N}^{\star}}$ generated from the stabilized gradient method \eqref{EqP1} verifies
\begin{align}\label{Reslt1}
    \|\bm{x}^{[k+1]}_{\gamma} - \bm{x}^{\star}\|_{n} \leq \left\|(\mathbf{I} + \gamma\mathbf{A}^{t}\mathbf{A})^{-1}\right\|\left\|\mathbf{I} - \alpha_{k}\mathbf{A}\right\|\|\bm{x}^{[k]}_{\gamma} - \bm{x}^{\star}\|_{n}, \, \forall \, \gamma > 0, \, \forall \, k \in \mathbb{N}^{\star},
\end{align}
\end{theorem}
\end{tcolorbox}
\begin{proof}
The gradient iterative method \eqref{Eq2} is consistent with the linear system \eqref{Eq1}, that is to say the exact solution $\bm{x}^{\star}$ satisfies the following equation (see \cite{QuarteroniSaccoSaleri07} - Definition $4.1$)
\begin{align}\label{Eq2Bis1}
    \bm{x}^{\star} = \bm{x}^{\star} - \alpha_{k}(\mathbf{A}\bm{x}^{\star} - \bm{b}), \; k = 0, 1, 2, \cdots,
\end{align}  
 Otherwise from the system \eqref{Eq1}, we also obtain  the normal equation satisfied by the exact solution 
 \begin{align}\label{Eq2Bis2}
    \gamma\mathbf{A}^t\mathbf{A}\bm{x}^{\star} = \gamma\mathbf{A}^{t}\bm{b}, \; \forall \; \gamma > 0. 
\end{align}  
Then adding both equations \eqref{Eq2Bis1} and \eqref{Eq2Bis2}, we thus remark that the stabilized gradient method \eqref{EqP1} is also consistent with the linear system \eqref{Eq1}
\begin{align}\label{EqP2}
    \left(\mathbf{I} + \gamma\mathbf{A}^t\mathbf{A}\right)\bm{x}^{\star} = \left(\mathbf{I} - \alpha_{k}\mathbf{A}\right)\bm{x}^{\star} + \alpha_{k}\bm{b} + \gamma\mathbf{A}^{t}\bm{b},  \, \forall \, \gamma > 0,
\end{align}
Denoting the error at the $(k+1)^{\text{th}}$ iteration by $\bm{e}^{[k+1]}_{\gamma} = \bm{x}^{[k+1]}_{\gamma} - \bm{x}^{\star}$, we obtain by subtracting equations \eqref{EqP1} and \eqref{EqP2} the equality
\begin{align*}
    \left(\mathbf{I} + \gamma\mathbf{A}^t\mathbf{A}\right)\bm{e}^{[k+1]}_{\gamma} = \left(\mathbf{I} - \alpha_{k}\mathbf{A}\right)\bm{e}^{[k]}_{\gamma},  \, \forall \, \gamma > 0, \, \forall \, k = 0, 1, 2, 3, \cdots.
\end{align*}
The matrix $\left(\mathbf{I} + \gamma\mathbf{A}^t\mathbf{A}\right)$ being nonsingular for any $\gamma > 0$, we then obtain the recurrence relation 
\begin{align*}
    \bm{e}^{[k+1]}_{\gamma} = \left(\mathbf{I} + \gamma\mathbf{A}^t\mathbf{A}\right)^{-1}\left(\mathbf{I} - \alpha_{k}\mathbf{A}\right)\bm{e}^{[k]}_{\gamma},  \, \forall \, \gamma > 0, \, \forall \, k = 0, 1, 2, 3, \cdots,
\end{align*}
which established the result \eqref{Reslt1} by using the matrix norm sub-multiplicity, and its compatibility with the euclidean vector norm.
\end{proof}

\begin{tcolorbox}[breakable, colback=gray!5!]
\begin{corollary}\label{CorMSN12} Let $\mathbf{A} \in \mathbb{R}^{n\times n}$ be a nonsingular $n\times n$ matrix, and $\bm{x}^{[0]}_{\gamma}$ any initial datum. The sequence $\left\{\bm{x}^{[k]}_{\gamma}\right\}_{k \in \mathbb{N}^{\star}}$, computed from the stabilized gradient method \eqref{EqP1}, satisfies
\begin{itemize}
    \item If the stepsize $\alpha_{k}$ is a non-constant parameter:
    \begin{align}\label{Reslt1Fin2}
        \|\bm{x}^{[k]}_{\gamma} - \bm{x}^{\star}\|_{n} \leq \left(\frac{\kappa\left(\mathbf{V}\right)}{1 + \gamma\sigma^{2}_{n}}\right)^{k}\left(\prod_{i=0}^{k-1}\left\|\mathbf{I} - \alpha_{i}\mathbf{A}\right\|\right)\|\bm{x}^{[0]}_{\gamma} - \bm{x}^{\star}\|_{n}, \, \forall \, \gamma > 0, \, \forall \, k \in \mathbb{N}^{\star}.
    \end{align}
    \item If the stepsize is a constant parameter, namely $\alpha = \alpha_{k}$:
    \begin{align}\label{Reslt1FinBis1}
        \|\bm{x}^{[k]}_{\gamma} - \bm{x}^{\star}\|_{n} \leq \left(\left(\frac{\kappa\left(\mathbf{V}\right)}{1 + \gamma\sigma^{2}_{n}}\right)\left\|\mathbf{I} - \alpha\mathbf{A}\right\|\right)^{k}\|\bm{x}^{[0]}_{\gamma} - \bm{x}^{\star}\|_{n}, \, \forall \, \gamma > 0, \, \forall \, k \in \mathbb{N}^{\star}.
    \end{align}
\end{itemize}
\end{corollary}
\end{tcolorbox}
\begin{proof}
Multiplying term by term the inequality \eqref{Reslt1}, for $k = 0, \cdots, k-1$, by using the fact that the $n\times n$ matrix $\mathbf{A} \in \mathbb{R}^{n\times n}$ is nonsingular (and thus of rank $n$) in the Lemma \ref{Lem1}, we obtain the estimate \eqref{Reslt1Fin2} after simplification. If the stepsize is constant, we thus arrive from \eqref{Reslt1Fin2} to the estimate \eqref{Reslt1FinBis1}.
\end{proof}

\begin{remark} \text{}  
\begin{itemize}
    \item We may see from estimates \eqref{Reslt1Fin2} and \eqref{Reslt1FinBis1}, how our approach stabilizes the gradient method by damping the smallest singular value $\sigma^{2}_{n}$ by the stabilization parameter $\gamma > 0$, but also by divising the whole right-hand side term by this parameter (which can be chosen larger).  
    \item The convergence of the stabilized gradient method \eqref{EqP1} is faster than the parameter $\gamma$ increases.   
    \item We thus obtain, across the stabilization term $\gamma$, a total control of the convergence of the approximate solution $\bm{x}^{[k]}_{\gamma}$ independently to the stepsize $\alpha_{k}$ and the structure of the matrix $\mathbf{A}$ whose the only required criteria is to be nonsingular (for the obtaining of a unique solution).
    \item With respect to the estimate \eqref{Reslt1FinBis1}, the error $\bm{x}^{[k]}_{\gamma} - \bm{x}^{\star}$ converges toward zero if and only if 
    \begin{align}\label{Rm31}
        \left(\frac{\kappa\left(\mathbf{V}\right)}{1 + \gamma\sigma^{2}_{n}}\right)\left\|\mathbf{I} - \alpha\mathbf{A}\right\| < 1. 
    \end{align}
    If $\mathbf{A}$ is symmetric, then its singular values are equal to the absolute values of its eigenvalues $\sigma_{i} = |\lambda_{i}|$, where $\lambda_{i}$, $i = 1, \cdots, n$, are eigenvalues of $\mathbf{A}$. We thus obtain from \eqref{Rm31}
    \begin{align}\label{Err1}
        \left\|\mathbf{I} - \alpha\mathbf{A}\right\| & < \frac{1 + \gamma\sigma^{2}_{n}}{\kappa\left(\mathbf{V}\right)}   \nonumber\\
        \max_{i=1, \cdots, n}|1 - \alpha\lambda_{i}| & < \frac{1 + \gamma\lambda^{2}_{n}}{\kappa\left(\mathbf{V}\right)}   \nonumber\\
        1 - \frac{1 + \gamma\lambda^{2}_{n}}{\kappa\left(\mathbf{V}\right)} < \alpha\lambda_{i} & < 1 + \frac{1 + \gamma\lambda^{2}_{n}}{\kappa\left(\mathbf{V}\right)}, \; \forall \; \gamma > 0, \; \forall \; i = 1, \cdots, n.
    \end{align}
    We note from \eqref{Err1} that, the stepsize $\alpha$ and eigenvalues $\lambda_{i}$ may be of different signs when 
    \begin{align}\label{Err1Bis}
        1 - \frac{1 + \gamma\lambda^{2}_{n}}{\kappa\left(\mathbf{V}\right)} < 0 \iff \gamma > \frac{\kappa\left(\mathbf{V}\right) - 1}{\lambda^{2}_{n}} > 0.
    \end{align}
    Thus the convergence of the stabilized method \eqref{EqP1} is no longer based on the idea that the eigenvalues of the matrix $\mathbf{A}$ and the parameter $\alpha$ should be of the same sign, as it is the case for the gradient method \eqref{Eq2} (see \cite{AllaireKaber08} - Chap. 9, Theorem $9.1.1$). Moreover, it is also no longer based on the choice of the stepsize $\alpha$.
    \item The rate of convergence of the stabilized gradient method \eqref{EqP1} is not governed by the eigenvalues of the matrix $\mathbf{A}$ (\cite{Hackbusch10} - Chap. 9, Theorem 9.10). This convergence is guided by the choice of the stabilization parameter $\gamma$, independently of the linear system data $\mathbf{A}$ and $\bm{b}$. If a very large value of the stabilization term is chosen, the stabilized gradient method \eqref{EqP1} should behave like a direct method by making one iteration. 
\end{itemize}
\end{remark}

The convergence criterion of the sequence $\left\{\bm{x}^{[k]}_{\gamma}\right\}_{k \in \mathbb{N}}$ is determined by the relative criterion 
\begin{align}\label{CTR1}
    \frac{\|\bm{b} - \mathbf{A}\bm{x}^{[k]}_{\gamma}\|_{n}}{\|\bm{b} - \mathbf{A}\bm{x}^{[0]}_{\gamma}\|_{n}} \leq \varepsilon, \; \forall \; k = 1, 2, 3, \cdots.
\end{align}
In the building of the algorithm, an initial value $\bm{x}$, the matrix $\mathbf{A}$, the right hand side $\bm{b}$, a number of maximum iterations $k_{\text{max}}$, a stepsize parameter $\alpha$ and the stabilization parameter $\gamma$ are used as inputs. The algorithm returns as output the final value $\bm{x}$ at the end of the $k_{\text{max}}$ iterations (when the stopping condition is not verified) or before. We thus present the following general schematic algorithm of the stabilized gradient method \eqref{EqP1}.
\begin{algorithm}
\caption{Stabilized Gradient (SG) Method for Solving $\mathbf{A}\bm{x} = \bm{b}$}\label{alg1}
\begin{algorithmic}[1]
\Require $\mathbf{A} \in \mathbb{R}^{n\times n}$, $\bm{b} \in \mathbb{R}^n$, $\bm{x} \in \mathbb{R}^n$, $\varepsilon > 0$, $k_{\text{max}}$, $\gamma$, $\alpha$ 
\State $\mathbf{M} = \left(\mathbf{I} + \gamma\mathbf{A}^{t}\mathbf{A}\right)$,\; $\mathbf{N} = \left(\mathbf{I} - \alpha\mathbf{A}\right)$,\; $\bm{c} = \alpha\bm{b} + \gamma\mathbf{A}^{t}\bm{b}$\; 
\State $\bm{r}_{0} = \mathbf{A}\bm{x} - \bm{b}$,\; $\bm{r} = \bm{r}_{0}$,\; $k = 0$\;
\While{($\|\bm{r}\|_{n} \geq \varepsilon\|\bm{r}_{0}\|_{n}$ \text{ and } $k \leq k_{\text{max}}$)}
    \State $k = k + 1$ 
    \State \text{Solve} $\mathbf{M}\bm{x} = \mathbf{N}\bm{x} + \bm{c}$  \Comment{The matrix $\mathbf{M}$ is symmetric and nonsingular.}  
    \State $\bm{r} = \mathbf{A}\bm{x} - \bm{b}$  
\EndWhile
\Ensure $\bm{x}$
\end{algorithmic}
\end{algorithm}

\section{Numerical Examples}\label{sec5}
In this section, we use three sets of numerical examples formulated under the linear system \eqref{Eq1}, in order to test the stabilized gradient method \eqref{EqP1}. The first example in the sub-section \ref{NEx1} focus on well-conditioned linear systems, whose the matrices are non-symmetric and not positive definite. Contrary to the first example, the second numerical example investigated in the sub-section \ref{NEx2} is with symmetric and positive definite matrices. However, these matrices are defined with respect to a parameter that makes them increasingly ill-conditioned. As for the last set of examples investigated in the sub-section \ref{NEx3}, where we distinguish small and large scale problems, we mainly focus on ill-conditioned linear systems. 

For each of these examples, we start the algorithm \ref{alg1} by the initial vector $\bm{x}^{[0]}_{\gamma} = \bm{0}$, we fix the stopping tolerance to $\varepsilon = 10^{-5}$, and the largest iteration step to $k_{\text{max}} = n$ (any change of these values will be mentioned accordingly). We only face through the forthcoming examples, quadrature errors from the evaluation of coefficients of the matrix $\mathbf{A}$ and the right-hand side $\mathbf{b}$ (in subsection \ref{NEx3}), and rounding errors. 

\subsection{\bf Well-conditioned linear systems}\label{NEx1}
We consider here the two non-symmetric and not positive definite $4\times 4$ matrices $\mathbf{A}_{1}$ and $\mathbf{A}_{2}$, provided in the Chapter $8$ in \cite{AllaireKaber08} as an exercise (Exercise $8.2$, in section $8.7$) 
\begin{align}\label{Mats}
    \mathbf{A}_{1} = \left(\begin{array}{llll}
 1 & 2 & 3 & 4 \\
 4 & 5 & 6 & 7 \\
 4 & 3 & 2 & 0 \\
 0 & 2 & 3 & 4 
\end{array}\right), \; \mathbf{A}_{2} = \left(\begin{array}{llll}
 2 & 4 & -4 & 1 \\
 2 & 2 &  2 & 0 \\
 2 & 2 & 1  & 0 \\
 2 & 0 & 0  & 2
\end{array}\right),
\end{align}
where $\textbf{cond}(\mathbf{A}_{1}) = 171.62$, $\textbf{cond}(\mathbf{A}_{2}) = 35.37$, $\textbf{rank}(\mathbf{A}_{1}) = 4$ and $\textbf{rank}(\mathbf{A}_{2}) = 4$. Let us first note that for the matrix $\mathbf{A}_{1}$ the Gauss-Seidel method converges, but the Jacobi method diverges. Whereas for the matrix $\mathbf{A}_{2}$, it is the Jacobi method that converges and the Gauss-Seidel method diverges (see \cite{AllaireKaber07} - Exercise $2$ of Chap. $8$). 

\begin{figure}[h]
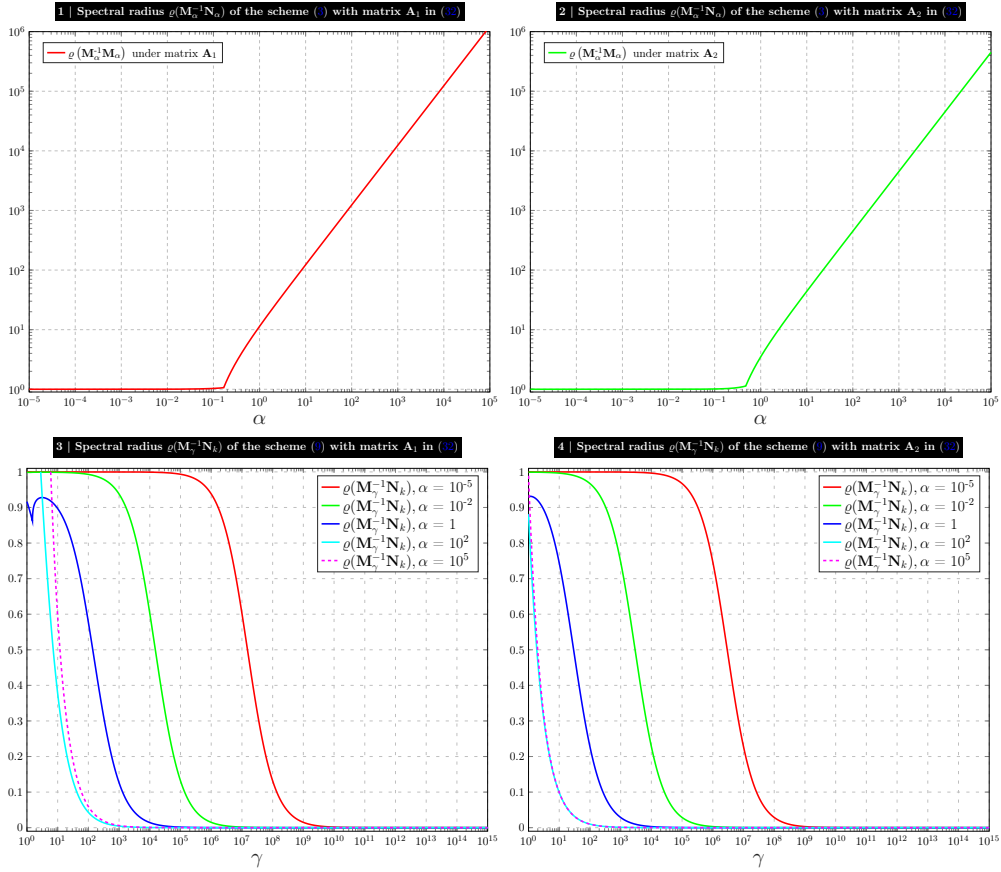
  
    \centering  
    \scalebox{0.4}{
%
%
\begin{tikzpicture}

\begin{axis}[%
width=6.028in,
height=4.719in,
at={(1.011in,0.676in)},
scale only axis,
xmode=log,
xmin=1e-05,
xmax=100000,
xlabel style={font=\color{white!15!black}, scale=1.8},
xlabel={$\alpha$},
ymode=log,
ymin=0.9,
ymax=1000000,
axis background/.style={fill=white},
title style={font=\bfseries, scale=1.},
title={\normalsize\colorbox{black}{\color{white}1 | Spectral radius $\varrho(\mathbf{M}^{-1}_{\alpha}\mathbf{N}_{\alpha})$ of the scheme \eqref{LIM1} with matrix $\mathbf{A}_{1}$ in \eqref{Mats}}},
xmajorgrids,
ymajorgrids,
grid style={dashed},
legend style={at={(0.03,0.97)}, anchor=north west, legend cell align=left, align=left, draw=white!15!black}
]
\addplot [color=red, line width=1.5pt]
  table[row sep=crcr]{%
1e-05	1.00000260432333\\
0.0170447812649629	1.00463792024583\\
0.0491778045881989	1.01444907118394\\
0.0904585328363252	1.02904447038087\\
0.138091307406656	1.04855318513779\\
0.167180816602493	1.06513133441345\\
0.186925281229698	1.30902841191015\\
0.211803456507603	1.61634058057595\\
0.243706349204311	2.01042684421141\\
0.285487991070698	2.52654214724093\\
0.341614931903498	3.21986035547879\\
0.419389300227613	4.18058233485659\\
0.531401831943166	5.56423743233577\\
0.700856085184714	7.65745180483651\\
0.974323662274939	11.0355096099851\\
1.456535424616	16.992117793094\\
2.42375687431713	28.9398960351632\\
4.80432475892626	58.3462922480367\\
13.3411098841158	163.798477585195\\
91.46287941198	1128.81179329006\\
100000	1235267.12249262\\
};
\addlegendentry{$\varrho\left(\mathbf{M}^{\text{-1}}_{\alpha}\mathbf{M}_{\alpha}\right) \text{ under matrix } \mathbf{A}_{1}$}


\end{axis}
\end{tikzpicture}%
%
%
\begin{tikzpicture}

\begin{axis}[%
width=6.028in,
height=4.719in,
at={(1.011in,0.676in)},
scale only axis,
xmode=log,
xmin=1e-05,
xmax=100000,
xminorticks=true,
xlabel style={font=\color{white!15!black}, scale=1.8},
xlabel={$\alpha$},
ymode=log,
ymin=0.9,
ymax=1000000,
axis background/.style={fill=white},
title style={font=\bfseries, scale=1.},
title={\normalsize\colorbox{black}{\color{white}2 | Spectral radius $\varrho(\mathbf{M}^{-1}_{\alpha}\mathbf{N}_{\alpha})$ of the scheme \eqref{LIM1} with matrix $\mathbf{A}_{2}$ in \eqref{Mats}}},
xmajorgrids,
ymajorgrids,
grid style={dashed},
legend style={at={(0.03,0.97)}, anchor=north west, legend cell align=left, align=left, draw=white!15!black}
]

\addplot [color=green, line width=1.5pt]
  table[row sep=crcr]{%
1e-05	1.00000053795718\\
0.056153526744846	1.00462100679573\\
0.122498660805739	1.01415641699306\\
0.191474278985062	1.02862270066726\\
0.263204763072762	1.04838466007946\\
0.338384876889366	1.07397244532121\\
0.417981101578928	1.10612057632236\\
0.474913097203419	1.13494089602378\\
0.533093082721273	1.39648523148972\\
0.606776353704442	1.72772357699146\\
0.701844738444652	2.1550973085141\\
0.827269378551249	2.71893560884622\\
0.997269070577719	3.48315815178295\\
1.2354476593386	4.55387448430186\\
1.5832851994288	6.11755549011921\\
2.11908634637669	8.52621464795568\\
3.00505728129102	12.5090392800311\\
4.62212681308062	19.7784700358077\\
8.0410073493267	35.147824804992\\
17.2507214404136	76.5494946463399\\
57.2257070878942	256.254438938947\\
770.852311214111	3464.31635723776\\
100000	449542.486712754\\
};
\addlegendentry{$\varrho\left(\mathbf{M}^{\text{-1}}_{\alpha}\mathbf{M}_{\alpha}\right) \text{ under matrix } \mathbf{A}_{2}$}

\end{axis}
\end{tikzpicture}%
}
    \scalebox{0.4}{
        \input{RayonSpectralA1.tex}    
        \input{RayonSpectralA2.tex}}   
    \caption{\footnotesize The spectral radii $\varrho(\mathbf{M}^{-1}_{\alpha}\mathbf{N}_{\alpha})$ and $\varrho(\mathbf{M}^{-1}_{\gamma}\mathbf{N}_{k})$ of iterative matrices $\mathbf{M}^{-1}_{\alpha}\mathbf{N}_{\alpha}$ and $\mathbf{M}^{-1}_{\gamma}\mathbf{N}_{k}$ from the gradient method \eqref{LIM1} and the stabilized gradient method \eqref{EqP1}, respectively. The former is represented with respect to the stepsize $\alpha \in [10^{-5}, \, 10^{5}]$, and the latter with respect to the stabilization parameter $\gamma \in [10^{0}, \, 10^{15}]$ for the different values of the fixed stepsize $\alpha_{k} = \alpha = 10^{i}, i = -5, -2, 0, 2, 5$.}
    \label{FigWellC1}
\end{figure}

The spectral radii of the iteration matrix $\mathbf{M}^{-1}_{\alpha}\mathbf{N}_{\alpha}$ (from the gradient method \eqref{LIM1}) with respect to matrices $\mathbf{A}_{1}$ and $\mathbf{A}_{2}$ in \eqref{Mats} are respectively represented at Figure \ref{FigWellC1} ({\small\colorbox{black}{\color{white}1}} and {\small\colorbox{black}{\color{white}2}}), as functions of the stepsize $\alpha$. In particular, for the fixed stepsize $ \alpha = 1$ we note
\begin{align*}
    \varrho\left(\mathbf{M}^{-1}_{\alpha}\mathbf{N}_{\alpha}\right) = 11.3527 > 1, \text{ and } \varrho\left(\mathbf{M}^{-1}_{\alpha}\mathbf{N}_{\alpha}\right) = 3.4954 > 1.
\end{align*}
Whereas, we see from the Figure \ref{FigWellC1} ({\small\colorbox{black}{\color{white}3}} and {\small\colorbox{black}{\color{white}4}}) how the spectral radii of the iteration matrix $\mathbf{M}^{-1}_{\gamma}\mathbf{N}_{k}$ (computed from \eqref{EqP1Gen}) of the stabilized method \eqref{EqP1} is less than one and is decreasing with respect to the stabilization parameter $\gamma \in [10^{0}, \, 10^{15}]$, for the different values of the fixed stepsizes $\alpha_{k} = \alpha = 10^{i}, i = -5, -2, 0, 2, 5$. The stabilized iterative method defined by \eqref{EqP1} thus converges for any $\gamma > 10^{1}$ (see \cite{AllaireKaber08} - Theorem $8.1.1$). This convergence is faster than the stabilization parameter is large.
\begin{table}[h]
\begin{center} 
\begin{tikzpicture}[scale=0.87]
\draw[gray] [white, fill=gray!15.5]
    (0.0,6).. controls (15.5,6) ..(15.5,6)
    .. controls (15.5, 7.3) ..(15.5, 7.3)
    .. controls (0.0, 7.3) ..(0.0, 7.3)
    .. controls (0.0,6) ..(0.0,6);

\draw[gray] [white, fill=gray!10]
    (0.0,4.0).. controls (15.5,4.0) ..(15.5,4.0)
    .. controls (15.5, 5) .. (15.5, 5)
    .. controls (0.0, 5)  .. (0.0, 5)
    .. controls (0.0,4.0) .. (0.0,4.0);
\draw[gray] [white, fill=gray!10]
    (0.0,2.0).. controls (15.5,2.0) ..(15.5,2.0)
    .. controls (15.5, 3) .. (15.5, 3)
    .. controls (0.0, 3)  .. (0.0, 3)
    .. controls (0.0,2.0) .. (0.0,2.0);
\draw[gray] [white, fill=gray!10]
    (0.0,0.0).. controls (15.5,0.0) ..(15.5,0.0)
    .. controls (15.5, 1) .. (15.5, 1)
    .. controls (0.0, 1)  .. (0.0, 1)
    .. controls (0.0,0.0) .. (0.0,0.0);

\draw[, very thin] (1.5,0) -- (1.5,7.3);
\draw[, very thin] (4.5,0) -- (4.5,7.3);
\draw[, very thin] (7.5,0) -- (7.5,7.3);
\draw[, very thin] (10.5,0) -- (10.5,7.3);
\draw[, very thin] (13.5,0) -- (13.5,7.3);

\draw[, very thin] (0,0) -- (15.5,0);
\draw[, very thin] (0,1) -- (15.5,1);
\draw[, very thin] (0,2)   -- (15.5,2);
\draw[, very thin] (0,3)   -- (15.5,3);
\draw[, very thin] (0,4)   -- (15.5,4);
\draw[, very thin] (0,5) -- (15.5,5);

\draw[black](0.56, 7.)  node [rotate=40, scale=0.9] {\text{Stab.}};
\draw[black](0.75, 6.65)  node [rotate=40, scale=0.86] {\text{Parameter}};
\draw[black](1., 6.3)  node [rotate=40, scale=1] {$(\gamma)$};
\draw[black](0.75, 5.5)  node [rotate=0, scale=1] {$10^{3}$};
\draw[black](0.75, 4.5)  node [rotate=0, scale=1] {$10^{4}$};
\draw[black](0.75, 3.5)  node [rotate=0, scale=1] {$10^{5}$};
\draw[black](0.75, 2.5)  node [rotate=0, scale=1] {$10^{6}$};
\draw[black](0.75, 1.5)  node [rotate=0, scale=1] {$10^{10}$};
\draw[black](0.75, 0.5)  node [rotate=0, scale=1] {$10^{12}$};

\draw[black](3, 6.7)   node [rotate=0, scale=1] {$\|\mathbf{A}_{1}\bm{x}^{[k]}_{\gamma} - \bm{b}_{1}\|_{n}$};
\draw[black](3, 5.5)  node [rotate=0, scale=1] {$3.501759e^{-07}$};
\draw[black](3, 4.5)  node [rotate=0, scale=1] {$2.883407e^{-08}$};
\draw[black](3, 3.5)  node [rotate=0, scale=1] {$1.444926e^{-06}$};
\draw[black](3, 2.5)  node [rotate=0, scale=1] {$1.448875e^{-08}$};
\draw[black](3, 1.5)  node [rotate=0, scale=1] {$9.995320e^{-09}$};
\draw[black](3, 0.5)  node [rotate=0, scale=1] {$9.987286e^{-11}$};

\draw[black](6, 6.6)   node [rotate=0, scale=1] {$\frac{\|\mathbf{A}_{1}\bm{x}^{[k]}_{\gamma} - \bm{b}_{1}\|_{n}}{\|\mathbf{A}_{1}\bm{x}^{[0]}_{\gamma} - \bm{b}_{1}\|_{n}}$};
\draw[black](6, 5.5)  node [rotate=0, scale=1] {$1.282085e^{-08}$};
\draw[black](6, 4.5)  node [rotate=0, scale=1] {$1.055690e^{-09}$};
\draw[black](6, 3.5)  node [rotate=0, scale=1] {$5.290250e^{-08}$};
\draw[black](6, 2.5)  node [rotate=0, scale=1] {$5.304710e^{-10}$};
\draw[black](6, 1.5)  node [rotate=0, scale=1] {$3.659547e^{-10}$};
\draw[black](6, 0.5)  node [rotate=0, scale=1] {$3.656605e^{-12}$};

\draw[black](9., 6.7)   node [rotate=0, scale=1] {$\|\bm{x}^{\star} - \bm{x}^{[k]}_{\gamma}\|_{n}$};
\draw[black](9, 5.5)  node [rotate=0, scale=1] {$4.311856e^{-06}$};
\draw[black](9, 4.5)  node [rotate=0, scale=1] {$3.550480e^{-07}$};
\draw[black](9, 3.5)  node [rotate=0, scale=1] {$1.779234e^{-05}$};
\draw[black](9, 2.5)  node [rotate=0, scale=1] {$1.784097e^{-07}$};
\draw[black](9, 1.5)  node [rotate=0, scale=1] {$1.223421e^{-07}$};
\draw[black](9, 0.5)  node [rotate=0, scale=1] {$1.222441e^{-09}$};

\draw[black](12., 6.6)   node [rotate=0, scale=1] {$\frac{\|\bm{x}^{\star} - \bm{x}^{[k]}_{\gamma}\|_{n}}{\|\bm{x}^{\star}\|_{n}}$};
\draw[black](12, 5.5)  node [rotate=0, scale=1] {$2.155928e^{-06}$};
\draw[black](12, 4.5)  node [rotate=0, scale=1] {$1.775240e^{-07}$};
\draw[black](12, 3.5)  node [rotate=0, scale=1] {$8.896171e^{-06}$};
\draw[black](12, 2.5)  node [rotate=0, scale=1] {$8.920487e^{-08}$};
\draw[black](12, 1.5)  node [rotate=0, scale=1] {$6.117104e^{-08}$};
\draw[black](12, 0.5)  node [rotate=0, scale=1] {$6.112207e^{-10}$};

\draw[black](14.5, 6.85)   node [rotate=0, scale=1] {\text{$\#$ Iters.}};
\draw[black](14.5, 6.3)   node [rotate=0, scale=1] {$[k]$};
\draw[black](14.5, 5.5)  node [rotate=0, scale=1] {$7$};
\draw[black](14.5, 4.5)  node [rotate=0, scale=1] {$4$};
\draw[black](14.5, 3.5)  node [rotate=0, scale=1] {$2$};
\draw[black](14.5, 2.5)  node [rotate=0, scale=1] {$2$};
\draw[black](14.5, 1.5)  node [rotate=0, scale=1] {$1$};
\draw[black](14.5, 0.5)  node [rotate=0, scale=1] {$1$};

\end{tikzpicture} 
\caption{Numerical results obtained from the stabilized method \eqref{EqP1} with the matrix $\mathbf{A}_{1}$ in \eqref{Mats} and the right-hand vector computed as $\bm{b}_{1} = \mathbf{A}_{1}\bm{x}^{\star}$, where $\bm{x}^{\star} = \bm{1}^{t}$ ($\textbf{cond}(\mathbf{A}_{1}) = 171.62$, $\textbf{rank}(\mathbf{A}_{1}) = 4$). These computations are made with respect to the different values of the stabilization parameter $\gamma = 10^{i}$, $i = 3$, $4$, $5$, $6$, $10$, $12$, under the stopping criterion \eqref{CTR2}, the initial value $\bm{x}^{[0]}_{\gamma} = \bm{0}$, and the parameters $\alpha = 1$, $k_{\text{max}} = 100$ and $\varepsilon = 10^{-5}$.}
\label{TabWellC1}
\end{center}
\begin{center} 
\begin{tikzpicture}[scale=0.87]
\draw[gray] [white, fill=gray!15.5]
    (0.0,6).. controls (15.5,6) ..(15.5,6)
    .. controls (15.5, 7.3) ..(15.5, 7.3)
    .. controls (0.0, 7.3) ..(0.0, 7.3)
    .. controls (0.0,6) ..(0.0,6);

\draw[gray] [white, fill=gray!10]
    (0.0,4.0).. controls (15.5,4.0) ..(15.5,4.0)
    .. controls (15.5, 5) .. (15.5, 5)
    .. controls (0.0, 5)  .. (0.0, 5)
    .. controls (0.0,4.0) .. (0.0,4.0);
\draw[gray] [white, fill=gray!10]
    (0.0,2.0).. controls (15.5,2.0) ..(15.5,2.0)
    .. controls (15.5, 3) .. (15.5, 3)
    .. controls (0.0, 3)  .. (0.0, 3)
    .. controls (0.0,2.0) .. (0.0,2.0);
\draw[gray] [white, fill=gray!10]
    (0.0,0.0).. controls (15.5,0.0) ..(15.5,0.0)
    .. controls (15.5, 1) .. (15.5, 1)
    .. controls (0.0, 1)  .. (0.0, 1)
    .. controls (0.0,0.0) .. (0.0,0.0);

\draw[, very thin] (1.5,0) -- (1.5,7.3);
\draw[, very thin] (4.5,0) -- (4.5,7.3);
\draw[, very thin] (7.5,0) -- (7.5,7.3);
\draw[, very thin] (10.5,0) -- (10.5,7.3);
\draw[, very thin] (13.5,0) -- (13.5,7.3);

\draw[, very thin] (0,0) -- (15.5,0);
\draw[, very thin] (0,1) -- (15.5,1);
\draw[, very thin] (0,2)   -- (15.5,2);
\draw[, very thin] (0,3)   -- (15.5,3);
\draw[, very thin] (0,4)   -- (15.5,4);
\draw[, very thin] (0,5) -- (15.5,5);

\draw[black](0.56, 7.)  node [rotate=40, scale=0.9] {\text{Stab.}};
\draw[black](0.75, 6.65)  node [rotate=40, scale=0.86] {\text{Parameter}};
\draw[black](1., 6.3)  node [rotate=40, scale=1] {$(\gamma)$};
\draw[black](0.75, 5.5)  node [rotate=0, scale=1] {$10^{3}$};
\draw[black](0.75, 4.5)  node [rotate=0, scale=1] {$10^{4}$};
\draw[black](0.75, 3.5)  node [rotate=0, scale=1] {$10^{5}$};
\draw[black](0.75, 2.5)  node [rotate=0, scale=1] {$10^{6}$};
\draw[black](0.75, 1.5)  node [rotate=0, scale=1] {$10^{10}$};
\draw[black](0.75, 0.5)  node [rotate=0, scale=1] {$10^{12}$};

\draw[black](3, 6.7)   node [rotate=0, scale=1] {$\|\mathbf{A}_{2}\bm{x}^{[k]}_{\gamma} - \bm{b}_{2}\|_{n}$};
\draw[black](3, 5.5)  node [rotate=0, scale=1] {$4.580501e^{-07}$};
\draw[black](3, 4.5)  node [rotate=0, scale=1] {$1.730743e^{-08}$};
\draw[black](3, 3.5)  node [rotate=0, scale=1] {$5.920622e^{-08}$};
\draw[black](3, 2.5)  node [rotate=0, scale=1] {$5.923764e^{-10}$};
\draw[black](3, 1.5)  node [rotate=0, scale=1] {$2.028547e^{-09}$};
\draw[black](3, 0.5)  node [rotate=0, scale=1] {$2.028403e^{-11}$};

\draw[black](6, 6.6)   node [rotate=0, scale=1] {$\frac{\|\mathbf{A}_{2}\bm{x}^{[k]}_{\gamma} - \bm{b}_{2}\|_{n}}{\|\mathbf{A}_{2}\bm{x}^{[0]}_{\gamma} - \bm{b}_{2}\|_{n}}$};
\draw[black](6, 5.5)  node [rotate=0, scale=1] {$4.939281e^{-08}$};
\draw[black](6, 4.5)  node [rotate=0, scale=1] {$1.866308e^{-09}$};
\draw[black](6, 3.5)  node [rotate=0, scale=1] {$6.384370e^{-09}$};
\draw[black](6, 2.5)  node [rotate=0, scale=1] {$6.387759e^{-11}$};
\draw[black](6, 1.5)  node [rotate=0, scale=1] {$2.187438e^{-10}$};
\draw[black](6, 0.5)  node [rotate=0, scale=1] {$2.187283e^{-12}$};

\draw[black](9., 6.7)   node [rotate=0, scale=1] {$\|\bm{x}^{\star} - \bm{x}^{[k]}_{\gamma}\|_{n}$};
\draw[black](9, 5.5)  node [rotate=0, scale=1] {$2.487421e^{-06}$};
\draw[black](9, 4.5)  node [rotate=0, scale=1] {$9.398736e^{-08}$};
\draw[black](9, 3.5)  node [rotate=0, scale=1] {$3.215171e^{-07}$};
\draw[black](9, 2.5)  node [rotate=0, scale=1] {$3.216879e^{-09}$};
\draw[black](9, 1.5)  node [rotate=0, scale=1] {$1.099173e^{-08}$};
\draw[black](9, 0.5)  node [rotate=0, scale=1] {$1.099093e^{-10}$};

\draw[black](12., 6.6)   node [rotate=0, scale=1] {$\frac{\|\bm{x}^{\star} - \bm{x}^{[k]}_{\gamma}\|_{n}}{\|\bm{x}^{\star}\|_{n}}$};
\draw[black](12, 5.5)  node [rotate=0, scale=1] {$1.243710e^{-06}$};
\draw[black](12, 4.5)  node [rotate=0, scale=1] {$4.699368e^{-08}$};
\draw[black](12, 3.5)  node [rotate=0, scale=1] {$1.607585e^{-07}$};
\draw[black](12, 2.5)  node [rotate=0, scale=1] {$1.608440e^{-09}$};
\draw[black](12, 1.5)  node [rotate=0, scale=1] {$5.495865e^{-09}$};
\draw[black](12, 0.5)  node [rotate=0, scale=1] {$5.495463e^{-11}$};

\draw[black](14.5, 6.85)   node [rotate=0, scale=1] {\text{$\#$ Iters.}};
\draw[black](14.5, 6.3)   node [rotate=0, scale=1] {$[k]$};
\draw[black](14.5, 5.5)  node [rotate=0, scale=1] {$4$};
\draw[black](14.5, 4.5)  node [rotate=0, scale=1] {$3$};
\draw[black](14.5, 3.5)  node [rotate=0, scale=1] {$2$};
\draw[black](14.5, 2.5)  node [rotate=0, scale=1] {$2$};
\draw[black](14.5, 1.5)  node [rotate=0, scale=1] {$1$};
\draw[black](14.5, 0.5)  node [rotate=0, scale=1] {$1$};

\end{tikzpicture} 
\caption{Numerical results obtained from the stabilized method \eqref{EqP1} with the matrix $\mathbf{A}_{2}$ in \eqref{Mats} and the right-hand vector computed as $\bm{b}_{2} = \mathbf{A}_{2}\bm{x}^{\star}$, where $\bm{x}^{\star} = \bm{1}^{t}$ ($\textbf{cond}(\mathbf{A}_{2}) = 35.37$, $\textbf{rank}(\mathbf{A}_{2}) = 4$). The computations are made with respect to the different values of the stabilization parameter $\gamma = 10^{i}$, $i = 3$, $4$, $5$, $6$, $10$, $12$, under the stopping criterion \eqref{CTR2}, the initial value $\bm{x}^{[0]}_{\gamma} = \bm{0}$, and the parameters $\alpha = 1$, $k_{\text{max}} = 100$ and $\varepsilon = 10^{-5}$.}
\label{TabWellC2}
\end{center}
\end{table}

In the solving of linear systems $\mathbf{A}_{1}\bm{x}^{\star} = \bm{b}_{1}$ and $\mathbf{A}_{2}\bm{x}^{\star} = \bm{b}_{2}$, the right-hand sides $\bm{b}_{1}$ and $\bm{b}_{2}$ are computed from the exact solution $\bm{x}^{\star} = (1, 1, 1, 1)^{t}$. Instead the normalized residual in \eqref{CTR1} as a stopping criterion of the iterations, we use a criterion based on the relative error
\begin{align}\label{CTR2}
      \frac{\|\bm{x}^{\star} - \bm{x}^{[k]}_{\gamma}\|_{n}}{\|\bm{x}^{\star}\|_{n}} < \varepsilon, \; \forall \; k = 1, 2, 3, \cdots.
\end{align}
From the stabilized gradient method \eqref{EqP1} with the initial vector $\bm{x}^{[0]}_{\gamma} = \bm{0}$ and the tolerance $\varepsilon = 10^{-5}$, we compute approximations of the optimal solution $\bm{x}^{\star}$ and present results in Table \ref{TabWellC1} and Table \ref{TabWellC2}. We note that the stabilized method not only converges but behaves like a direct method by making one iteration, whenever the stabilization parameter becomes larger (most often when $\gamma > \frac{1}{\varepsilon}$). The stabilized method provides very accurate approximations.

\subsection{\bf Increasingly ill-conditioned linear systems}\label{NEx2}
In this second part, we take back a sequence of examples analyzed in \cite{Beck14} - Chap. $4$, namely {\bf Example $4.6$} (exact line search), {\bf Example $4.8$} (constant stepsize) and {\bf Example $4.9$} (stepsize selection by backtracking), where the gradient method is applied to the problem
\begin{align}\label{Function1}
    \min_{(x, y) \in \mathbb{R}^2}\left(x^2 + ay^2\right).
\end{align}
This problem is often used to give a more formalized analysis of the zigzagging phenomenon and the empirically observed slow convergence rate of the steepest descent algorithm (see \cite{Beck14} - Chap. $4$, \cite{GrivaNashSofer09}  - Chap. $12$, and \cite{BazaraaSheraliShetty06}  - Chap. $8.6$). Problem \eqref{Function1} can be formulated under the quadratic problem 
\begin{align}\label{Function2}
    \min_{\bm{x} \in \mathbb{R}^2}\left<\mathbb{A}_{a}\bm{x}, \,\bm{x}\right>_{2} + 2\left<\bm{b}, \, \bm{x}\right>_{2},
\end{align}
\clearpage
\begin{figure}[h]
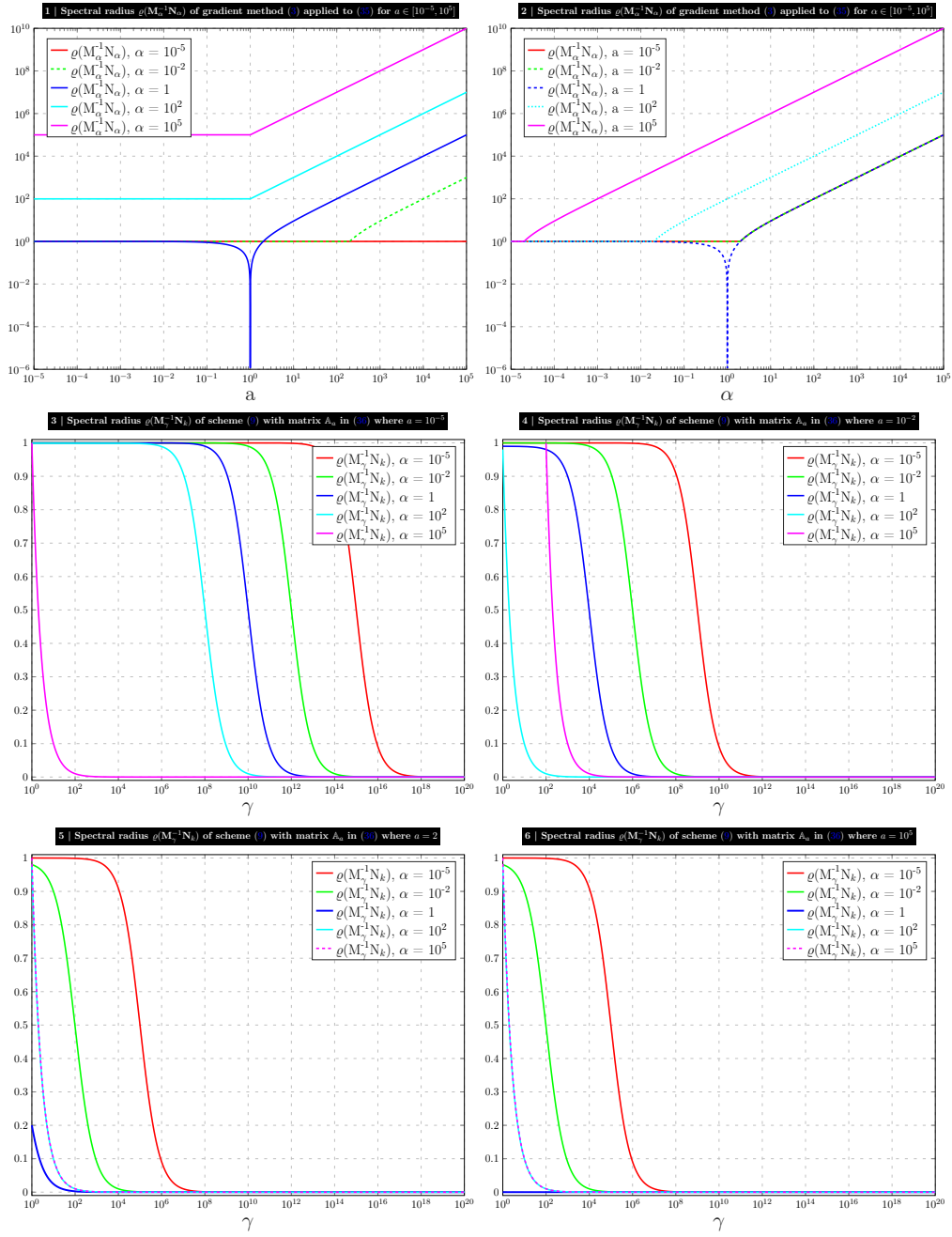
  
    \centering  
    \scalebox{0.4}{
        \definecolor{mycolor1}{rgb}{0.00000,1.00000,1.00000}%
\definecolor{mycolor2}{rgb}{1.00000,0.00000,1.00000}%
\begin{tikzpicture}

\begin{axis}[%
width=6.028in,
height=4.754in,
at={(1.011in,0.642in)},
scale only axis,
xmode=log,
xmin=1e-05,
xmax=100000,
xlabel style={font=\color{white!15!black}, scale=1.9},
xlabel={\text{a}},
ymode=log,
ymin=1.0e-06,
ymax=10000000000,
axis background/.style={fill=white},
title style={font=\bfseries, scale=0.9},
title={\normalsize\colorbox{black}{\color{white}1 | Spectral radius $\varrho(\mathbf{M}^{-1}_{\alpha}\mathbf{N}_{\alpha})$ of gradient method \eqref{LIM1} applied to \eqref{Function2} for $a \in [10^{-5}, 10^{5}]$}},
xmajorgrids,
ymajorgrids,
grid style={dashed},
legend style={legend cell align=left, align=left, draw=white!15!black, font=\Large},
ytick={0.000001, 0.0001, 0.01, 1, 100, 10000, 1000000, 100000000, 10000000000},
grid style = loosely dashed,
legend pos=north west,
]
\addplot [color=red, line width=1.5pt]
  table[row sep=crcr]{%
1e-05	0.9999999999\\
100000	0.99999\\
};
\addlegendentry{$\varrho\text{(M}^{\text{-1}}_{\alpha}\text{N}_{\alpha}\text{), }\alpha\text{ = 10}^{\text{-5}}$}

\addplot [color=green, dashed, line width=1.5pt]
  table[row sep=crcr]{%
1e-05 0.9999999\\
0.492313187422511 0.995076868125775\\
3.73407947600289  0.989999999999999\\
199.489356765898  0.994893567658976\\
223.038448216881  1.23038448216881\\
252.806793511704  1.52806793511704\\
291.117790490581  1.91117790490581\\
341.496890320409  2.41496890320408\\
409.497977587446  3.09497977587446\\
504.27422460669 4.0427422460669\\
641.761678824844  5.41761678824844\\
851.694402579699  7.51694402579698\\
1194.7191024838 10.947191024838\\
1810.24695559764  17.1024695559764\\
3077.99456635026  29.7799456635026\\
6341.33767965106  62.4133767965106\\
19179.7052302084  190.797052302084\\
100000  999.000000000001\\
};
\addlegendentry{$\varrho\text{(M}^{\text{-1}}_{\alpha}\text{N}_{\alpha}\text{), }\alpha\text{ = 10}^{\text{-2}}$}

\addplot [color=blue, line width=1.5pt]
  table[row sep=crcr]{%
1e-05	0.99999\\
0.0049938819265627	0.995006118073437\\
0.0155258476087133	0.984474152391287\\
0.0310463357040413	0.968953664295959\\
0.0511708193414097	0.94882918065859\\
0.0754920793619511	0.924507920638049\\
0.103568193605116	0.896431806394884\\
0.134932452172618	0.865067547827382\\
0.169108494342676	0.830891505657325\\
0.205626664730531	0.794373335269469\\
0.244040237370828	0.755959762629172\\
0.283936644296703	0.716063355703297\\
0.324946939682887	0.675053060317113\\
0.366748593044273	0.633251406955727\\
0.409067956421804	0.590932043578195\\
0.451671268685387	0.548328731314613\\
0.494360166009515	0.505639833990486\\
0.536955719837899	0.463044280162101\\
0.579286591637851	0.420713408362149\\
0.621174300096666	0.378825699903334\\
0.662418544547808	0.337581455452192\\
0.702791474847326	0.297208525152674\\
0.742024913624524	0.257975086375476\\
0.779814828464624	0.220185171535375\\
0.815820269050694	0.184179730949306\\
0.849666494843715	0.150333505156285\\
0.880962646736207	0.119037353263794\\
0.90930901314553	0.0906909868544703\\
0.934315378090417	0.0656846219095826\\
0.955620588791235	0.0443794112087663\\
0.972918583442071	0.027081416557929\\
0.985971812761792	0.0140281872382077\\
0.994659648479957	0.00534035152004164\\
0.999069037851145	0.000930962148854463\\
0.99999654612798	3.45387202027592e-06\\
0.999998848708	1.15129199884656e-06\\
1.00498018795078	0.00498018795078247\\
1.01574221090801	0.0157422109080092\\
1.03215875660332	0.0321587566033185\\
1.05449251725203	0.0544925172520301\\
1.0832021253505	0.0832021253504968\\
1.1189366648977	0.1189366648977\\
1.16256570275812	0.162565702758124\\
1.21523846781968	0.215238467819679\\
1.27845224657647	0.278452246576468\\
1.35416165872057	0.35416165872057\\
1.44493568228439	0.444935682284392\\
1.5541699913794	0.554169991379394\\
1.68641137421465	0.686411374214651\\
1.84783643729866	0.847836437298663\\
2.04696553808427	1.04696553808427\\
2.29579198890161	1.29579198890161\\
2.61160530239845	1.61160530239845\\
3.02001116118474	2.02001116118474\\
3.56021770120712	2.56021770120711\\
4.29459845267924	3.29459845267924\\
5.32721926756573	4.32721926756573\\
6.84203960118609	5.84203960118609\\
9.18898613700222	8.18898613700222\\
13.1002810117585	12.1002810117585\\
20.3196222881222	19.3196222881222\\
35.8510314848901	34.8510314848901\\
78.9953794895156	77.9953794895156\\
279.667335451587	278.667335451587\\
5063.54729120761	5062.54729120761\\
100000	99999\\
};
\addlegendentry{$\varrho\text{(M}^{\text{-1}}_{\alpha}\text{N}_{\alpha}\text{), }\alpha\text{ = 1}$}

\addplot [color=mycolor1, line width=1.5pt]
  table[row sep=crcr]{%
1e-05	98.9999999999999\\
1.00489457165056	99.4894571650559\\
4.35578629296372	434.578629296373\\
243.583989191972	24357.3989191971\\
100000	9999999\\
};
\addlegendentry{$\varrho\text{(M}^{\text{-1}}_{\alpha}\text{N}_{\alpha}\text{), }\alpha\text{ = 10}^{\text{2}}$}

\addplot [color=mycolor2, line width=1.5pt]
  table[row sep=crcr]{%
1e-05	99999\\
1.00494547776074	100493.547776075\\
100000	9999999998.99998\\
};
\addlegendentry{$\varrho\text{(M}^{\text{-1}}_{\alpha}\text{N}_{\alpha}\text{), }\alpha\text{ = 10}^{\text{5}}$}

\end{axis}
\end{tikzpicture}%
     %
        \definecolor{mycolor1}{rgb}{0.00000,1.00000,1.00000}%
\definecolor{mycolor2}{rgb}{1.00000,0.00000,1.00000}%
\begin{tikzpicture}

\begin{axis}[%
width=6.028in,
height=4.754in,
at={(1.011in,0.642in)},
scale only axis,
xmode=log,
xmin=1e-05,
xmax=100000,
xlabel style={font=\color{white!15!black}, scale=1.9},
xlabel={$\alpha$},
ymode=log,
ymin=1.0e-06,
ymax=10000000000,
axis background/.style={fill=white},
title style={font=\bfseries, scale=0.9},
title={\normalsize\colorbox{black}{\color{white}2 | Spectral radius $\varrho(\mathbf{M}^{-1}_{\alpha}\mathbf{N}_{\alpha})$ of gradient method \eqref{LIM1} applied to \eqref{Function2} for $\alpha \in [10^{-5}, 10^{5}]$}},
xmajorgrids,
ymajorgrids,
grid style={dashed},
legend style={legend cell align=left, align=left, draw=white!15!black, font=\Large},
ytick={0.000001, 0.0001, 0.01, 1, 100, 10000, 1000000, 100000000, 10000000000},
grid style = loosely dashed,
legend pos=north west,
]
\addplot [color=red, line width=1.5pt]
  table[row sep=crcr]{%
1e-05 0.9999999999\\
2.00494530807401  1.00494530807401\\
2.24299102118876  1.24299102118876\\
2.54414313046339  1.54414313046339\\
2.93207818475396  1.93207818475396\\
3.44279016533574  2.44279016533574\\
4.13309665485049  3.13309665485049\\
5.09684669298946  4.09684669298946\\
6.49792330824323  5.49792330824322\\
8.64326032546605  7.64326032546604\\
12.1619050505313  11.1619050505313\\
18.5097531411352  17.5097531411352\\
31.69148939665  30.6914893966501\\
66.1050566927257  65.1050566927257\\
205.618757884451  204.618757884451\\
2127.64336200666  2126.64336200666\\
100000  99999\\
};
\addlegendentry{$\varrho\text{(M}^{\text{-1}}_{\alpha}\text{N}_{\alpha}\text{), }\text{a = 10}^{\text{-5}}$}

\addplot [color=green, dashed, line width=1.5pt]
  table[row sep=crcr]{%
1e-05 0.9999999\\
0.492313187422511 0.995076868125775\\
1.54084700580941  0.984591529941906\\
1.98553034949825  0.985530349498244\\
2.2186486640251 1.2186486640251\\
2.51312181486273  1.51312181486273\\
2.89176199772529  1.89176199772529\\
3.38915594929068  2.38915594929068\\
4.05967774464774  3.05967774464774\\
4.99273562387356  3.99273562387356\\
6.34356771630639  5.34356771630639\\
8.40081831760449  7.40081831760449\\
11.7505882093452  10.7505882093452\\
17.7317918628665  16.7317918628665\\
29.9581263188651  28.9581263188651\\
61.0270923730858  60.0270923730858\\
179.991387791602  178.991387791602\\
1523.5644070719 1522.5644070719\\
100000  99999\\
};
\addlegendentry{$\varrho\text{(M}^{\text{-1}}_{\alpha}\text{N}_{\alpha}\text{), }\text{a = 10}^{\text{-2}}$}

\addplot [color=blue, dashed, line width=1.5pt]
  table[row sep=crcr]{%
1e-05 0.99999\\
0.0049938819265627  0.995006118073437\\
0.0155258476087133  0.984474152391287\\
0.0310463357040413  0.968953664295959\\
0.0511708193414097  0.94882918065859\\
0.0754920793619511  0.924507920638049\\
0.103568193605116 0.896431806394884\\
0.134932452172618 0.865067547827382\\
0.169108494342676 0.830891505657325\\
0.205626664730531 0.794373335269469\\
0.244040237370828 0.755959762629172\\
0.283936644296703 0.716063355703297\\
0.324946939682887 0.675053060317113\\
0.366748593044273 0.633251406955727\\
0.409067956421804 0.590932043578197\\
0.451671268685387 0.548328731314613\\
0.494360166009515 0.505639833990486\\
0.536955719837899 0.463044280162101\\
0.579286591637851 0.420713408362149\\
0.621174300096666 0.378825699903334\\
0.662418544547808 0.337581455452192\\
0.702791474847326 0.297208525152674\\
0.742024913624524 0.257975086375476\\
0.779814828464624 0.220185171535375\\
0.815820269050694 0.184179730949306\\
0.849666494843715 0.150333505156285\\
0.880962646736207 0.119037353263794\\
0.90930901314553  0.0906909868544703\\
0.934315378090417 0.0656846219095825\\
0.955620588791235 0.0443794112087663\\
0.972918583442071 0.0270814165579288\\
0.985971812761792 0.0140281872382077\\
0.994659648479957 0.00534035152004161\\
0.999069037851145 0.000930962148854382\\
0.99999654612798  3.45387202023816e-06\\
0.999998848708  1.15129199875601e-06\\
1.00498018795078  0.00498018795078242\\
1.01574221090801  0.0157422109080092\\
1.03215875660332  0.0321587566033185\\
1.05449251725203  0.0544925172520301\\
1.0832021253505 0.0832021253504968\\
1.1189366648977 0.1189366648977\\
1.16256570275812  0.162565702758123\\
1.21523846781968  0.215238467819679\\
1.27845224657647  0.278452246576468\\
1.35416165872057  0.35416165872057\\
1.44493568228439  0.444935682284392\\
1.5541699913794 0.554169991379394\\
1.68641137421465  0.686411374214651\\
1.84783643729866  0.847836437298663\\
2.04696553808427  1.04696553808427\\
2.29579198890161  1.29579198890161\\
2.61160530239845  1.61160530239845\\
3.02001116118474  2.02001116118474\\
3.56021770120712  2.56021770120711\\
4.29459845267924  3.29459845267924\\
5.32721926756573  4.32721926756573\\
6.84203960118609  5.84203960118609\\
9.18898613700222  8.18898613700222\\
13.1002810117585  12.1002810117585\\
20.3196222881222  19.3196222881222\\
35.8510314848901  34.8510314848901\\
78.9953794895156  77.9953794895156\\
279.667335451587  278.667335451587\\
5063.54729120761  5062.54729120761\\
100000  99999\\
};
\addlegendentry{$\varrho\text{(M}^{\text{-1}}_{\alpha}\text{N}_{\alpha}\text{), }\text{a = 1}$}

\addplot [color=mycolor1, dotted, line width=1.5pt]
  table[row sep=crcr]{%
1e-05 0.99999\\
0.0049938819265627  0.995006118073437\\
0.0155258476087133  0.984474152391287\\
0.0198550200421171  0.985502004211715\\
0.0221861188220947  1.21861188220947\\
0.0251307436455596  1.51307436455596\\
0.0289170073990977  1.89170073990978\\
0.0338906854762532  2.38906854762532\\
0.0405955435616492  3.05955435616492\\
0.0499256088469493  3.99256088469494\\
0.0634330188164884  5.34330188164884\\
0.0840040824327304  7.40040824327304\\
0.117498793443504 10.7498793443504\\
0.177304363930824 16.7304363930824\\
0.299551464817389 28.955146481739\\
0.610184931496559 60.0184931496559\\
1.79949450790488  178.949450790489\\
15.2270445465382  1521.70445465381\\
100000  9999999\\
};
\addlegendentry{$\varrho\text{(M}^{\text{-1}}_{\alpha}\text{N}_{\alpha}\text{), }\text{a = 10}^{\text{2}}$}

\addplot [color=mycolor2, line width=1.5pt]
  table[row sep=crcr]{%
1e-05 0.99999\\
2.00492915018737e-05  1.00492915018737\\
2.24297294488787e-05  1.24297294488788\\
2.54411676912012e-05  1.54411676912012\\
2.93204105250802e-05  1.93204105250801\\
3.44273863815365e-05  2.44273863815364\\
4.13302527940906e-05  3.13302527940906\\
5.09674693860834e-05  4.09674693860833\\
6.49778117058084e-05  5.49778117058084\\
8.64303145728598e-05  7.64303145728597\\
0.000121615270047821  11.1615270047821\\
0.000185090499197762  17.5090499197762\\
0.000316899205297795  30.6899205297795\\
0.000661002621718045  65.1002621718045\\
0.00205593903013241 204.593903013241\\
0.0212697474261501  2125.97474261501\\
100000  9999999998.99998\\
};
\addlegendentry{$\varrho\text{(M}^{\text{-1}}_{\alpha}\text{N}_{\alpha}\text{), }\text{a = 10}^{\text{5}}$}

\end{axis}
\end{tikzpicture}%
}  \\
    \scalebox{0.4}{
        \input{RayonSpectral3.tex}     %
        \input{RayonSpectral4.tex}}  \\
    \scalebox{0.4}{
        \input{RayonSpectral5.tex}     
        \input{RayonSpectral6.tex}}
    \caption{\footnotesize The spectral radius $\varrho(\mathbf{M}^{-1}_{\alpha}\mathbf{N}_{\alpha})$ of the iteration matrix $\mathbf{M}^{-1}_{\alpha}\mathbf{N}_{\alpha}$ from the gradient method \eqref{LIM1} applied to the problem \eqref{Function2} is represented in {\tiny\colorbox{black}{\color{white}1}} and {\tiny\colorbox{black}{\color{white}2}}, with respect to the parameters $a \in [10^{-5}, 10^{5}]$ (for $\alpha = 10^{i}, i = -5, -2, 0, 2, 5$) and $\alpha \in [10^{-5}, 10^{5}]$ (for $a = 10^{i}, i = -5, -2, 0, 2, 5$), respectively. The spectral radius $\varrho(\mathbf{M}^{-1}_{\gamma}\mathbf{N}_{k})$ of the iteration matrix $\mathbf{M}^{-1}_{\gamma}\mathbf{N}_{k}$ from the stabilized gradient method \eqref{EqP1} applied to the problem \eqref{Function2} is represented in {\tiny\colorbox{black}{\color{white}3}}, {\tiny\colorbox{black}{\color{white}4}}, {\tiny\colorbox{black}{\color{white}5}} and {\tiny\colorbox{black}{\color{white}6}}, with respect to the parameter $\gamma$ (for $\alpha = 10^{i}, i = -5, -2, 0, 2, 5$) and fixed $a = 10^{-5}, 10^{-2}, 2, 10^{5}$.}
    \label{Fig461}
\end{figure}
where $\bm{x} = (x, y)^{t}$, $\mathbb{A}_{a} \in \mathbb{R}^{2\times 2}$ is a symmetric positive definite matrix and $\bm{b} \in \mathbb{R}^{2}$ is a vector
\begin{align}\label{MatV1}
    \mathbb{A}_{a} = \left(\begin{array}{ll}
1 & 0 \\
0 & a
\end{array}\right), \text{ and } \, \bm{b} = \left(\begin{array}{l}
0  \\
0 
\end{array}\right).
\end{align}

In order to have an insight on the convergence of the gradient method \eqref{LIM1} applied to the quadratic problem \eqref{Function2}, we represent in Figure \ref{Fig461} - {\small\colorbox{black}{\color{white}1}} the spectral radius $\varrho(\mathbf{M}^{-1}_{\alpha}\mathbf{N}_{\alpha})$ of the iterative matrix $\mathbf{M}^{-1}_{\alpha}\mathbf{N}_{\alpha}$, with respect to the parameter $a \in [10^{-5}, 10^{5}]$ and for each of the different stepsizes $\alpha = 10^{i}, i = -5, -2, 0, 2, 5$. We thus observe increasing spectral radii when the parameter $a \geq 1$ and $\alpha \ne 10^{-5}$. Larger is the parameter $a$, slower thus is the convergence of the gradient method when it does not diverge. We also represent in Figure \ref{Fig461} - {\small\colorbox{black}{\color{white}2}} the same spectral radius, but with respect to the stepsizes $\alpha \in [10^{-5}, 10^{5}]$ for each of the parameters $a = 10^{i}, i = -5, -2, 0, 2, 5$. We thus observe a more unstable gradient method with respect to the stepsize $\alpha$.

The speed of convergence of iterative methods is also related to the conditioning of the Hessian matrix $\mathbb{A}_{a}$ of the quadratic problem \eqref{Function2}. For a symmetric positive definite matrix, the conditioning is defined as the ratio of the largest to the smallest eigenvalues. We then obtain the conditioning
\begin{align*}
    \textbf{cond}(\mathbb{A}_{a}) = \frac{\lambda_{\text{max}}}{\lambda_{\text{min}}} = 
    \left\{ 
        \begin{array}{rcl} 
            a   & \mbox{if} & a \geq 1 \\ 
            1/a & \mbox{if} & 0 < a < 1 
        \end{array}\right..
\end{align*}

For insights on the convergence of the stabilized gradient method \eqref{EqP1}, we present in Figure \ref{Fig461} - {\small\colorbox{black}{\color{white}3}}, {\small\colorbox{black}{\color{white}4}}, {\small\colorbox{black}{\color{white}5}}, {\small\colorbox{black}{\color{white}6}} the spectral radius $\varrho(\mathbf{M}^{-1}_{\gamma}\mathbf{N}_{k})$ from the scheme \eqref{EqP1Gen}, with respect to the stabilization parameter $\gamma$, for the different values of the stepsizes $\alpha = 10^{i}, i = -5, -2, 0, 2, 5$, and for each of the parameters $a = 10^{-5}$, $a = 10^{-2}$, $a = 2$ and $a = 10^{5}$. We thus observe decreasing spectral radius when $\gamma$ increases. The stabilized iterative method \eqref{EqP1} thus converges from some $\gamma_{0} > 0$ (i.e. for all $\gamma > \gamma_{0}$), whatever the fixed stepsize $\alpha$ among these values (which is coherent with Corollary \ref{CorMSN12} where this convergence is acquired for any $\alpha$ whenever $\gamma$ is chosen large enough). The convergence is faster than the stabilization parameter $\gamma$ is large. \\
In the forthcoming subsub-sections \ref{beck46}, \ref{beck48} and \ref{beck49}, we go into detail on the efficiency and the accuracy of the stabilized gradient method \eqref{EqP1} by re-investigating examples already analyzed in the chapter $4$ in \cite{Beck14} and based on the quadratic problem \eqref{Function2}. Through these examples, the following stopping criterion of the iterations is used  
\begin{align}\label{SCI2}
    \|\mathbb{A}_{a}\bm{x}^{[k]}_{\gamma} + \bm{b}\|_{2} \leq \varepsilon, \; \forall \; k = 1, 2, 3, \cdots.
\end{align}

\subsubsection{\bf Example $4.6$ (exact line stepsize) from Chapter $4$ in \cite{Beck14}:}\label{beck46}
In order to solve problem \eqref{Function2} by using the stabilized gradient method \eqref{EqP1} and under the same context than {\bf Example $4.6$} in \cite{Beck14}, we fixe the tolerance parameter to $\varepsilon = 10^{-5}$, the initial vector to $\bm{x}^{[0]}_{\gamma} = (2, 1)^{t}$, the parameter to $a = 2$, and the exact line search defined at the $k^{th}$ iteration by
\begin{align}\label{Exa461}
   \alpha_{k} = \frac{\left\|\mathbb{A}_{2}\bm{x}^{[k]}_{\gamma} - \bm{b}\right\|^{2}_{2}}{2\left(\mathbb{A}_{2}\bm{x}^{[k]}_{\gamma} - \bm{b}\right)^{t}\mathbb{A}_{2}\left(\mathbb{A}_{2}\bm{x}^{[k]}_{\gamma} - \bm{b}\right)}.
\end{align}
We thus obtain results provided in the following Table \ref{Tab461}, where we present at the last column the stabilized solution $\bm{x}^{[k]}_{\gamma}$ which is pretty close to the optimal solution $\bm{x}^{\star} = (0,0)^{t}$.
\begin{table}[h]
\begin{center}
\begin{tikzpicture}[scale=0.87]
\draw[gray] [white, fill=gray!15.5]
    (0.0,6).. controls (15.5,6) ..(15.5,6)
    .. controls (15.5, 7.3) ..(15.5, 7.3)
    .. controls (0.0, 7.3) ..(0.0, 7.3)
    .. controls (0.0,6) ..(0.0,6);

\draw[gray] [white, fill=gray!10]
    (0.0,4.0).. controls (15.5,4.0) ..(15.5,4.0)
    .. controls (15.5, 5) .. (15.5, 5)
    .. controls (0.0, 5)  .. (0.0, 5)
    .. controls (0.0,4.0) .. (0.0,4.0);
\draw[gray] [white, fill=gray!10]
    (0.0,2.0).. controls (15.5,2.0) ..(15.5,2.0)
    .. controls (15.5, 3) .. (15.5, 3)
    .. controls (0.0, 3)  .. (0.0, 3)
    .. controls (0.0,2.0) .. (0.0,2.0);
\draw[gray] [white, fill=gray!10]
    (0.0,0.0).. controls (15.5,0.0) ..(15.5,0.0)
    .. controls (15.5, 1) .. (15.5, 1)
    .. controls (0.0, 1)  .. (0.0, 1)
    .. controls (0.0,0.0) .. (0.0,0.0);

\draw[, very thin] (1.5,0) -- (1.5,7.3);
\draw[, very thin] (4.5,0) -- (4.5,7.3);
\draw[, very thin] (7.5,0) -- (7.5,7.3);
\draw[, very thin] (10.5,0) -- (10.5,7.3);
\draw[, very thin] (12.,0) -- (12.,7.3);

\draw[, very thin] (0,0) -- (15.5,0);
\draw[, very thin] (0,1) -- (15.5,1);
\draw[, very thin] (0,2)   -- (15.5,2);
\draw[, very thin] (0,3)   -- (15.5,3);
\draw[, very thin] (0,4)   -- (15.5,4);
\draw[, very thin] (0,5) -- (15.5,5);

\draw[black](0.56, 7.)  node [rotate=40, scale=0.9] {\text{Stab.}};
\draw[black](0.7, 6.65)  node [rotate=40, scale=0.86] {\text{Parameter}};
\draw[black](1., 6.3)  node [rotate=40, scale=1] {$(\gamma)$};
\draw[black](0.75, 5.5)  node [rotate=0, scale=1] {$1$};
\draw[black](0.75, 4.5)  node [rotate=0, scale=1] {$10$};
\draw[black](0.75, 3.5)  node [rotate=0, scale=1] {$10^{2}$};
\draw[black](0.75, 2.5)  node [rotate=0, scale=1] {$10^{5}$};
\draw[black](0.75, 1.5)  node [rotate=0, scale=1] {$10^{7}$};
\draw[black](0.75, 0.5)  node [rotate=0, scale=1] {$10^{10}$};

\draw[black](3, 6.7)   node [rotate=0, scale=1] {$\|\mathbb{A}_{2}\bm{x}^{[k]}_{\gamma} - \bm{b}\|_{2}$};
\draw[black](3, 5.5)  node [rotate=0, scale=1] {$5.274815e^{-6}$};
\draw[black](3, 4.5)  node [rotate=0, scale=1] {$1.052846e^{-6}$};
\draw[black](3, 3.5)  node [rotate=0, scale=1] {$6.570074e^{-7}$};
\draw[black](3, 2.5)  node [rotate=0, scale=1] {$1.353509e^{-10}$};
\draw[black](3, 1.5)  node [rotate=0, scale=1] {$2.687419e^{-7}$};
\draw[black](3, 0.5)  node [rotate=0, scale=1] {$2.687419e^{-10}$};

\draw[black](6, 6.6)   node [rotate=0, scale=1] {$\frac{\|\mathbb{A}_{2}\bm{x}^{[k]}_{\gamma} - \bm{b}\|_{2}}{\|\mathbb{A}_{2}\bm{x}^{[0]}_{\gamma} + \bm{b}\|_{2}}$};
\draw[black](6, 5.5)  node [rotate=0, scale=1] {$9.324644e^{-7}$};
\draw[black](6, 4.5)  node [rotate=0, scale=1] {$1.861187e^{-7}$};
\draw[black](6, 3.5)  node [rotate=0, scale=1] {$1.161436e^{-7}$};
\draw[black](6, 2.5)  node [rotate=0, scale=1] {$2.392689e^{-11}$};
\draw[black](6, 1.5)  node [rotate=0, scale=1] {$4.750730e^{-8}$};
\draw[black](6, 0.5)  node [rotate=0, scale=1] {$4.750731e^{-11}$};

\draw[black](9., 6.7)   node [rotate=0, scale=1] {$\|\bm{x}^{\star} - \bm{x}^{[k]}_{\gamma}\|_{2}$};
\draw[black](9, 5.5)  node [rotate=0, scale=1] {$2.637408e^{-6}$};
\draw[black](9, 4.5)  node [rotate=0, scale=1] {$5.264231e^{-7}$};
\draw[black](9, 3.5)  node [rotate=0, scale=1] {$3.285037e^{-7}$};
\draw[black](9, 2.5)  node [rotate=0, scale=1] {$6.767544e^{-11}$};
\draw[black](9, 1.5)  node [rotate=0, scale=1] {$1.335935e^{-7}$};
\draw[black](9, 0.5)  node [rotate=0, scale=1] {$1.335935e^{-10}$};

\draw[black](11.25, 6.8)   node [rotate=0, scale=1] {\text{$\#$ Iters}};
\draw[black](11.25, 6.3)   node [rotate=0, scale=1] {$[k]$};
\draw[black](11.25, 5.5)  node [rotate=0, scale=1] {$10$};
\draw[black](11.25, 4.5)  node [rotate=0, scale=1] {$5$};
\draw[black](11.25, 3.5)  node [rotate=0, scale=1] {$3$};
\draw[black](11.25, 2.5)  node [rotate=0, scale=1] {$2$};
\draw[black](11.25, 1.5)  node [rotate=0, scale=1] {$1$};
\draw[black](11.25, 0.5)  node [rotate=0, scale=1] {$1$};

\draw[black](13.8, 6.7)    node [rotate=0, scale=1] {$\bm{x}^{[k]}_{\gamma}$};
\draw[black](13.8, 5.5)  node [rotate=0, scale=1] {$1e^{-5}\left(0.264, 0.000\right)$};
\draw[black](13.8, 4.5)  node [rotate=0, scale=1] {$1e^{-6}\left(0.526, 0.000\right)$};
\draw[black](13.8, 3.5)  node [rotate=0, scale=1] {$1e^{-6}\left(0.329, 0.000\right)$};
\draw[black](13.8, 2.5)  node [rotate=0, scale=0.95] {$1e^{-10}\left(0.677, 0.000\right)$};
\draw[black](13.8, 1.5)  node [rotate=0, scale=1] {$1e^{-6}\left(0.133, 0.008\right)$};
\draw[black](13.8, 0.5)  node [rotate=0, scale=1] {$1e^{-9}\left(0.133, 0.008\right)$};

\end{tikzpicture} 
\caption{Numerical results obtained from the stabilized gradient method \eqref{EqP1} with the data $\mathbb{A}_{a}$ and $\bm{b}$ given in \eqref{MatV1}, where $a = 2$ ($\textbf{cond}(\mathbb{A}_{2}) = 2$, $\textbf{rank}(\mathbb{A}_{2}) = 2$). The computations are made with respect to the different values of the stabilization parameter $\gamma = 10^{i}$, $i = 0$, $1$, $2$, $5$, $7$, $10$, the initial iteration fixed to $\bm{x}^{[0]}_{\gamma} = (2, 1)^{t}$, and the stepsize computed from \eqref{Exa461}. The maximum iterations is $k_{\text{max}} = 100$ and the stopping criterion \eqref{SCI2} is used with $\varepsilon = 10^{-5}$.}
\label{Tab461}
\begin{tikzpicture}[scale=0.87]
\draw[gray] [white, fill=gray!15.5]
    (0.0,6).. controls (15.5,6) ..(15.5,6)
    .. controls (15.5, 7.3) ..(15.5, 7.3)
    .. controls (0.0, 7.3) ..(0.0, 7.3)
    .. controls (0.0,6) ..(0.0,6);

\draw[gray] [white, fill=gray!10]
    (0.0,4.0).. controls (15.5,4.0) ..(15.5,4.0)
    .. controls (15.5, 5) .. (15.5, 5)
    .. controls (0.0, 5)  .. (0.0, 5)
    .. controls (0.0,4.0) .. (0.0,4.0);
\draw[gray] [white, fill=gray!10]
    (0.0,2.0).. controls (15.5,2.0) ..(15.5,2.0)
    .. controls (15.5, 3) .. (15.5, 3)
    .. controls (0.0, 3)  .. (0.0, 3)
    .. controls (0.0,2.0) .. (0.0,2.0);
\draw[gray] [white, fill=gray!10]
    (0.0,0.0).. controls (15.5,0.0) ..(15.5,0.0)
    .. controls (15.5, 1) .. (15.5, 1)
    .. controls (0.0, 1)  .. (0.0, 1)
    .. controls (0.0,0.0) .. (0.0,0.0);

\draw[, very thin] (1.5,0) -- (1.5,7.3);
\draw[, very thin] (4.5,0) -- (4.5,7.3);
\draw[, very thin] (7.5,0) -- (7.5,7.3);
\draw[, very thin] (10.5,0) -- (10.5,7.3);
\draw[, very thin] (12.,0) -- (12.,7.3);

\draw[, very thin] (0,0) -- (15.5,0);
\draw[, very thin] (0,1) -- (15.5,1);
\draw[, very thin] (0,2)   -- (15.5,2);
\draw[, very thin] (0,3)   -- (15.5,3);
\draw[, very thin] (0,4)   -- (15.5,4);
\draw[, very thin] (0,5) -- (15.5,5);

\draw[black](0.56, 7.)  node [rotate=40, scale=0.9] {\text{Stab.}};
\draw[black](0.7, 6.65)  node [rotate=40, scale=0.86] {\text{Parameter}};
\draw[black](1., 6.3)  node [rotate=40, scale=1] {$(\gamma)$};
\draw[black](0.75, 5.5)  node [rotate=0, scale=1] {$1$};
\draw[black](0.75, 4.5)  node [rotate=0, scale=1] {$10$};
\draw[black](0.75, 3.5)  node [rotate=0, scale=1] {$10^{2}$};
\draw[black](0.75, 2.5)  node [rotate=0, scale=1] {$10^{5}$};
\draw[black](0.75, 1.5)  node [rotate=0, scale=1] {$10^{7}$};
\draw[black](0.75, 0.5)  node [rotate=0, scale=1] {$10^{10}$};

\draw[black](3, 6.7)   node [rotate=0, scale=1] {$\|\mathbb{A}_{2}\bm{x}^{[k]}_{\gamma} - \bm{b}\|_{2}$};
\draw[black](3, 5.5)  node [rotate=0, scale=1] {$5.089472e^{-6}$};
\draw[black](3, 4.5)  node [rotate=0, scale=1] {$1.199938e^{-6}$};
\draw[black](3, 3.5)  node [rotate=0, scale=1] {$2.830419e^{-6}$};
\draw[black](3, 2.5)  node [rotate=0, scale=1] {$3.243883e^{-10}$};
\draw[black](3, 1.5)  node [rotate=0, scale=1] {$3.687817e^{-7}$};
\draw[black](3, 0.5)  node [rotate=0, scale=1] {$3.687818e^{-10}$};

\draw[black](6, 6.6)   node [rotate=0, scale=1] {$\frac{\|\mathbb{A}_{2}\bm{x}^{[k]}_{\gamma} - \bm{b}\|_{2}}{\|\mathbb{A}_{2}\bm{x}^{[0]}_{\gamma} + \bm{b}\|_{2}}$};
\draw[black](6, 5.5)  node [rotate=0, scale=1] {$8.997000e^{-7}$};
\draw[black](6, 4.5)  node [rotate=0, scale=1] {$2.121211e^{-7}$};
\draw[black](6, 3.5)  node [rotate=0, scale=1] {$5.003521e^{-7}$};
\draw[black](6, 2.5)  node [rotate=0, scale=1] {$5.734430e^{-11}$};
\draw[black](6, 1.5)  node [rotate=0, scale=1] {$6.519202e^{-8}$};
\draw[black](6, 0.5)  node [rotate=0, scale=1] {$6.519202e^{-11}$};

\draw[black](9., 6.7)   node [rotate=0, scale=1] {$\|\bm{x}^{\star} - \bm{x}^{[k]}_{\gamma}\|_{2}$};
\draw[black](9, 5.5)  node [rotate=0, scale=1] {$2.544736e^{-6}$};
\draw[black](9, 4.5)  node [rotate=0, scale=1] {$5.999692e^{-7}$};
\draw[black](9, 3.5)  node [rotate=0, scale=1] {$1.415143e^{-6}$};
\draw[black](9, 2.5)  node [rotate=0, scale=1] {$1.620461e^{-10}$};
\draw[black](9, 1.5)  node [rotate=0, scale=1] {$1.811077e^{-7}$};
\draw[black](9, 0.5)  node [rotate=0, scale=1] {$1.811077e^{-10}$};

\draw[black](11.25, 6.8)   node [rotate=0, scale=1] {\text{$\#$ Iters}};
\draw[black](11.25, 6.3)   node [rotate=0, scale=1] {$[k]$};
\draw[black](11.25, 5.5)  node [rotate=0, scale=1] {$17$};
\draw[black](11.25, 4.5)  node [rotate=0, scale=1] {$6$};
\draw[black](11.25, 3.5)  node [rotate=0, scale=1] {$3$};
\draw[black](11.25, 2.5)  node [rotate=0, scale=1] {$2$};
\draw[black](11.25, 1.5)  node [rotate=0, scale=1] {$1$};
\draw[black](11.25, 0.5)  node [rotate=0, scale=1] {$1$};

\draw[black](13.8, 6.7)    node [rotate=0, scale=1] {$\bm{x}^{[k]}_{\gamma}$};
\draw[black](13.8, 5.5)  node [rotate=0, scale=1] {$1e^{-5}\left(0.255, 0.000\right)$};
\draw[black](13.8, 4.5)  node [rotate=0, scale=1] {$1e^{-6}\left(0.600, 0.000\right)$};
\draw[black](13.8, 3.5)  node [rotate=0, scale=1] {$1e^{-5}\left(0.142, 0.001\right)$};
\draw[black](13.8, 2.5)  node [rotate=0, scale=1] {$1e^{-9}\left(0.162, 0.004\right)$};
\draw[black](13.8, 1.5)  node [rotate=0, scale=1] {$1e^{-6}\left(0.180, 0.020\right)$};
\draw[black](13.8, 0.5)  node [rotate=0, scale=1] {$1e^{-9}\left(0.180, 0.020\right)$};
\end{tikzpicture} 
\caption{Numerical results obtained from the stabilized gradient method \eqref{EqP1} with the data $\mathbb{A}_{a}$ and $\bm{b}$ given in \eqref{MatV1}, where $a = 2$ ($\textbf{cond}(\mathbb{A}_{2}) = 2$, $\textbf{rank}(\mathbb{A}_{2}) = 2$). The computations are made with respect to the different values of the stabilization parameter $\gamma = 10^{i}$, $i = 0$, $1$, $2$, $5$, $7$, $10$, the initial iteration fixed to $\bm{x}^{[0]}_{\gamma} = (2, 1)^{t}$ and the constant stepsize $\alpha = 0.1$. The maximum iterations is $k_{\text{max}} = 100$, and the stopping criterion \eqref{SCI2} is used with $\varepsilon = 10^{-5}$.}
\label{Tab481}
\end{center} 
\end{table}
We also note that the stabilized gradient method \eqref{EqP1} behaves like a direct method by making one iteration, when the stabilization parameter becomes large (most often when $\gamma > \frac{1}{\varepsilon}$).\\
When compared with the {\bf Example $4.6$} in \cite{Beck14} (Chapter $4$) where the gradient method stopped after $13$ iterations with the solution $\bm{x}^{[k]} = 1e^{-5}(0.1254, \,-0.0627)$, the stabilized gradient method we propose here seem to be very efficient and accurate.

\subsubsection{\bf Example $4.8$ (constant stepsize) from Chapter $4$ in \cite{Beck14}:}\label{beck48}
We revise here the optimization problem in \eqref{Function2}, but under the same context than {\bf Example} $4.8$ in \cite{Beck14}, where a constant stepsize $\alpha = 0.1$ is used. By executing the stabilized gradient algorithm \eqref{alg1} with the initial vector $\bm{x}^{[0]}_{\gamma} = (2, 1)^{t}$ and the tolerance $\varepsilon = 10^{-5}$, we obtain results provided in Table \ref{Tab481}, where the solution $\bm{x}^{[k]}_{\gamma}$ remains close to the optimal solution $\bm{x}^{\star} = (0,0)^{t}$. 

We note the same behavioural than in the previous example, that is to say the stabilized gradient method that behaves like a direct method when the stabilization parameter becomes large (most often when $\gamma > \frac{1}{\varepsilon}$). While the gradient method presented in {\bf Example $4.8$} (\cite{Beck14} Chapter $4$) makes a number of $58$ iterations before reaching the same tolerance $\varepsilon = 10^{-5}$ in \eqref{SCI2}.\\ 
Otherwise taking a large stepsize (as in {\bf Example} $4.8$ in \cite{Beck14}, where the choice $\alpha = 100$ is also made) the stabilized gradient method \eqref{EqP1} always converges with the same behaviour, while the gradient method doesn't converge anymore. 

\subsubsection{\bf Example $4.9$ (backtracking stepsize) from Chapter $4$ in \cite{Beck14}:}\label{beck49}
We finally analyse the stabilized gradient method \eqref{EqP1} with the backtracking stepsize by again considering the problem \eqref{Function2} under the choices $a = \frac{1}{100}$, $\bm{x}^{[0]}_{\gamma} = \left(\frac{1}{100}, 1\right)^{t}$ and $\varepsilon = 10^{-5}$, exactly as in the second part of the {\bf Example} $4.9$ in \cite{Beck14} where $201$ iterations have been made before reaching convergence. The backtracking stepsize is defined with the same parameters than in \cite{Beck14}, namely $s = 2$, $\text{alpha} = 0.25$, and $\text{beta} = 0.5$, where the stopping criterion \eqref{SCI2} is used.

We thus present in Table \ref{Tab491} results computed from the stabilized gradient method \eqref{EqP1}. Compared to the gradient method \eqref{Eq2} in the {\bf Example} $4.9$ in \cite{Beck14} that made a number of $201$ iterations before converging, the stabilized method keeps the same behaviour than the previous examples by providing an accurate approximation of the exact solution $\bm{x}^{\star} = (0, 0)^t$, with few iterations as long as $\gamma > \frac{1}{\varepsilon}$. 
\begin{table}[h]
\begin{center}
\begin{tikzpicture}
\draw[gray] [white, fill=gray!15.5]
    (0.0,6).. controls (15.5,6) ..(15.5,6)
    .. controls (15.5, 7.3) ..(15.5, 7.3)
    .. controls (0.0, 7.3) ..(0.0, 7.3)
    .. controls (0.0,6) ..(0.0,6);

\draw[gray] [white, fill=gray!10]
    (0.0,4.0).. controls (15.5,4.0) ..(15.5,4.0)
    .. controls (15.5, 5) .. (15.5, 5)
    .. controls (0.0, 5)  .. (0.0, 5)
    .. controls (0.0,4.0) .. (0.0,4.0);
\draw[gray] [white, fill=gray!10]
    (0.0,2.0).. controls (15.5,2.0) ..(15.5,2.0)
    .. controls (15.5, 3) .. (15.5, 3)
    .. controls (0.0, 3)  .. (0.0, 3)
    .. controls (0.0,2.0) .. (0.0,2.0);
\draw[gray] [white, fill=gray!10]
    (0.0,0.0).. controls (15.5,0.0) ..(15.5,0.0)
    .. controls (15.5, 1) .. (15.5, 1)
    .. controls (0.0, 1)  .. (0.0, 1)
    .. controls (0.0,0.0) .. (0.0,0.0);

\draw[, very thin] (1.5,0) -- (1.5,7.3);
\draw[, very thin] (4.5,0) -- (4.5,7.3);
\draw[, very thin] (7.5,0) -- (7.5,7.3);
\draw[, very thin] (10.5,0) -- (10.5,7.3);
\draw[, very thin] (12.,0) -- (12.,7.3);

\draw[, very thin] (0,0) -- (15.5,0);
\draw[, very thin] (0,1) -- (15.5,1);
\draw[, very thin] (0,2)   -- (15.5,2);
\draw[, very thin] (0,3)   -- (15.5,3);
\draw[, very thin] (0,4)   -- (15.5,4);
\draw[, very thin] (0,5) -- (15.5,5);

\draw[black](0.56, 7.)  node [rotate=40, scale=0.9] {\text{Stab.}};
\draw[black](0.7, 6.65)  node [rotate=40, scale=0.9] {\text{Parameter}};
\draw[black](1., 6.3)  node [rotate=40, scale=1] {$(\gamma)$};
\draw[black](0.75, 5.5)  node [rotate=0, scale=1] {$1$};
\draw[black](0.75, 4.5)  node [rotate=0, scale=1] {$10$};
\draw[black](0.75, 3.5)  node [rotate=0, scale=1] {$10^{2}$};
\draw[black](0.75, 2.5)  node [rotate=0, scale=1] {$10^{5}$};
\draw[black](0.75, 1.5)  node [rotate=0, scale=1] {$10^{7}$};
\draw[black](0.75, 0.5)  node [rotate=0, scale=1] {$10^{10}$};

\draw[black](3, 6.8)   node [rotate=0, scale=0.9] {$\left\|\mathbb{A}_{\frac{1}{100}}\bm{x}^{[k]}_{\gamma} - \bm{b}\right\|_{2}$};
\draw[black](3, 5.5)  node [rotate=0, scale=1] {$9.975445e^{-6}$};
\draw[black](3, 4.5)  node [rotate=0, scale=1] {$9.994624e^{-6}$};
\draw[black](3, 3.5)  node [rotate=0, scale=1] {$9.824913e^{-6}$};
\draw[black](3, 2.5)  node [rotate=0, scale=1] {$1.272837e^{-6}$};
\draw[black](3, 1.5)  node [rotate=0, scale=1] {$1.936525e^{-8}$};
\draw[black](3, 0.5)  node [rotate=0, scale=1] {$1.979998e^{-8}$};

\draw[black](6, 6.7)   node [rotate=0, scale=0.8] {$\frac{\left\|\mathbb{A}_{\frac{1}{100}}\bm{x}^{[k]}_{\gamma} - \bm{b}\right\|_{2}}{\left\|\mathbb{A}_{\frac{1}{100}}\bm{x}^{[0]}_{\gamma} + \bm{b}\right\|_{2}}$};
\draw[black](6, 5.5)  node [rotate=0, scale=1] {$3.526852e^{-4}$};
\draw[black](6, 4.5)  node [rotate=0, scale=1] {$3.533633e^{-4}$};
\draw[black](6, 3.5)  node [rotate=0, scale=1] {$3.473631e^{-4}$};
\draw[black](6, 2.5)  node [rotate=0, scale=1] {$4.500157e^{-5}$};
\draw[black](6, 1.5)  node [rotate=0, scale=1] {$6.846650e^{-7}$};
\draw[black](6, 0.5)  node [rotate=0, scale=1] {$7.000350e^{-7}$};

\draw[black](9., 6.7)   node [rotate=0, scale=1] {$\|\bm{x}^{\star} - \bm{x}^{[k]}_{\gamma}\|_{2}$};
\draw[black](9, 5.5)  node [rotate=0, scale=1] {$4.987722e^{-4}$};
\draw[black](9, 4.5)  node [rotate=0, scale=1] {$4.997312e^{-4}$};
\draw[black](9, 3.5)  node [rotate=0, scale=1] {$4.912456e^{-4}$};
\draw[black](9, 2.5)  node [rotate=0, scale=1] {$6.364183e^{-5}$};
\draw[black](9, 1.5)  node [rotate=0, scale=1] {$9.682625e^{-7}$};
\draw[black](9, 0.5)  node [rotate=0, scale=1] {$9.899990e^{-7}$};

\draw[black](11.25, 6.8)   node [rotate=0, scale=1] {\text{$\#$ Iters}};
\draw[black](11.25, 6.3)   node [rotate=0, scale=1] {$[k]$};
\draw[black](11.25, 5.5)  node [rotate=0, scale=1] {$375$};
\draw[black](11.25, 4.5)  node [rotate=0, scale=1] {$359$};
\draw[black](11.25, 3.5)  node [rotate=0, scale=1] {$253$};
\draw[black](11.25, 2.5)  node [rotate=0, scale=1] {$4$};
\draw[black](11.25, 1.5)  node [rotate=0, scale=1] {$2$};
\draw[black](11.25, 0.5)  node [rotate=0, scale=1] {$1$};

\draw[black](13.8, 6.7)  node [rotate=0, scale=1] {$\bm{x}^{[k]}_{\gamma}$};
\draw[black](13.8, 5.5)  node [rotate=0, scale=1] {$1e^{-3}\left(0.000, 0.499\right)$};
\draw[black](13.8, 4.5)  node [rotate=0, scale=1] {$1e^{-3}\left(0.000, 0.499\right)$};
\draw[black](13.8, 3.5)  node [rotate=0, scale=1] {$1e^{-3}\left(0.000, 0.491\right)$};
\draw[black](13.8, 2.5)  node [rotate=0, scale=1] {$1e^{-4}\left(0.000, 0.636\right)$};
\draw[black](13.8, 1.5)  node [rotate=0, scale=1] {$1e^{-6}\left(0.000, 0.968\right)$};
\draw[black](13.8, 0.5)  node [rotate=0, scale=1] {$1e^{-6}\left(0.000, 0.990\right)$};
\end{tikzpicture} 
\caption{Numerical results obtained from the stabilized method \eqref{EqP1} with the data $\mathbb{A}_{a}$ and $\bm{b}$ given in \eqref{MatV1}, where $a = \frac{1}{100}$ ($\textbf{cond}\left(\mathbb{A}_{\frac{1}{100}}\right) = 100$ and $\textbf{rank}\left(\mathbb{A}_{\frac{1}{100}}\right) = 2$). The computations are made with respect to the different values of the stabilization parameter $\gamma = 10^{i}$, $i = 0$, $1$, $2$, $5$, $7$, $10$, the initial iteration fixed to $\bm{x}^{[0]}_{\gamma} = \left(\frac{1}{100}, 1\right)^{t}$, and the backtracking stepsize provided in the {\bf Example} $4.9$ in \cite{Beck14}. The maximum iterations is $k_{\text{max}} = 1\,000$, and the stopping criterion \eqref{SCI2} is used with $\varepsilon = 10^{-5}$.}
\label{Tab491}
\end{center} 
\end{table}

\subsection{\bf Ill-conditioned linear systems}\label{NEx3}
In this third part, we investigated the ill-conditioned linear systems under both small and large scale problems.  
\subsubsection{\bf Small-scale inverse problems}\label{NEx3SSP} 
In the analyzing the efficiency and the accuracy of the stabilized gradient method \eqref{EqP1}, we also investigated severely ill-conditioned numerical examples. The problems studied throughout this part are standard small-scale inverse problems known as the {\bf shaw}, the {\bf heat} and the {\bf gravity} problems \cite{HansenBook10,HansenBook98}. These test problems are classical examples of ill-posed problems, formulated throught the following Fredholm integral equation of the first kind (see \cite{Groetsch84})
\begin{align}\label{Fredholm1}
    \int_{a}^{b}K(s, t)x^{\star}(t)dt = f(s), \;\, c \leq s \leq d,
\end{align}
where the kernel $K$ and the right hand side $f$ are given, and $x^{\star}$ is the desired unknown solution. \\
To apply a quadrature formula in the approximation of the integral equation in \eqref{Fredholm1}, we use as in \cite{HansenBook98} the midpoint quadrature-collocation rule with $n$ equidistantly spaced points $s_{i}$ and $t_{i}$
\begin{align}\label{Fredholm5}
   s_{i} = c + ih, \;\, t_{i} = a + ih, \;\, h = \frac{b - a}{n}, \;\, i = 1, \cdots, n.
\end{align}
Rewriting the integral equation \eqref{Fredholm1} through the discretizations in \eqref{Fredholm5}, we obtain the linear system 
\begin{align}\label{Fredholm7}
    \sum_{j=1}^{n}w_jK(s_i, t^{\star}_{j})x^{\star}(t^{\star}_j) = f(s_i), \;\, i = 1, \cdots, n, 
\end{align}
where each midpoint is given by $t^{\star}_j = t_j - \frac{h}{2}$, and $w_j = h$ for $j = 1, \cdots, n$ are the weights of the quadrature formula. We thus compute the $n\times n$ matrix $\mathbf{A} = (a_{ij})_{1 \leq i, j \leq n}$, the $n$-vector right hand $\bm{b} = (b_{i})_{1 \leq i \leq n}$, and the $n$-vector exact solution $\bm{x}^{\star} = (x^{\star}_{j})_{1 \leq j \leq n}$ as follows
\begin{align}\label{Fredholm8}
    a_{ij} = w_jK(s_i, t^{\star}_{j}), \;\, b_{i} = f(s_i), \;\, x^{\star}_{j} = x^{\star}(t^{\star}_j), \;\, i = 1, \cdots, n,  \;\, j = 1, \cdots, n.
\end{align}

As we may seen in \eqref{Reslt1Fin2} - \eqref{Reslt1FinBis1} and \eqref{Err1} - \eqref{Err1Bis}, and through the previous numerical examples \ref{NEx1} and \ref{NEx2}, the value of the stepsize $\alpha$ has a limit influence on the convergence of the stabilized gradient method \eqref{EqP1}. What really matter is the value of the stabilization parameter $\gamma$, which should be fixed as large as possible for reducing the residual. 

$\bullet$ \,{\bf The shaw test problem}\label{shaw}

The desired solution considered here is defined by the following analytical expression
\begin{align}\label{Fredholm2}
   \bm{x}^{\star}(t) = a_1\exp\left(-c_{1}\left(t - t_1\right)^2\right) + a_{2}\exp\left(-c_{2}\left(t - t_{2}\right)^2\right), \;\, -\frac{\pi}{2} \leq t \leq \frac{\pi}{2},
\end{align}
where the parameters $a_1$, $a_2$, $c_1$, $c_2$, $x_1$ and $x_2$ are constants defined as follows 
\begin{align}\label{SetPar2}
   \begin{aligned}
        a_1 & = 2, \;\, c_1 = 6, \;\, x_1 =  0.8, \\
        a_2 & = 1, \;\, c_2 = 2, \;\, x_2 = -0.5,
   \end{aligned}
\end{align}
and where parameters on the Fredholm integral equation \eqref{Fredholm1} are given as $a = c = -\frac{\pi}{2}$ $b = d = \frac{\pi}{2}$. As for the exact kernel $K$, it is analytically defined by the following analytical expression 
\begin{align}\label{Fredholm3}
   K(s, t) = \left(\cos(s) + \cos(t)\right)^2\left(\frac{\sin(\psi(s, t))}{\psi(s, t)}\right)^2, \;\, -\frac{\pi}{2} \leq s \leq \frac{\pi}{2},
\end{align}
with the given function $\psi(s, t) = \pi\left(\sin(s) + \cos(t)\right)$. Throughout this numerical example, we use the Matlab code \textsf{shaw(n)} provided in {\it Regularization Tools} \cite{HansenBook07} to obtain the data $\mathbf{A}$, $\bm{b}$, and $\bm{x}^{\star}$. 

For having an insight on the convergence of the stabilized gradient method \eqref{EqP1} applied to the {\bf shaw} test problem, we present in Figure \ref{RayonShawHeatGravity}-{\small\colorbox{black}{\color{white}1}} the spectral radius $\varrho(\mathbf{M}^{-1}_{\gamma}\mathbf{N}_{k})$ with respect to the stabilization parameter $\gamma$, for the different values of the stepsize $\alpha = 10^{i}, i = -3, -2, 0, 2, 3$, and under the dimension $n = 1\,000$. We may thus note how the stabilized method should perform less good under the {\bf shaw} test problem (or generally under Ill-conditioned linear systems). 
\begin{figure}[h]  
    \centering  
    \scalebox{0.43}{
%
%
\definecolor{mycolor1}{rgb}{0.00000,1.00000,1.00000}%
\definecolor{mycolor2}{rgb}{1.00000,0.00000,1.00000}%
\begin{tikzpicture}

\begin{axis}[%
width=4.521in,
height=3.566in,
at={(0.758in,0.481in)},
scale only axis,
xmode=log,
xmin=1,
xmax=10000000000,
xlabel style={font=\color{white!15!black}, scale=2},
xlabel={$\gamma$},
ymin=0.9,
ymax=2,
axis background/.style={fill=white},
title style={font=\bfseries, scale=0.9},
title={\normalsize\colorbox{black}{\color{white}1 | Spectral radius $\varrho(\mathbf{M}^{-1}_{\gamma}\mathbf{N}_{k})$ under the shaw test problem}},
xmajorgrids,
ymajorgrids,
legend style={legend cell align=left, align=left, draw=white!15!black, font=\Large},
grid style = loosely dashed,
]
\addplot [color=red, line width=1.85pt]
  table[row sep=crcr]{%
1 1.00000007232485\\
10000000000 1.00000394942935\\
};
\addlegendentry[scale=0.8]{$\varrho(\mathbf{M}^{-1}_{\gamma}\mathbf{N}_{k}), \alpha = 10^{-3}$}

\addplot [color=green, dashed, line width=1.85pt]
  table[row sep=crcr]{%
1 1.00002458567661\\
19338917.5045523  1.00000000728083\\
10000000000 1.00000394942939\\
};
\addlegendentry[scale=0.8]{$\varrho(\mathbf{M}^{-1}_{\gamma}\mathbf{N}_{k}), \alpha = 10^{-2}$}

\addplot [color=blue, dotted, line width=1.85pt]
  table[row sep=crcr]{%
1 1.20665357407072\\
1.12266777351081  1.18713755722381\\
1.26038292967973  1.16596633895148\\
1.41499129743458  1.14308018797178\\
1.58856512942805  1.11843407735793\\
1.78343087693191  1.09200110333312\\
2.00220037181558  1.06377590238846\\
2.24780583354873  1.03377789128271\\
2.52353917043477  1.02294321244801\\
4.50055767570051  1.02173337207901\\
7.14942898659759  1.02011686875869\\
10.1163797976621  1.01831231632549\\
14.3145893752348  1.01576977143443\\
18.0418640939208  1.01352305349466\\
22.7396575235793  1.01070544669048\\
28.6606761694826  1.00717642101283\\
32.1764175025074  1.00509261742542\\
36.1234269970944  1.00366852024733\\
81.1984499318401  1.00281114811307\\
144.811822767454  1.00160364731843\\
258.261876068267  1.00007689371844\\
10000000000 1.0000039494262\\
};
\addlegendentry[scale=0.8]{$\varrho(\mathbf{M}^{-1}_{\gamma}\mathbf{N}_{k}), \alpha = 1$}

\addplot [color=mycolor1, dashed, line width=1.85pt]
  table[row sep=crcr]{%
1162.32246867985  2.0322574308224\\
1304.9019780144 1.93476165563073\\
1464.97139830729  1.83588283146568\\
2072.92177959537  1.53745576261034\\
2327.20247896041  1.43958017718546\\
2612.67522556333  1.36783676956275\\
2933.16627839005  1.35995484116919\\
3292.97125509715  1.35121359255031\\
3696.91270719503  1.34153300751504\\
4150.40475785047  1.33082888494261\\
4659.52566866469  1.31901343711133\\
5231.09930805626  1.30599613107507\\
5872.78661318949  1.29168481069119\\
6593.18827133354  1.27598713452156\\
7401.95999691565  1.25881235907447\\
8309.94194935339  1.24007348682043\\
9329.3040262847 1.2196897832672\\
10473.7089795945  1.1975896464962\\
11758.4955405216  1.17371378574654\\
13200.8840083142  1.14801863334249\\
14820.2070579886  1.12047987783831\\
16638.1688607613  1.09109596815641\\
18679.1359902078  1.05989140232821\\
20970.4640132323  1.02691958491404\\
23542.8641432242  1.00770233200548\\
105956.017927762  1.00719110312251\\
238168.555197616  1.00637203921998\\
424757.155253691  1.00521837667962\\
674754.405311071  1.00367680315886\\
954771.611420807  1.00195572006862\\
1203377.84077759  1.00043263874872\\
1516716.88847092  1.00023621916405\\
10000000000 1.00000394916166\\
};
\addlegendentry[scale=0.8]{$\varrho(\mathbf{M}^{-1}_{\gamma}\mathbf{N}_{k}), \alpha = 10^{2}$}

\addplot [color=mycolor2, dotted, line width=1.85pt]
  table[row sep=crcr]{%
84066.5288561832  2.06372303546155\\
94378.7827777539  1.91901292107719\\
105956.017927762  1.778968148215\\
118953.406737032  1.64425515870736\\
133545.15629299 1.51542241259685\\
149926.843278605  1.39289670527198\\
168318.035333096  1.2769842000398\\
188965.233969121  1.16787554506111\\
212145.178491064  1.07707718951104\\
378346.261713193  1.075977070542\\
601027.678207039  1.07450661451147\\
850448.93418027 1.07286434916814\\
1203377.84077759  1.07054911671861\\
1516716.88847092  1.06850195385465\\
1911644.0753857 1.06593287067337\\
2409403.56023953  1.06271238701582\\
3036771.11803546  1.0586809594407\\
3827494.47851632  1.05364317955205\\
4297004.70432085  1.05067450376198\\
4824108.70416537  1.04736153696055\\
5415871.37807949  1.04366698234378\\
6080224.26164943  1.03955015951709\\
6826071.83427241  1.03496686325061\\
7663410.86800746  1.02986926228108\\
8603464.41668453  1.02420585358107\\
9658832.24115871  1.01792148945333\\
10843659.6868961  1.01095749775534\\
12173827.2773966  1.00325191981001\\
13667163.5646201  1.0023710477374\\
87035913.6148518  1.00192857122604\\
219638537.241654  1.00112985405679\\
439760360.930271  1.00040146548441\\
10000000000 1.00000394831141\\
};
\addlegendentry[scale=0.8]{$\varrho(\mathbf{M}^{-1}_{\gamma}\mathbf{N}_{k}), \alpha = 10^{3}$}

\end{axis}
\end{tikzpicture}%
     %
%
%
\definecolor{mycolor1}{rgb}{0.00000,1.00000,1.00000}%
\definecolor{mycolor2}{rgb}{1.00000,0.00000,1.00000}%
\begin{tikzpicture}

\begin{axis}[%
width=4.521in,
height=3.566in,
at={(0.758in,0.481in)},
scale only axis,
xmode=log,
xmin=1,
xmax=10000000000,
xlabel style={font=\color{white!15!black}, scale=2},
xlabel={$\gamma$},
ymin=0.9,
ymax=2,
axis background/.style={fill=white},
title style={font=\bfseries, scale=0.9},
title={\normalsize\colorbox{black}{\color{white}2 | Spectral radius $\varrho(\mathbf{M}^{-1}_{\gamma}\mathbf{N}_{k})$ under the heat test problem}},
xmajorgrids,
ymajorgrids,
legend style={legend cell align=left, align=left, draw=white!15!black, font=\Large},
grid style = loosely dashed,
]
\addplot [color=red, line width=1.85pt]
  table[row sep=crcr]{%
1 1.00000003932289\\
10000000000 1.00000004095022\\
};
\addlegendentry[scale=0.8]{$\varrho(\mathbf{M}^{-1}_{\gamma}\mathbf{N}_{k}), \alpha = 10^{-3}$}

\addplot [color=green, dashed, line width=1.85pt]
  table[row sep=crcr]{%
1 1.00000044426473\\
10000000000 1.00000004095023\\
};
\addlegendentry[scale=0.8]{$\varrho(\mathbf{M}^{-1}_{\gamma}\mathbf{N}_{k}), \alpha = 10^{-2}$}

\addplot [color=blue, dotted, line width=1.85pt]
  table[row sep=crcr]{%
1 1.02103855713061\\
1.58856512942805  1.02033939368505\\
2.52353917043477  1.01895887437155\\
4.00880632889846  1.01702998313885\\
5.67242606849199  1.01508211530608\\
6.36824994471859  1.01422530920902\\
7.14942898659759  1.0144232597689\\
8.02643352225717  1.01322217248222\\
11.3573335834311  1.01162535862431\\
14.3145893752348  1.01000254911374\\
20.2550193923067  1.00666977357328\\
28.6606761694826  1.00263402378236\\
36.1234269970944  1.00004417949049\\
10000000000 1.00000004095189\\
};
\addlegendentry[scale=0.8]{$\varrho(\mathbf{M}^{-1}_{\gamma}\mathbf{N}_{k}), \alpha = 1$}

\addplot [color=mycolor1, dashed, line width=1.85pt]
  table[row sep=crcr]{%
517.092024289676  2.03163641979487\\
580.522551609491  1.98501545272943\\
651.733960488242  1.93036786758179\\
731.680714342721  1.87208514894527\\
821.434358491942  1.81258295189678\\
922.197882333434  1.75364887082788\\
1035.32184329566  1.69636114441989\\
1162.32246867985  1.64071091458385\\
1304.9019780144 1.58606317678739\\
1644.67617799466  1.47829781277451\\
1846.42494289555  1.42516199161013\\
2072.92177959537  1.37267516642243\\
2327.20247896041  1.32092921714303\\
2612.67522556333  1.26997521520605\\
2933.16627839005  1.21983981527076\\
3292.97125509715  1.17054247284745\\
3696.91270719503  1.12210356770681\\
4150.40475785047  1.07454514829568\\
4659.52566866469  1.02788892648723\\
5231.09930805626  1.0044448133498\\
47137.5313411673  1.00423868409663\\
212145.178491064  1.00341419549883\\
601027.678207039  1.00230059304611\\
1350993.52119803  1.00102337943101\\
2409403.56023953  1.00000890791497\\
10000000000 1.0000000458116\\
};
\addlegendentry[scale=0.8]{$\varrho(\mathbf{M}^{-1}_{\gamma}\mathbf{N}_{k}), \alpha = 10^{2}$}

\addplot [color=mycolor2, dotted, line width=1.85pt]
  table[row sep=crcr]{%
29673.0240818887  2.10106099858406\\
33312.9478793467  1.97183537784294\\
37399.373024788 1.84726807647121\\
41987.0708444391  1.72781151510173\\
47137.5313411673  1.61413736275697\\
52919.7873595844  1.50670026383946\\
59411.3398496504  1.40530938783298\\
66699.1966303011  1.30931828245271\\
74881.0385759003  1.21815393623126\\
84066.5288561832  1.1315106997856\\
94378.7827777539  1.04928273570323\\
105956.017927762  1.04531613433474\\
212145.178491064  1.04482809911982\\
378346.261713193  1.04375725049104\\
601027.678207039  1.042191866959\\
1203377.84077759  1.03874698939681\\
1702769.1722259 1.0362015794002\\
3036771.11803546  1.03078941853986\\
5415871.37807949  1.02401417663122\\
6826071.83427241  1.0211430373759\\
8603464.41668453  1.01815706043651\\
10843659.6868961  1.01431296401454\\
13667163.5646201  1.00982225871063\\
21711179.4569451  1.00008838402581\\
10000000000 1.00001111102493\\
};
\addlegendentry[scale=0.8]{$\varrho(\mathbf{M}^{-1}_{\gamma}\mathbf{N}_{k}), \alpha = 10^{3}$}

\end{axis}
\end{tikzpicture}%
%
%
\definecolor{mycolor1}{rgb}{0.00000,1.00000,1.00000}%
\definecolor{mycolor2}{rgb}{1.00000,0.00000,1.00000}%
\begin{tikzpicture}

\begin{axis}[%
width=4.521in,
height=3.566in,
at={(0.758in,0.481in)},
scale only axis,
xmode=log,
xmin=1,
xmax=10000000000,
xlabel style={font=\color{white!15!black}, scale=2},
xlabel={$\gamma$},
ymin=0.9,
ymax=2,
axis background/.style={fill=white},
title style={font=\bfseries, scale=0.9},
title={\normalsize\colorbox{black}{\color{white}3 | Spectral radius $\varrho(\mathbf{M}^{-1}_{\gamma}\mathbf{N}_{k})$ under the gravity test problem}},
xmajorgrids,
ymajorgrids,
legend style={legend cell align=left, align=left, draw=white!15!black, font=\Large},
grid style = loosely dashed,
]
\addplot [color=red, line width=1.85pt]
  table[row sep=crcr]{%
1 1.00000000000001\\
10000000000 1.00000019253827\\
};
\addlegendentry[scale=0.8]{$\varrho(\mathbf{M}^{-1}_{\gamma}\mathbf{N}_{k}), \alpha = 10^{-3}$}

\addplot [color=green, dashed, line width=1.85pt]
  table[row sep=crcr]{%
1 1.00000000000001\\
10000000000 1.00000019253826\\
};
\addlegendentry[scale=0.8]{$\varrho(\mathbf{M}^{-1}_{\gamma}\mathbf{N}_{k}), \alpha = 10^{-2}$}

\addplot [color=blue, dotted, line width=1.85pt]
  table[row sep=crcr]{%
1 1.00000000000001\\
10000000000 1.00000019253815\\
};
\addlegendentry[scale=0.8]{$\varrho(\mathbf{M}^{-1}_{\gamma}\mathbf{N}_{k}), \alpha = 1$}

\addplot [color=mycolor1, dashed, line width=1.85pt]
  table[row sep=crcr]{%
258.465822845482  2.11\\
289.942285388288  1.91367834503374\\
325.508859983506  1.73166759329595\\
365.438307095726  1.56460331349027\\
410.265810582719  1.41170131914029\\
460.592204114511  1.27213131130569\\
517.092024289676  1.14503859442351\\
580.522551609491  1.02956235494547\\
651.733960488242  1.00000000000002\\
10000000000 1.00000019253056\\
};
\addlegendentry[scale=0.8]{$\varrho(\mathbf{M}^{-1}_{\gamma}\mathbf{N}_{k}), \alpha = 10^{2}$}

\addplot [color=mycolor2, dotted, line width=1.85pt]
  table[row sep=crcr]{%
19436.6863233119  2.11\\
20970.4640132323  1.968982161379\\
23542.8641432242  1.77315352572968\\
26430.814869741 1.59505462694525\\
29673.0240818887  1.43341825910281\\
33312.9478793467  1.3826490368643\\
37399.373024788 1.3498208856507\\
41987.0708444391  1.31477499025559\\
47137.5313411673  1.27753703948832\\
52919.7873595844  1.23816709050543\\
59411.3398496504  1.1967623077409\\
66699.1966303011  1.15345865396588\\
74881.0385759003  1.10843127508415\\
84066.5288561832  1.06189337297744\\
94378.7827777539  1.0140934428891\\
105956.017927762  1.00000000000247\\
10000000000 1.00000019251151\\
};
\addlegendentry[scale=0.8]{$\varrho(\mathbf{M}^{-1}_{\gamma}\mathbf{N}_{k}), \alpha = 10^{3}$}

\end{axis}
\end{tikzpicture}%
}    
    \caption{\footnotesize The spectral radius $\varrho(\mathbf{M}^{-1}_{\gamma}\mathbf{N}_{k})$ of the stabilized gradient method \eqref{EqP1} applied to the discret linear system \eqref{Fredholm7}, for the Shaw (in {\tiny\colorbox{black}{\color{white}1}}), the Heat (in {\tiny\colorbox{black}{\color{white}2}}), and the Gravity (in {\tiny\colorbox{black}{\color{white}3}}) test problems. These are represented with respect to the stabilization parameter $\gamma$, for the different values of the stepsize $\alpha = 10^{i}, i = -3, -2, 0, 2, 3$, and under the dimension $n = 1\,000$.}
    \label{RayonShawHeatGravity}
\end{figure}
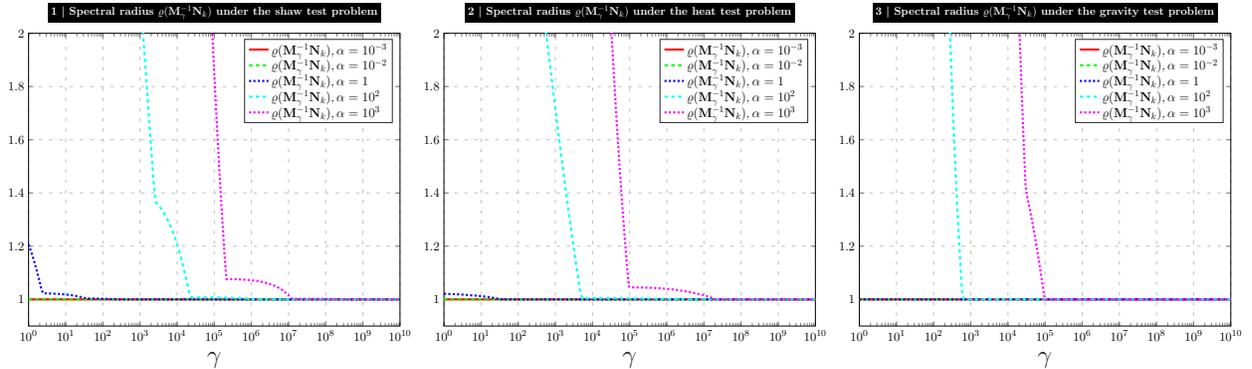
\begin{table}[h]
\begin{center} 
\begin{tikzpicture}
\draw[gray] [white, fill=gray!15.5]
    (0.0,6).. controls (15.5,6) ..(15.5,6)
    .. controls (15.5, 7.25) ..(15.5, 7.25)
    .. controls (0.0, 7.25) ..(0.0, 7.25)
    .. controls (0.0,6) ..(0.0,6);

\draw[gray] [white, fill=gray!10]
    (0.0,4.0).. controls (15.5,4.0) ..(15.5,4.0)
    .. controls (15.5, 5) .. (15.5, 5)
    .. controls (0.0, 5)  .. (0.0, 5)
    .. controls (0.0,4.0) .. (0.0,4.0);
\draw[gray] [white, fill=gray!10]
    (0.0,2.0).. controls (15.5,2.0) ..(15.5,2.0)
    .. controls (15.5, 3) .. (15.5, 3)
    .. controls (0.0, 3)  .. (0.0, 3)
    .. controls (0.0,2.0) .. (0.0,2.0);
\draw[gray] [white, fill=gray!10]
    (0.0,0.0).. controls (15.5,0.0) ..(15.5,0.0)
    .. controls (15.5, 1) .. (15.5, 1)
    .. controls (0.0, 1)  .. (0.0, 1)
    .. controls (0.0,0.0) .. (0.0,0.0);

\draw[, very thin] (1.5,-1) -- (1.5,7.25);
\draw[, very thin] (4.5,-1) -- (4.5,7.25);
\draw[, very thin] (7.5,-1) -- (7.5,7.25);
\draw[, very thin] (10.5,-1) -- (10.5,7.25);
\draw[, very thin] (13.5,-1) -- (13.5,7.25);

\draw[, very thin] (0,-1) -- (15.5,-1);
\draw[, very thin] (0,0) -- (15.5,0);
\draw[, very thin] (0,1) -- (15.5,1);
\draw[, very thin] (0,2)   -- (15.5,2);
\draw[, very thin] (0,3)   -- (15.5,3);
\draw[, very thin] (0,4)   -- (15.5,4);
\draw[, very thin] (0,5) -- (15.5,5);

\draw[black](0.56, 7.)  node [rotate=40, scale=0.9] {\text{Stab.}};
\draw[black](0.7, 6.65)  node [rotate=40, scale=0.9] {\text{Parameter}};
\draw[black](1., 6.3)  node [rotate=40, scale=1] {$(\gamma)$};
\draw[black](0.75, 5.5)  node [rotate=0, scale=1] {$10^{3}$};
\draw[black](0.75, 4.5)  node [rotate=0, scale=1] {$10^{4}$};
\draw[black](0.75, 3.5)  node [rotate=0, scale=1] {$10^{5}$};
\draw[black](0.75, 2.5)  node [rotate=0, scale=1] {$10^{6}$};
\draw[black](0.75, 1.5)  node [rotate=0, scale=1] {$10^{10}$};
\draw[black](0.75, 0.5)  node [rotate=0, scale=1] {$10^{12}$};
\draw[black](0.75, -0.5)  node [rotate=0, scale=1] {$10^{14}$};

\draw[black](3, 6.7)   node [rotate=0, scale=1] {$\|\mathbf{A}\bm{x}^{[k]}_{\gamma} - \bm{b}\|_{n}$};
\draw[black](3, 5.5)  node [rotate=0, scale=1] {$7.355131e^{-04}$};
\draw[black](3, 4.5)  node [rotate=0, scale=1] {$7.275623e^{-04}$};
\draw[black](3, 3.5)  node [rotate=0, scale=1] {$6.441531e^{-04}$};
\draw[black](3, 2.5)  node [rotate=0, scale=1] {$5.556602e^{-04}$};
\draw[black](3, 1.5)  node [rotate=0, scale=1] {$1.839275e^{-06}$};
\draw[black](3, 0.5)  node [rotate=0, scale=1] {$2.058136e^{-07}$};
\draw[black](3, -0.5)  node [rotate=0, scale=1] {$4.931643e^{-09}$};

\draw[black](6, 6.6)   node [rotate=0, scale=1] {$\frac{\|\mathbf{A}\bm{x}^{[k]}_{\gamma} - \bm{b}\|_{n}}{\|\mathbf{A}\bm{x}^{[0]}_{\gamma} - \bm{b}\|_{n}}$};
\draw[black](6, 5.5)  node [rotate=0, scale=1] {$9.977567e^{-06}$};
\draw[black](6, 4.5)  node [rotate=0, scale=1] {$9.869711e^{-06}$};
\draw[black](6, 3.5)  node [rotate=0, scale=1] {$8.738228e^{-06}$};
\draw[black](6, 2.5)  node [rotate=0, scale=1] {$7.537782e^{-06}$};
\draw[black](6, 1.5)  node [rotate=0, scale=1] {$2.495059e^{-08}$};
\draw[black](6, 0.5)  node [rotate=0, scale=1] {$2.791954e^{-09}$};
\draw[black](6, -0.5)  node [rotate=0, scale=1] {$6.689996e^{-11}$};

\draw[black](9., 6.7)   node [rotate=0, scale=1] {$\|\bm{x}^{\star} - \bm{x}^{[k]}_{\gamma}\|_{n}$};
\draw[black](9, 5.5)  node [rotate=0, scale=1] {$1.161398e^{+00}$};
\draw[black](9, 4.5)  node [rotate=0, scale=1] {$1.149915e^{+00}$};
\draw[black](9, 3.5)  node [rotate=0, scale=1] {$1.119442e^{+00}$};
\draw[black](9, 2.5)  node [rotate=0, scale=1] {$1.082791e^{+00}$};
\draw[black](9, 1.5)  node [rotate=0, scale=1] {$5.634969e^{-01}$};
\draw[black](9, 0.5)  node [rotate=0, scale=1] {$2.179780e^{-01}$};
\draw[black](9, -0.5)  node [rotate=0, scale=1] {$4.101602e^{+00}$};

\draw[black](12., 6.6)   node [rotate=0, scale=1] {$\frac{\|\bm{x}^{\star} - \bm{x}^{[k]}_{\gamma}\|_{n}}{\|\bm{x}^{\star}\|_{n}}$};
\draw[black](12, 5.5)  node [rotate=0, scale=1] {$3.679279e^{-02}$};
\draw[black](12, 4.5)  node [rotate=0, scale=1] {$3.642899e^{-02}$};
\draw[black](12, 3.5)  node [rotate=0, scale=1] {$3.546363e^{-02}$};
\draw[black](12, 2.5)  node [rotate=0, scale=1] {$3.430252e^{-02}$};
\draw[black](12, 1.5)  node [rotate=0, scale=1] {$1.785143e^{-02}$};
\draw[black](12, 0.5)  node [rotate=0, scale=1] {$6.905484e^{-03}$};
\draw[black](12, -0.5)  node [rotate=0, scale=1] {$1.299376e^{-01}$};

\draw[black](14.5, 6.85)   node [rotate=0, scale=1] {\text{$\#$ Iters}};
\draw[black](14.5, 6.3)   node [rotate=0, scale=1] {$[k]$};
\draw[black](14.5, 5.5)  node [rotate=0, scale=1] {$222$};
\draw[black](14.5, 4.5)  node [rotate=0, scale=1] {$37$};
\draw[black](14.5, 3.5)  node [rotate=0, scale=1] {$5$};
\draw[black](14.5, 2.5)  node [rotate=0, scale=1] {$1$};
\draw[black](14.5, 1.5)  node [rotate=0, scale=1] {$1$};
\draw[black](14.5, 0.5)  node [rotate=0, scale=1] {$1$};
\draw[black](14.5, -0.5)  node [rotate=0, scale=1] {$1$};

\end{tikzpicture} 
\caption{Numerical results from the stabilized method \eqref{EqP1} of the {\bf shaw} problem for different values of the stabilization parameter $\gamma$, under the fixed dimension $n=1\,000$, stepsize $\alpha = 1$, maximum iterations number $k_{\text{max}} = n$ and stopping tolerance $\varepsilon = 10^{-5}$ in \eqref{CTR1}. The matrix $\mathbf{A}$, the vector $\bm{b}$ and the exact solution $\bm{x}^{\star}$ are computed from the Matlab code \textsf{shaw(n)} provided in {\it Regularization Tools} \cite{HansenBook07}, where $\textbf{cond}(\mathbf{A}) = 7.697e^{+20}$ and $\textbf{rank}(\mathbf{A}) = 20$.}
\label{Tabshaw1}
\end{center}
\end{table}
\begin{figure}[h]
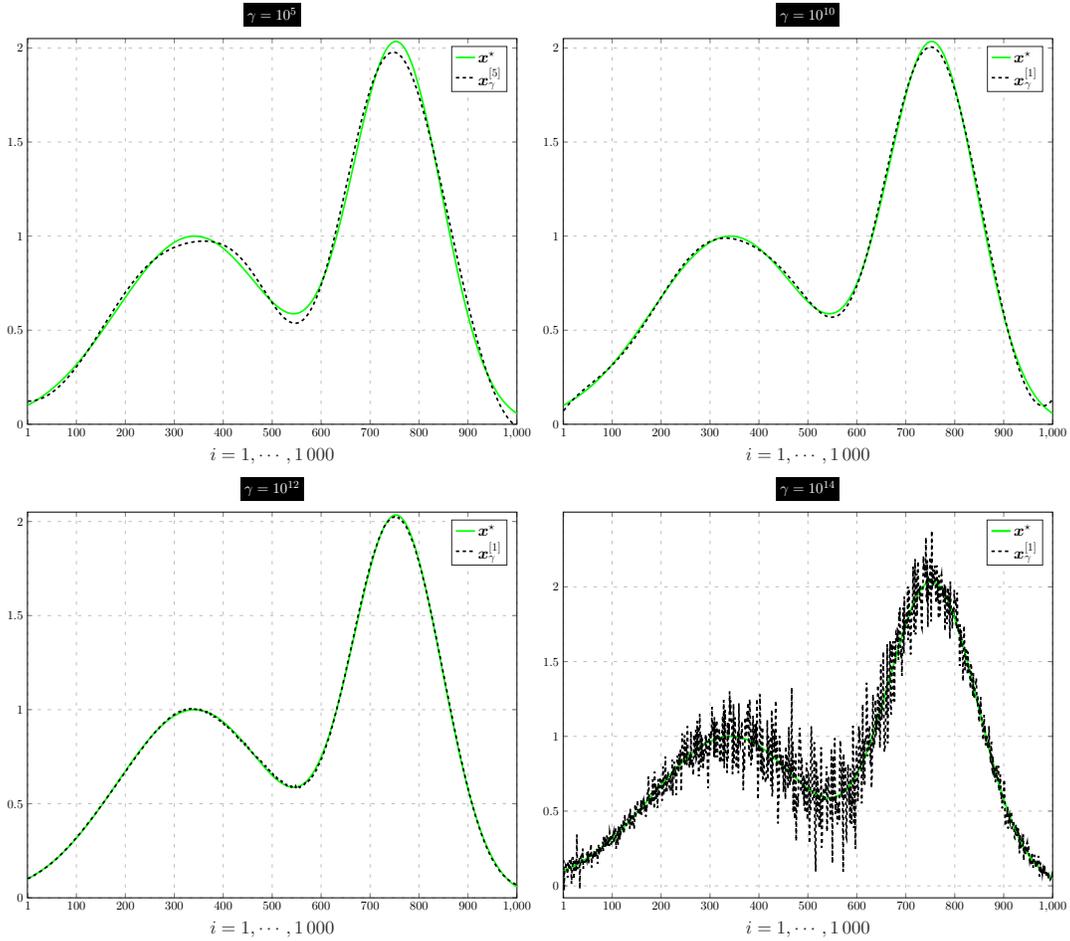
  
    \centering  
    \scalebox{0.425}{
        \begin{tikzpicture}

\begin{axis}[%
width=6.028in,
height=4.754in,
at={(1.011in,0.642in)},
scale only axis,
xmin=0,
xmax=1000,
xlabel style={font=\color{white!15!black}, scale=1.5},
xlabel={$i = 1, \cdots, 1\,000$},
ymin=0,
ymax=2.05,
axis background/.style={fill=white},
title style={font=\bfseries, scale=1.2},
title={\normalsize\colorbox{black}{\color{white}$\gamma = 10^{5}$}},
xmajorgrids,
ymajorgrids,
grid style={dashed},
legend style={legend cell align=left, align=left, draw=white!15!black, font=\Large},
xtick={1, 100, 200, 300, 400, 500, 600, 700, 800, 900, 1000},
ytick={0, 0.5, 1, 1.5, 2},
grid style = loosely dashed,
]
\addplot [color=green, line width=1.5pt]
  table[row sep=crcr]{%
1	0.101622890399199\\
12	0.11752833951391\\
23	0.135275479236952\\
33	0.153088617997582\\
43	0.172564794229629\\
53	0.193752344974769\\
63	0.216684183111852\\
73	0.241375356851108\\
83	0.267820687603376\\
93	0.295992554424629\\
103	0.325838896876235\\
113	0.357281509818222\\
124	0.393585091233035\\
136	0.435032224574229\\
149	0.481772135913275\\
164	0.537501802126258\\
189	0.632466767005212\\
209	0.707699006442226\\
222	0.754914931057783\\
233	0.793174965110097\\
243	0.826205074325685\\
252	0.85421848086196\\
260	0.87755681689157\\
268	0.899258028728809\\
276	0.919170868476272\\
283	0.935016420667694\\
290	0.94929751546772\\
297	0.961934875547513\\
304	0.972857783071731\\
311	0.982004779482509\\
318	0.989324296191853\\
325	0.994775211731167\\
332	0.998327332785152\\
339	0.999961798646723\\
346	0.999671410938845\\
353	0.99746089293285\\
360	0.993347085388791\\
367	0.987359088503126\\
374	0.97953836216675\\
381	0.969938799221154\\
388	0.958626788609308\\
395	0.945681287092157\\
403	0.929004356996302\\
411	0.910471278552791\\
419	0.890261648613318\\
428	0.865774908098274\\
438	0.83677618930119\\
450	0.800182029003054\\
486	0.689049606611775\\
495	0.663833924558162\\
502	0.645784976299637\\
509	0.629479345284949\\
515	0.617164743330022\\
521	0.606634333231341\\
526	0.599398340736911\\
531	0.593715922713841\\
536	0.589733499107865\\
541	0.587598204555206\\
546	0.587456671721952\\
550	0.588876323472505\\
554	0.591737280914913\\
558	0.596110172588055\\
562	0.60206305674194\\
566	0.609660780227387\\
570	0.618964301562869\\
574	0.630029982465771\\
578	0.642908853022277\\
582	0.65764585654631\\
586	0.674279081046734\\
590	0.692838985057165\\
594	0.713347626370251\\
598	0.73581790293656\\
602	0.760252815821559\\
606	0.78664476464337\\
610	0.814974886319987\\
614	0.845212448221559\\
618	0.877314306935205\\
622	0.911224443792548\\
626	0.946873588074595\\
630	0.984178938382229\\
634	1.0230439920399\\
638	1.06335849158484\\
643	1.1155997957427\\
648	1.16963961440581\\
654	1.23643212952459\\
661	1.31629733285456\\
683	1.5693146070538\\
688	1.6246306312928\\
693	1.6781785796336\\
697	1.71943878925003\\
701	1.75904715141803\\
705	1.79677834167012\\
709	1.83241324560959\\
712	1.85763639604761\\
715	1.8814759200892\\
718	1.90384980400847\\
721	1.92467993991761\\
724	1.94389250739266\\
727	1.96141833701722\\
730	1.97719325329979\\
733	1.99115839454203\\
736	2.0032605073809\\
739	2.01345221389602\\
742	2.0216922493529\\
745	2.02794566886507\\
748	2.03218402146865\\
751	2.0343854903449\\
754	2.03453499816715\\
757	2.03262427680579\\
760	2.02865190088937\\
763	2.02262328498682\\
766	2.0145506444444\\
769	2.00445292018424\\
772	1.99235566803566\\
775	1.97829091343317\\
778	1.9622969725674\\
781	1.94441824132093\\
784	1.92470495354894\\
787	1.9032129104819\\
790	1.88000318322713\\
793	1.85514179052564\\
796	1.82869935408246\\
799	1.80075073392823\\
802	1.77137464638531\\
805	1.7406532673084\\
809	1.6977464454601\\
813	1.6528089157913\\
817	1.60605631120882\\
821	1.5577089887463\\
826	1.49537255453231\\
831	1.43133148782363\\
837	1.35284512443252\\
847	1.21996556196962\\
858	1.07413444482813\\
864	0.996079083795962\\
870	0.919786911144797\\
875	0.857908201893565\\
880	0.797837224925274\\
885	0.739787911687813\\
889	0.694927238604009\\
893	0.651555592232739\\
897	0.609738330062555\\
901	0.569528223327211\\
905	0.530965869616921\\
909	0.494080190552381\\
913	0.458889003436639\\
917	0.42539965559854\\
921	0.393609710107853\\
925	0.363507671688353\\
929	0.335073741942665\\
933	0.30828059343537\\
937	0.283094152728609\\
941	0.259474383119255\\
945	0.237376058558425\\
949	0.2167495210374\\
953	0.197541414568718\\
958	0.175439874984363\\
963	0.15535143010004\\
968	0.137157320916799\\
973	0.120736764645358\\
978	0.105968528575204\\
983	0.0927323245842899\\
989	0.0787050493099741\\
995	0.0665160673247556\\
1000	0.0576260334245262\\
};
\addlegendentry{$\bm{x}^{\star}$}

\addplot [color=black, dashed, line width=1.5pt]
  table[row sep=crcr]{%
1	0.123083758447706\\
7	0.123490764692178\\
13	0.125373401013803\\
20	0.129420277301847\\
27	0.135438329716635\\
34	0.143399009801783\\
41	0.153267397887589\\
48	0.165001816270092\\
55	0.178553405998855\\
62	0.193865680925114\\
69	0.210874081167276\\
76	0.229505548412817\\
83	0.249678148141697\\
91	0.274502327708888\\
99	0.301069023632408\\
108	0.33283055210461\\
117	0.366323100154204\\
127	0.405217769707406\\
139	0.453643214811109\\
158	0.532369389206337\\
178	0.614751837601261\\
190	0.66237182438158\\
200	0.700326864308749\\
209	0.73279208767633\\
217	0.760093416076415\\
225	0.785779362043968\\
233	0.809738133606402\\
241	0.831891775837107\\
249	0.852197262493178\\
257	0.870646272859062\\
265	0.88726351615594\\
273	0.902103518567969\\
281	0.915245871565048\\
289	0.926789041086181\\
298	0.938000522950574\\
307	0.947487033719995\\
317	0.956183570840267\\
327	0.963092536078648\\
337	0.968303614206889\\
347	0.971813965026399\\
357	0.973513303310142\\
366	0.973312373601743\\
374	0.97156223834088\\
382	0.968120232015167\\
389	0.963538644161531\\
396	0.9573199268915\\
403	0.949300209399098\\
409	0.940875555766638\\
415	0.930928769780166\\
421	0.919390232933665\\
427	0.906211277252964\\
433	0.891369142355416\\
439	0.874871665540695\\
445	0.856761543165362\\
451	0.837120043711138\\
458	0.812435718612505\\
465	0.786106194718741\\
474	0.750414371113038\\
503	0.633507860453165\\
509	0.611776104996352\\
515	0.591910499578489\\
520	0.577123939260105\\
525	0.564237210541251\\
529	0.555476967871755\\
533	0.548239307436006\\
537	0.542653597479898\\
541	0.538843179856599\\
545	0.536924160999206\\
549	0.537004198478826\\
553	0.539181374005921\\
557	0.543543101496198\\
560	0.548294043018245\\
563	0.554344023048429\\
566	0.561716205652033\\
569	0.570429408109703\\
572	0.580497915122351\\
575	0.591931355435008\\
578	0.60473458828028\\
581	0.618907592294363\\
584	0.634445412279774\\
587	0.65133808213875\\
590	0.669570646594025\\
593	0.689123108260105\\
597	0.71720235581347\\
601	0.74751027491925\\
605	0.779958718530793\\
609	0.814444687716332\\
613	0.85085091694134\\
617	0.889046552171408\\
621	0.928887954898642\\
625	0.970219627562528\\
630	1.02372552217787\\
635	1.07895322658874\\
641	1.14698219066338\\
649	1.23957151697402\\
664	1.41361615455185\\
670	1.48141110929396\\
675	1.53634664035201\\
680	1.5894520184753\\
684	1.63036913003089\\
688	1.66969621377734\\
692	1.7072604391318\\
696	1.74289837940091\\
700	1.77645683880576\\
704	1.8077935923543\\
707	1.82975912565144\\
710	1.85035056402637\\
713	1.86952120472154\\
716	1.88722839060506\\
719	1.90343356918572\\
722	1.91810235618789\\
725	1.931204575214\\
728	1.94271427708668\\
731	1.95260974603116\\
734	1.96087347529988\\
737	1.96749216822059\\
740	1.97245667668449\\
743	1.97576194306032\\
746	1.97740694678339\\
749	1.97739461591812\\
752	1.97573173064802\\
755	1.97242883097988\\
758	1.96750007464436\\
761	1.96096315390855\\
764	1.95283913394928\\
767	1.94315231751693\\
770	1.93193010737218\\
773	1.91920284760772\\
776	1.90500366869912\\
779	1.88936832232696\\
782	1.87233502075753\\
785	1.85394427059362\\
788	1.8342386994243\\
792	1.80599664027807\\
796	1.77560711396677\\
800	1.74318557800871\\
804	1.70885171309112\\
808	1.6727287275105\\
812	1.63494269488376\\
816	1.59562188681809\\
821	1.54451051107355\\
826	1.49145960047281\\
831	1.43672538754913\\
837	1.36918165034376\\
844	1.28836297967121\\
852	1.19416658456475\\
881	0.850817426373055\\
888	0.770592138452002\\
894	0.703461724800604\\
900	0.638136431359385\\
905	0.585265632700043\\
910	0.533968819885558\\
915	0.48437616792819\\
920	0.43660863202058\\
925	0.390778255793066\\
930	0.346988523034838\\
934	0.313490303688354\\
938	0.281404418984152\\
942	0.250772997359036\\
946	0.221635245456582\\
950	0.194027588270501\\
954	0.167983800714978\\
958	0.143535143110512\\
962	0.120710491315322\\
966	0.0995364575109079\\
970	0.0800375149594856\\
974	0.0622361073051252\\
978	0.046152752260241\\
982	0.031806142668529\\
986	0.0192132272412664\\
990	0.00838929416988776\\
994	-0.000651965059660142\\
998	-0.0078983975150777\\
1000	-0.0108451661299114\\
};
\addlegendentry{$\bm{x}^{[5]}_{\gamma}$}

\end{axis}
\end{tikzpicture}%
        \begin{tikzpicture}

\begin{axis}[%
width=6.028in,
height=4.754in,
at={(1.011in,0.642in)},
scale only axis,
xmin=0,
xmax=1000,
xlabel style={font=\color{white!15!black}, scale=1.5},
xlabel={$i = 1, \cdots, 1\,000$},
ymin=0,
ymax=2.05,
axis background/.style={fill=white},
title style={font=\bfseries, scale=1.2},
title={\normalsize\colorbox{black}{\color{white}$\gamma = 10^{10}$}},
xmajorgrids,
ymajorgrids,
grid style={dashed},
legend style={legend cell align=left, align=left, draw=white!15!black, font=\Large},
xtick={1, 100, 200, 300, 400, 500, 600, 700, 800, 900, 1000},
ytick={0, 0.5, 1, 1.5, 2},
grid style = loosely dashed,
]
\addplot [color=green, line width=1.5pt]
  table[row sep=crcr]{%
1	0.101622890399199\\
11	0.116007834887\\
21	0.131907221914275\\
30	0.147572617796186\\
39	0.164571346637558\\
48	0.182942200406956\\
57	0.202714496961789\\
66	0.223906629057751\\
75	0.246524637742255\\
84	0.270560841485576\\
93	0.295992554424629\\
102	0.322780928428301\\
112	0.354068208240392\\
122	0.38685790373097\\
132	0.421018800182196\\
143	0.459987710272799\\
156	0.507582055809053\\
172	0.567748888726555\\
210	0.711392209890278\\
222	0.754914931057783\\
232	0.789773969455268\\
241	0.819749658137994\\
249	0.845076895771285\\
257	0.868988379896336\\
264	0.888621735168726\\
271	0.906942734530389\\
278	0.923852868408517\\
285	0.939259991808285\\
292	0.953079177101472\\
298	0.963601794178885\\
304	0.972857783071731\\
310	0.980808860166803\\
316	0.987422081971204\\
322	0.992670105435991\\
328	0.996531418137579\\
334	0.998990537992313\\
340	1.00003818284074\\
346	0.999671410938845\\
352	0.997893734139666\\
358	0.994715206308683\\
364	0.990152490297987\\
370	0.984228907567967\\
376	0.976974475279917\\
382	0.968425936359267\\
388	0.958626788609308\\
395	0.945681287092157\\
402	0.931193919365228\\
409	0.915269126738735\\
416	0.898024383824122\\
424	0.876868269117949\\
433	0.851483164489423\\
443	0.821718846674457\\
456	0.781437012187894\\
479	0.709821918734292\\
488	0.683281995523998\\
496	0.661162019544122\\
503	0.643340606012316\\
509	0.629479345284949\\
515	0.617164743330022\\
520	0.608255189315855\\
525	0.600725913770702\\
530	0.594721109628608\\
535	0.590386922798075\\
539	0.588221255788881\\
543	0.587294176117894\\
547	0.587680404172829\\
551	0.589453643372394\\
555	0.592686052750878\\
559	0.597447676283082\\
563	0.603805830937745\\
567	0.611824456257182\\
571	0.621563429103276\\
575	0.633077848079438\\
579	0.646417293021386\\
583	0.66162506582566\\
586	0.674279081046734\\
589	0.6880170258014\\
592	0.702848683975276\\
595	0.718781000291187\\
598	0.73581790293656\\
601	0.753960131025451\\
604	0.773205068309039\\
607	0.793546584590104\\
610	0.814974886319987\\
614	0.845212448221559\\
618	0.877314306935205\\
622	0.911224443792548\\
626	0.946873588074595\\
630	0.984178938382229\\
634	1.0230439920399\\
638	1.06335849158484\\
642	1.10499849638359\\
647	1.1587018857449\\
652	1.21396355619595\\
658	1.2818711225882\\
666	1.37408993440863\\
679	1.52409441993791\\
685	1.5916240436643\\
690	1.64628747056236\\
694	1.6886359881845\\
698	1.72950463124516\\
702	1.76866461311761\\
705	1.79677834167012\\
708	1.82371274466607\\
711	1.84937826047883\\
714	1.87368828753199\\
717	1.89655959770596\\
720	1.91791273794615\\
723	1.93767241733383\\
726	1.95576787694529\\
729	1.97213323991991\\
732	1.98670783927207\\
735	1.99943652111915\\
738	2.01026992115601\\
741	2.01916471238928\\
744	2.02608382233723\\
747	2.03099661811973\\
750	2.03387905809052\\
753	2.03471380890778\\
756	2.03349032718779\\
759	2.03020490515496\\
762	2.02486067995812\\
765	2.01746760660296\\
768	2.00804239471279\\
771	1.99660840960246\\
774	1.98319553841156\\
777	1.96784002230299\\
780	1.95058425597358\\
783	1.93147655596692\\
786	1.91057089949129\\
789	1.88792663565755\\
792	1.86360817123386\\
795	1.83768463318518\\
798	1.81022951041007\\
801	1.78132027721347\\
804	1.75103800115392\\
807	1.71946693798452\\
810	1.68669411645681\\
814	1.64128238706826\\
818	1.59411029966543\\
822	1.5453991105461\\
826	1.49537255453231\\
831	1.43133148782363\\
837	1.35284512443252\\
846	1.23329596287897\\
857	1.087277642137\\
863	1.00897986292534\\
868	0.944990532585393\\
873	0.882455818651806\\
877	0.833650181437292\\
881	0.786059350877622\\
885	0.739787911687813\\
889	0.694927238604009\\
893	0.651555592232739\\
897	0.609738330062555\\
901	0.569528223327211\\
905	0.530965869616921\\
909	0.494080190552381\\
913	0.458889003436639\\
917	0.42539965559854\\
921	0.393609710107853\\
925	0.363507671688353\\
929	0.335073741942665\\
933	0.30828059343537\\
937	0.283094152728609\\
941	0.259474383119255\\
945	0.237376058558425\\
949	0.2167495210374\\
953	0.197541414568718\\
957	0.179695389768995\\
961	0.163152773935167\\
965	0.147853202394003\\
969	0.133735207768268\\
974	0.117654802860329\\
979	0.103202543527914\\
984	0.0902584819544927\\
989	0.0787050493099741\\
994	0.0684280545677893\\
1000	0.0576260334245262\\
};
\addlegendentry{$\bm{x}^{\star}$}

\addplot [color=black, dashed, line width=1.5pt]
  table[row sep=crcr]{%
1	0.0722076075072664\\
6	0.0892480195819871\\
11	0.105036492396948\\
17	0.122545792788515\\
23	0.138719021650672\\
30	0.156224424655989\\
38	0.174850427663046\\
49	0.198911838919571\\
78	0.261181686027498\\
89	0.286768014413724\\
98	0.309206820372083\\
106	0.33051207957908\\
113	0.350256665929351\\
121	0.374190803972397\\
128	0.396334752235475\\
136	0.423000693677636\\
144	0.451079391639496\\
153	0.484181250353799\\
169	0.546531523629483\\
181	0.595330426073929\\
216	0.738850178027405\\
227	0.781455819437383\\
236	0.814441750386436\\
244	0.841980865013511\\
250	0.861388778192463\\
256	0.879632917334447\\
262	0.896593795176159\\
268	0.912206054068179\\
274	0.926463250979282\\
280	0.939256717460353\\
286	0.950558814218311\\
292	0.960414467695159\\
298	0.968782519475894\\
304	0.975683954535953\\
310	0.981131466018837\\
317	0.985733619698635\\
323	0.988220432089747\\
329	0.98942702357806\\
335	0.989430898285946\\
342	0.987958183514365\\
348	0.98552400502308\\
355	0.981490156549057\\
366	0.97249685210852\\
374	0.964202010828558\\
384	0.951813945589265\\
392	0.940408346116442\\
401	0.925981011204954\\
410	0.90987686995777\\
417	0.896136801017178\\
424	0.881268628846897\\
430	0.867654323177248\\
438	0.848258463287493\\
449	0.819092806153094\\
456	0.799163793632715\\
468	0.762872929924356\\
476	0.737560951473824\\
501	0.657772903178625\\
507	0.63990206530525\\
514	0.620573050820781\\
519	0.608087766672497\\
524	0.59692376506132\\
528	0.589060032586872\\
532	0.582304163473964\\
537	0.575692043895174\\
541	0.571896746802395\\
544	0.570004233086934\\
548	0.568967226501286\\
551	0.569212708574469\\
555	0.571146089697322\\
559	0.574902397739379\\
562	0.579032429515792\\
565	0.584278514548487\\
568	0.590688897138421\\
571	0.598341230933215\\
574	0.607090429873438\\
578	0.620809151923709\\
581	0.632544669879167\\
584	0.645552324673986\\
587	0.659800359864334\\
590	0.675381774048674\\
594	0.698135312654586\\
597	0.716635894634351\\
600	0.736406587116448\\
603	0.757413295304332\\
606	0.779633345447337\\
610	0.811069272646023\\
613	0.836030066172043\\
617	0.87092452213335\\
621	0.907679920590795\\
625	0.946252213622643\\
629	0.986353628248366\\
633	1.02782192116183\\
638	1.08152317253632\\
643	1.13687559295806\\
648	1.1935476223274\\
656	1.28596179475301\\
671	1.45992390576316\\
677	1.52763905789902\\
682	1.58253022186852\\
686	1.6250894971829\\
690	1.66625929776887\\
694	1.70580916990627\\
698	1.7435795655208\\
701	1.77064998902756\\
704	1.79648267788116\\
707	1.82110301408557\\
710	1.84441037877355\\
713	1.86631961378748\\
716	1.88679597372129\\
719	1.9057912978418\\
722	1.92323109058577\\
725	1.93903852441542\\
728	1.95327459664429\\
731	1.965826955289\\
734	1.9766739827528\\
737	1.985794989429\\
740	1.99319427423745\\
743	1.99879962395312\\
746	2.00262362773708\\
749	2.00464300342207\\
752	2.00488175033604\\
755	2.00333139979602\\
758	1.99996104836589\\
761	1.99484640677179\\
764	1.98796250890564\\
767	1.97931274025188\\
770	1.96895051630588\\
773	1.95686539625865\\
776	1.94311030279198\\
779	1.92772133284961\\
782	1.91072566365938\\
785	1.89215984754969\\
788	1.87208950736158\\
791	1.8505337520769\\
794	1.82754221347545\\
797	1.8031753089856\\
800	1.77749256438767\\
803	1.75056627938261\\
806	1.72238170248977\\
809	1.69308143221053\\
813	1.65232022303746\\
817	1.60977717824585\\
821	1.56559577490532\\
825	1.51993964548046\\
829	1.47295669521429\\
834	1.41262072849679\\
839	1.35079188093471\\
844	1.28778385571013\\
852	1.18527174446876\\
872	0.927329621828903\\
878	0.851547324697435\\
884	0.777360467193148\\
888	0.728964170745598\\
893	0.66994979416188\\
897	0.624035379056522\\
901	0.579441982464573\\
905	0.536210382996501\\
909	0.494486803209611\\
913	0.45439170326847\\
917	0.415983390628639\\
921	0.379366992286805\\
925	0.344652586989014\\
928	0.31989937699268\\
931	0.29629341196221\\
934	0.273859734140842\\
937	0.252636842054017\\
940	0.232658552649013\\
943	0.213952663325813\\
946	0.196531183431944\\
949	0.180442109285991\\
952	0.165710173994512\\
955	0.152350705781487\\
958	0.140400472483179\\
961	0.129886875505917\\
964	0.120815382878277\\
967	0.113217594857019\\
970	0.107116423235766\\
973	0.102537291614908\\
976	0.0994922430105589\\
979	0.0980057182501923\\
982	0.098105592974548\\
985	0.0997848649983553\\
988	0.103087701093614\\
991	0.108018922591555\\
994	0.114598247812069\\
997	0.122852101846775\\
1000	0.132775150786756\\
};
\addlegendentry{$\bm{x}^{[1]}_{\gamma}$}

\end{axis}
\end{tikzpicture}%
} 
        \\
        \scalebox{0.425}{
        \input{ShawSol4.tex}
        \input{ShawSol5.tex}} 
    \caption{\footnotesize The exact solution $\bm{x}^{\star}$ (in green solid line) computed from \eqref{Fredholm2}, and the stabilized solution $\bm{x}^{[k]}_{\gamma}$ (in black dashed line) obtained from \eqref{EqP1}, are presented with respect to the index $i = 1, \cdots, 1\, 000$, for different values of $\gamma = 10^{5}, 10^{10}, 10^{12}, 10^{14}$. Data $\mathbf{A}$ and $\bm{b}$ are computed from the Matlab code \textsf{shaw(n)} in \cite{HansenBook07}, the dimension $n = 1\,000$, stepsize $\alpha = 1$, $k_{\text{max}} = n$ and tolerance $\varepsilon = 10^{-5}$, with $\textbf{cond}(\mathbf{A}) = 7.697e^{+20}$, $\textbf{rank}(\mathbf{A}) = 20$.}
    \label{Figshaw1}
\end{figure}

We also presented in Table \ref{Tabshaw1}, results computed from the stabilized gradient method \eqref{EqP1} in the solving of the linear system \eqref{Fredholm7} under the dimension $n = 1\,000$. Despite the severely ill-conditioning $\textbf{cond}(\mathbf{A}) = 7.697e^{+20}$ and the rank-deficient $\textbf{rank}(\mathbf{A}) = 20$ faced, we obtain more and more accurated results with larger stabilization parameter $\gamma$. Moreover, the stabilized gradient method behaves as a direct method as long as $\gamma > \frac{1}{\varepsilon}$. However, as we mentioned through the SVD expansion at {\bf Remark \ref{RMK1}}, the solution becomes unstable as $\gamma \to +\infty$ since it coincide with the naïve solution, and this is what we observe in Table \ref{Tabshaw1} when $\gamma = 10^{14}$.

We finally present in Figure \ref{Figshaw1}, both the exact solution $\bm{x}^{\star}$ (in green solid line) computed from \eqref{Fredholm2} and the stabilized solution $\bm{x}^{[k]}_{\gamma}$ (presented in the Table \ref{Tabshaw1}), with respect to the mesh index $i = 1, \cdots, 1\, 000$ and for the different values of the parameter $\gamma = 10^{5}, 10^{10}, 10^{12}, 10^{14}$.

\clearpage
$\bullet$ \,{\bf The heat test problem}

The inverse heat equation that we investigate here is also defined by the Fredholm integral equation of the first kind \eqref{Fredholm1}, where the bounds are given by  $a = 0$, $b = s$, $c = 0$ and $d = 1$, and with the convolution kernel type $K(s, t) = \mathcal{K}(s - t)$ defined by 
\begin{align*}
   \mathcal{K}(s - t) = \frac{1}{2\kappa\sqrt{\pi}(s - t)^{3/2}}\exp\left(-\frac{1}{4\kappa^2(s-t)}\right).
\end{align*}
The parameter $\kappa$ controls the ill-conditioning of the linear system built from the discretization of the problem \eqref{Fredholm1}, by means of a quadrature rule (the midpoint quadrature rule described in \eqref{Fredholm5} is used here). Throughout this numerical example, we fix the parameter $\kappa = 1$ corresponding to the worst ill-conditioning system of the heat test problem. For the computing of the data $\mathbf{A}$, $\bm{b}$ and $\bm{x}^{\star}$ from \eqref{Fredholm8}, we use the Matlab code \textsf{heat(n, $\kappa$)} implemented in {\it Regularization Tools} \cite{HansenBook07}. The exact solution $\bm{x}^{\star} = (x^{\star}_{j})_{1 \leq j \leq n}$, where $x^{\star}_{j} = x^{\star}(t^{\star}_j)$ is defined as follows (see \cite{ELDÉN95} and \cite{Hanke95} - Section 6.7)\\
\begin{align}\label{Heat2}
    x^{\star}(t) = \left\{
    \begin{aligned}
        & \frac{300}{4}t^2,             & \text{if } & \; 0 \leq \; t \leq \frac{1}{10},      \\
        & \frac{3}{4}+(20 t-2)(3-20 t), & \text{if } & \; \frac{1}{10} < t \leq \frac{3}{20}, \\
        & \frac{3}{4}e^{2(3-20t)},      & \text{if } & \; \frac{3}{20} < t \leq \frac{1}{2},  \\
        & 0,                            & \text{if } & \; \frac{1}{2} < t \leq 1. 
    \end{aligned}
    \right.
\end{align}

We first remark in Figure \ref{RayonShawHeatGravity}-{\small\colorbox{black}{\color{white}2}}, the bad behavior of the spectral radius $\varrho(\mathbf{M}^{-1}_{\gamma}\mathbf{N}_{k})$ from the stabilized gradient iterative method \eqref{EqP1} applied to this {\bf heat} test problem. 

\begin{figure}[h]  
    \centering  
    \scalebox{0.425}{
        \begin{tikzpicture}

\begin{axis}[%
width=6.028in,
height=4.754in,
at={(1.011in,0.642in)},
scale only axis,
xmin=0,
xmax=1000,
xlabel style={font=\color{white!15!black}, scale=1.5},
xlabel={$i = 1, \cdots, 1\,000$},
ymin=-0.01,
ymax=1.01,
axis background/.style={fill=white},
title style={font=\bfseries, scale=1.2},
title={\normalsize\colorbox{black}{\color{white}$\gamma = 10^{5}$}},
xmajorgrids,
ymajorgrids,
grid style={dashed},
legend style={legend cell align=left, align=left, draw=white!15!black, font=\Large},
xtick={1, 100, 200, 300, 400, 500, 600, 700, 800, 900, 1000},
ytick={0, 0.2, 0.4, 0.6, 0.8, 1},
grid style = loosely dashed,
]
\addplot [color=green, line width=1.5pt]
  table[row sep=crcr]{%
1	7.50000000380169e-05\\
3	0.000675000000001091\\
5	0.00187500000004093\\
7	0.00367500000004384\\
9	0.00607500000000982\\
11	0.00907500000005257\\
13	0.0126749999999447\\
15	0.0168750000000273\\
17	0.0216749999999593\\
19	0.027074999999968\\
21	0.0330750000000535\\
23	0.0396749999999884\\
25	0.046875\\
27	0.0546749999999747\\
29	0.0630750000000262\\
31	0.0720750000000407\\
33	0.0816750000000184\\
35	0.0918749999999591\\
37	0.102674999999977\\
39	0.114074999999957\\
41	0.126075000000014\\
43	0.138675000000035\\
45	0.151875000000018\\
47	0.165674999999965\\
49	0.180074999999988\\
51	0.195074999999974\\
53	0.210675000000037\\
55	0.22687499999995\\
57	0.243675000000053\\
59	0.261075000000005\\
61	0.279075000000034\\
63	0.297675000000027\\
65	0.316874999999982\\
67	0.336675000000014\\
69	0.357075000000009\\
71	0.378074999999967\\
73	0.399675000000002\\
75	0.421875\\
77	0.444674999999961\\
79	0.468074999999999\\
81	0.492075\\
83	0.516674999999964\\
85	0.541875000000005\\
87	0.567675000000008\\
89	0.594074999999975\\
91	0.621075000000019\\
93	0.648675000000026\\
95	0.676874999999995\\
97	0.705675000000042\\
99	0.735075000000052\\
100	0.75\\
101	0.769599999999969\\
102	0.788400000000024\\
103	0.806400000000053\\
104	0.823600000000056\\
105	0.840000000000032\\
106	0.855599999999981\\
107	0.870400000000018\\
108	0.884400000000028\\
109	0.897600000000011\\
110	0.909999999999968\\
111	0.921600000000012\\
112	0.93240000000003\\
113	0.942400000000021\\
114	0.951599999999985\\
115	0.960000000000036\\
116	0.967599999999948\\
117	0.974399999999946\\
118	0.980400000000031\\
119	0.985599999999977\\
120	0.990000000000009\\
121	0.993600000000015\\
122	0.996399999999994\\
123	0.998399999999947\\
124	0.999599999999987\\
125	1\\
126	0.999599999999987\\
127	0.998399999999947\\
128	0.996399999999994\\
129	0.993600000000015\\
130	0.990000000000009\\
131	0.985599999999977\\
132	0.980400000000031\\
133	0.974399999999946\\
134	0.967599999999948\\
135	0.960000000000036\\
136	0.951599999999985\\
137	0.942400000000021\\
138	0.93240000000003\\
139	0.921600000000012\\
140	0.909999999999968\\
141	0.897600000000011\\
142	0.884400000000028\\
143	0.870400000000018\\
144	0.855599999999981\\
145	0.840000000000032\\
146	0.823600000000056\\
147	0.806400000000053\\
148	0.788400000000024\\
149	0.769599999999969\\
150	0.75\\
151	0.720592079364224\\
152	0.692337259789952\\
153	0.665190327537857\\
154	0.639107841724694\\
155	0.614048064808458\\
156	0.589970895799866\\
157	0.566837806091826\\
158	0.544611777805244\\
159	0.523257244553292\\
160	0.50274003452671\\
161	0.483027315812365\\
162	0.464087543854589\\
163	0.445890410977654\\
164	0.428406797886623\\
165	0.411608727070529\\
166	0.395469318032269\\
167	0.379962744274167\\
168	0.365064191969964\\
169	0.35074982025742\\
170	0.336996723087964\\
171	0.323782892571785\\
172	0.311087183761174\\
173	0.298889280813341\\
174	0.287169664481326\\
175	0.275909580878533\\
176	0.265091011469053\\
177	0.254696644233718\\
179	0.235114635661944\\
181	0.217038163454276\\
183	0.200351476474339\\
185	0.184947722956167\\
187	0.170728266287824\\
189	0.157602053400524\\
191	0.145485031718181\\
193	0.134299610933567\\
195	0.123974166166136\\
197	0.114442579317711\\
199	0.105643815690769\\
201	0.0975215331587833\\
203	0.0900237213835453\\
205	0.0831023687717334\\
208	0.0737051892032241\\
211	0.0653706385965052\\
214	0.057978555332511\\
217	0.0514223656157355\\
220	0.0456075469688813\\
223	0.040450265475215\\
227	0.0344694424867384\\
231	0.0293729213242386\\
235	0.0250299524702768\\
240	0.0204927918355224\\
245	0.0167780788921164\\
251	0.0131981043117548\\
258	0.00997491265684403\\
266	0.007243273220638\\
275	0.00505346024931441\\
286	0.00325461245301994\\
300	0.00185906413253178\\
319	0.000869421880452137\\
348	0.000272551744842531\\
408	2.4725336515985e-05\\
846	0\\
1000	0\\
};
\addlegendentry{$\bm{x}^{\star}$}

\addplot [color=black, dashed, line width=1.5pt]
  table[row sep=crcr]{%
1	-0.00144012358271084\\
4	0.000967205929669035\\
7	0.00418425213854334\\
10	0.00831975940695884\\
13	0.0134817276505146\\
16	0.0197741398800417\\
18	0.0246452596326208\\
20	0.0300882271224054\\
22	0.0361268664389627\\
24	0.0427822520588279\\
26	0.0500723866485941\\
28	0.0580119000286459\\
30	0.0666117780072\\
32	0.0758791312166522\\
34	0.0858170200868926\\
36	0.0964243557673399\\
38	0.10769589752158\\
40	0.119622372805225\\
42	0.132190742007765\\
44	0.145384634599509\\
46	0.159184974553796\\
48	0.17357080699071\\
50	0.18852032730706\\
53	0.211955068356474\\
56	0.236546859732698\\
59	0.26225245192677\\
62	0.289059863876105\\
65	0.316997608638644\\
67	0.336285533411001\\
69	0.356141364950417\\
71	0.376606736824556\\
73	0.397729828689535\\
75	0.419562928587197\\
77	0.442159209365627\\
79	0.465568779794467\\
81	0.48983412926566\\
83	0.51498513477145\\
85	0.541033858950982\\
87	0.567969402105405\\
89	0.595753103943025\\
91	0.624314395524152\\
93	0.653547588642596\\
96	0.698334393702112\\
102	0.788708895405648\\
104	0.818106807120103\\
106	0.846627460496507\\
108	0.873933174363742\\
109	0.887021844450942\\
110	0.899675768864654\\
111	0.911850918297432\\
112	0.923503470594028\\
113	0.934590022112957\\
114	0.945067795740101\\
115	0.954894843638272\\
116	0.964030245782965\\
117	0.972434304526132\\
118	0.980068735263558\\
119	0.986896853306462\\
120	0.992883759485039\\
121	0.997996521981463\\
122	1.00220435755716\\
123	1.00547881121838\\
124	1.00779393333414\\
125	1.00912645624135\\
126	1.00945596767394\\
127	1.00876508222802\\
128	1.00703960882595\\
129	1.00426871337118\\
130	1.00044507539712\\
131	0.995565036568337\\
132	0.989628739982322\\
133	0.982640256249397\\
134	0.974607696406565\\
135	0.965543305806136\\
136	0.955463538976915\\
137	0.944389110447332\\
138	0.932345018982232\\
139	0.919360542667619\\
140	0.905469201354208\\
141	0.890708684169681\\
142	0.875120739158433\\
143	0.858751024508024\\
144	0.841648919084832\\
145	0.823867291357374\\
146	0.805462228321062\\
147	0.786492723836432\\
148	0.767020329259253\\
150	0.726823517106368\\
152	0.685399940265256\\
159	0.538510786817142\\
161	0.498093230537393\\
162	0.478424214543679\\
163	0.459177389047227\\
164	0.440394262208429\\
165	0.422111886989569\\
166	0.404362709956899\\
167	0.387174472146853\\
168	0.370570162652257\\
169	0.354568020471675\\
170	0.339181587224857\\
171	0.324419804748231\\
172	0.31028715745083\\
173	0.296783855312924\\
174	0.283906055261127\\
175	0.271646114922873\\
176	0.259992878311778\\
177	0.248931985366085\\
178	0.238446203043395\\
179	0.228515771986054\\
180	0.219118764115478\\
181	0.210231445852742\\
182	0.201828641903148\\
183	0.193884094916257\\
184	0.186370816470685\\
186	0.17252846425481\\
188	0.160083432952547\\
190	0.148825431004184\\
192	0.138559590235673\\
194	0.129112194999493\\
196	0.120334504806124\\
199	0.108164680414347\\
202	0.0969363412219764\\
205	0.0864766169089535\\
208	0.0767143381257256\\
211	0.0676503925080851\\
214	0.0593265781798209\\
217	0.0517971180926224\\
220	0.0451057931118157\\
223	0.0392703388841937\\
226	0.0342745050788835\\
229	0.0300671433432171\\
232	0.0265669521712653\\
236	0.0228205189687287\\
240	0.0198779862188303\\
246	0.0164141649908061\\
256	0.0117182684258523\\
267	0.00724488639161791\\
275	0.00476986980004312\\
283	0.00319269544911549\\
292	0.0023105924834681\\
326	0.000666247569824918\\
345	0.000204051972673369\\
498	3.67106781595794e-06\\
852	0.000260424425505335\\
866	0.000307614273083345\\
895	-0.000443033333112908\\
906	0.00116795059102515\\
913	0.00179931399929956\\
918	0.00148166269082139\\
923	0.000361724932531615\\
932	-0.0028232528100034\\
937	-0.00404344192452299\\
940	-0.00410895393156352\\
943	-0.00347129881720321\\
946	-0.00203259818226798\\
949	0.000216959303656949\\
952	0.00318391454368339\\
957	0.00915175656552947\\
963	0.016289703390612\\
966	0.0191560847152914\\
969	0.0213045241466716\\
972	0.0227122127566872\\
976	0.0236354419424742\\
982	0.0238898474419784\\
993	0.0233891711562819\\
1000	0.022466544291774\\
};
\addlegendentry{$\bm{x}^{[73]}_{\gamma}$}

\end{axis}
\end{tikzpicture}%
        \begin{tikzpicture}

\begin{axis}[%
width=6.028in,
height=4.754in,
at={(1.011in,0.642in)},
scale only axis,
xmin=0,
xmax=1000,
xlabel style={font=\color{white!15!black}, scale=1.5},
xlabel={$i = 1, \cdots, 1\,000$},
ymin=-0.01,
ymax=1.01,
axis background/.style={fill=white},
title style={font=\bfseries, scale=1.2},
title={\normalsize\colorbox{black}{\color{white}$\gamma = 10^{10}$}},
xmajorgrids,
ymajorgrids,
grid style={dashed},
legend style={legend cell align=left, align=left, draw=white!15!black, font=\Large},
xtick={1, 100, 200, 300, 400, 500, 600, 700, 800, 900, 1000},
ytick={0, 0.2, 0.4, 0.6, 0.8, 1},
grid style = loosely dashed,
]
\addplot [color=green, line width=1.5pt]
  table[row sep=crcr]{%
1	7.50000000380169e-05\\
3	0.000675000000001091\\
5	0.00187500000004093\\
7	0.00367500000004384\\
9	0.00607500000000982\\
11	0.00907500000005257\\
13	0.0126749999999447\\
15	0.0168750000000273\\
17	0.0216749999999593\\
19	0.027074999999968\\
21	0.0330750000000535\\
23	0.0396749999999884\\
25	0.046875\\
27	0.0546749999999747\\
29	0.0630750000000262\\
31	0.0720750000000407\\
33	0.0816750000000184\\
35	0.0918749999999591\\
37	0.102674999999977\\
39	0.114074999999957\\
41	0.126075000000014\\
43	0.138675000000035\\
45	0.151875000000018\\
47	0.165674999999965\\
49	0.180074999999988\\
51	0.195074999999974\\
53	0.210675000000037\\
55	0.22687499999995\\
57	0.243675000000053\\
59	0.261075000000005\\
61	0.279075000000034\\
63	0.297675000000027\\
65	0.316874999999982\\
67	0.336675000000014\\
69	0.357075000000009\\
71	0.378074999999967\\
73	0.399675000000002\\
75	0.421875\\
77	0.444674999999961\\
79	0.468074999999999\\
81	0.492075\\
83	0.516674999999964\\
85	0.541875000000005\\
87	0.567675000000008\\
89	0.594074999999975\\
91	0.621075000000019\\
93	0.648675000000026\\
95	0.676874999999995\\
97	0.705675000000042\\
99	0.735075000000052\\
100	0.75\\
101	0.769599999999969\\
102	0.788400000000024\\
103	0.806400000000053\\
104	0.823600000000056\\
105	0.840000000000032\\
106	0.855599999999981\\
107	0.870400000000018\\
108	0.884400000000028\\
109	0.897600000000011\\
110	0.909999999999968\\
111	0.921600000000012\\
112	0.93240000000003\\
113	0.942400000000021\\
114	0.951599999999985\\
115	0.960000000000036\\
116	0.967599999999948\\
117	0.974399999999946\\
118	0.980400000000031\\
119	0.985599999999977\\
120	0.990000000000009\\
121	0.993600000000015\\
122	0.996399999999994\\
123	0.998399999999947\\
124	0.999599999999987\\
125	1\\
126	0.999599999999987\\
127	0.998399999999947\\
128	0.996399999999994\\
129	0.993600000000015\\
130	0.990000000000009\\
131	0.985599999999977\\
132	0.980400000000031\\
133	0.974399999999946\\
134	0.967599999999948\\
135	0.960000000000036\\
136	0.951599999999985\\
137	0.942400000000021\\
138	0.93240000000003\\
139	0.921600000000012\\
140	0.909999999999968\\
141	0.897600000000011\\
142	0.884400000000028\\
143	0.870400000000018\\
144	0.855599999999981\\
145	0.840000000000032\\
146	0.823600000000056\\
147	0.806400000000053\\
148	0.788400000000024\\
149	0.769599999999969\\
150	0.75\\
151	0.720592079364224\\
152	0.692337259789952\\
153	0.665190327537857\\
154	0.639107841724694\\
155	0.614048064808458\\
156	0.589970895799866\\
157	0.566837806091826\\
158	0.544611777805244\\
159	0.523257244553292\\
160	0.50274003452671\\
161	0.483027315812365\\
162	0.464087543854589\\
163	0.445890410977654\\
164	0.428406797886623\\
165	0.411608727070529\\
166	0.395469318032269\\
167	0.379962744274167\\
168	0.365064191969964\\
169	0.35074982025742\\
170	0.336996723087964\\
171	0.323782892571785\\
172	0.311087183761174\\
173	0.298889280813341\\
174	0.287169664481326\\
175	0.275909580878533\\
176	0.265091011469053\\
177	0.254696644233718\\
178	0.244709845967236\\
179	0.235114635661944\\
180	0.225895658934178\\
181	0.217038163454276\\
183	0.200351476474339\\
185	0.184947722956167\\
187	0.170728266287824\\
189	0.157602053400524\\
191	0.145485031718181\\
193	0.134299610933567\\
195	0.123974166166136\\
197	0.114442579317711\\
199	0.105643815690769\\
201	0.0975215331587833\\
203	0.0900237213835453\\
205	0.0831023687717334\\
207	0.0767131550366003\\
209	0.0708151673971997\\
212	0.0628074191942005\\
215	0.0557051836607343\\
218	0.0494060658197668\\
221	0.0438192494733585\\
224	0.0388641878795397\\
227	0.0344694424867384\\
231	0.0293729213242386\\
235	0.0250299524702768\\
239	0.0213291185356184\\
244	0.017462805281184\\
249	0.014297335718652\\
255	0.011246682615365\\
262	0.00850005986603719\\
270	0.00617231028672904\\
279	0.00430627476418977\\
290	0.00277339778733676\\
304	0.00158418995374632\\
323	0.000740872455367025\\
352	0.000232253276521988\\
412	2.10695418445539e-05\\
852	0\\
1000	0\\
};
\addlegendentry{$\bm{x}^{\star}$}

\addplot [color=black, dashed, line width=1.5pt]
  table[row sep=crcr]{%
1	1.04426578673156e-05\\
3	0.000681133289845093\\
5	0.00190580146693264\\
7	0.00370428566407099\\
9	0.00609069006316076\\
11	0.00907481949468547\\
13	0.0126623613119818\\
15	0.0168553750819456\\
17	0.0216544902954183\\
19	0.0270587412807117\\
21	0.0330667031312259\\
23	0.0396766201308765\\
25	0.0468867294658821\\
27	0.0546949431789017\\
29	0.0630992912725787\\
31	0.0720980717781003\\
33	0.0816902711286502\\
35	0.091876390132029\\
37	0.102658355601989\\
39	0.114040401163038\\
41	0.126027839614267\\
43	0.138626196188056\\
45	0.151839455573622\\
47	0.16566832291187\\
49	0.180108576335897\\
51	0.195151023627773\\
53	0.210782565681143\\
55	0.226989555146702\\
57	0.243761105988142\\
59	0.261094966678229\\
61	0.279000323349351\\
63	0.297498733466341\\
65	0.316621393879359\\
67	0.336401978426238\\
69	0.356866645543732\\
71	0.378022259866611\\
73	0.399847867866356\\
75	0.422291884358515\\
77	0.445278945117593\\
80	0.480600056001322\\
83	0.516783459176736\\
85	0.541423058950954\\
87	0.566615145860624\\
89	0.592558189301826\\
91	0.619486939684975\\
93	0.647617698984163\\
95	0.677081729708902\\
97	0.707861529771435\\
99	0.739744357012341\\
105	0.836898968872561\\
106	0.852378581857238\\
107	0.867394024138662\\
108	0.881852301954382\\
109	0.895667938014071\\
110	0.908764500025541\\
111	0.921075617826546\\
112	0.932546352750023\\
113	0.943133392417508\\
114	0.952805355610963\\
115	0.96154152459269\\
116	0.969332194720664\\
117	0.976176463941442\\
118	0.982082169330738\\
119	0.987064109997959\\
120	0.991142619881998\\
121	0.994342607662475\\
122	0.996692011351342\\
123	0.998220771786805\\
124	0.998959589235824\\
125	0.998938779726359\\
126	0.998186706720162\\
127	0.996729215851133\\
128	0.994587648771585\\
129	0.991778250543348\\
130	0.988310700243687\\
131	0.984186825090092\\
132	0.979400028509531\\
133	0.97393410219729\\
134	0.967763628928083\\
135	0.960853416016903\\
136	0.953159369575133\\
137	0.944630184549851\\
138	0.935208513879275\\
139	0.924834254699249\\
140	0.913447336550462\\
141	0.900991338989343\\
142	0.887417050528256\\
143	0.87268718330904\\
144	0.856780029905849\\
145	0.839692640305771\\
146	0.821444076860985\\
147	0.802077433289469\\
148	0.781661223857213\\
149	0.760289082954955\\
150	0.738078177390094\\
151	0.715167226034737\\
152	0.691712744354504\\
154	0.643859590544366\\
156	0.595936877738268\\
157	0.572382118333621\\
158	0.549306794868926\\
159	0.526844574326333\\
160	0.505107811035373\\
161	0.484183689900192\\
162	0.464134645044396\\
163	0.444997859009504\\
164	0.42678603732827\\
165	0.409490917979497\\
166	0.393084900870463\\
167	0.377525651968767\\
168	0.362759082185107\\
169	0.348724329252491\\
170	0.335356480389805\\
171	0.322590873141849\\
172	0.310365668880195\\
173	0.298624298018467\\
174	0.287317419119859\\
175	0.27640373871111\\
177	0.255632972438093\\
179	0.236145926934682\\
181	0.217882040696395\\
183	0.200844256402206\\
185	0.185054607090365\\
187	0.170519436054178\\
189	0.157208359922151\\
191	0.14505040573772\\
193	0.13394209717751\\
195	0.123762924521202\\
197	0.11439309297657\\
199	0.105727821500409\\
201	0.0976846798162114\\
203	0.0902058637747132\\
205	0.0832536675829942\\
207	0.0768039560387024\\
209	0.0708381873463395\\
212	0.062753878294302\\
215	0.0556289667154033\\
218	0.0493536188897679\\
221	0.0438086320436923\\
224	0.0388859993813639\\
228	0.0331482213370009\\
232	0.028231946650294\\
236	0.0240396881268907\\
240	0.0204800307765254\\
245	0.0167758906380868\\
250	0.0137421071493691\\
256	0.0108077763140955\\
263	0.00816443693838664\\
271	0.00593059421100861\\
280	0.00413779985240126\\
291	0.00266443050213638\\
305	0.00152215858111049\\
324	0.000711842128566786\\
354	0.000214380584452556\\
418	1.65868115118428e-05\\
898	3.39958685344754e-05\\
931	-0.000150515007931062\\
937	-0.000562487796855748\\
941	-0.000116731066896136\\
949	0.00136924922196613\\
951	0.00107543096203244\\
954	-0.000211135723134248\\
959	-0.00316339576500013\\
961	-0.00360280790562229\\
963	-0.00303742576033983\\
965	-0.00124257905622471\\
967	0.00177015080191723\\
969	0.00572530809654381\\
974	0.0163334892030207\\
976	0.0195079664074456\\
978	0.0216670062307003\\
980	0.0228901628421454\\
982	0.0234319190735732\\
986	0.0235771960577722\\
1000	0.0231307389640278\\
};
\addlegendentry{$\bm{x}^{[1]}_{\gamma}$}

\end{axis}
\end{tikzpicture}%
} 
        \\
        \scalebox{0.425}{
        \begin{tikzpicture}

\begin{axis}[%
width=6.028in,
height=4.754in,
at={(1.011in,0.642in)},
scale only axis,
xmin=0,
xmax=1000,
xlabel style={font=\color{white!15!black}, scale=1.5},
xlabel={$i = 1, \cdots, 1\,000$},
ymin=-0.01,
ymax=1.01,
axis background/.style={fill=white},
title style={font=\bfseries, scale=1.2},
title={\normalsize\colorbox{black}{\color{white}$\gamma = 10^{12}$}},
xmajorgrids,
ymajorgrids,
grid style={dashed},
legend style={legend cell align=left, align=left, draw=white!15!black, font=\Large},
xtick={1, 100, 200, 300, 400, 500, 600, 700, 800, 900, 1000},
ytick={0, 0.2, 0.4, 0.6, 0.8, 1},
grid style = loosely dashed,
]
\addplot [color=green, line width=1.5pt]
  table[row sep=crcr]{%
1	7.50000000380169e-05\\
3	0.000675000000001091\\
5	0.00187500000004093\\
7	0.00367500000004384\\
9	0.00607500000000982\\
11	0.00907500000005257\\
13	0.0126749999999447\\
15	0.0168750000000273\\
17	0.0216749999999593\\
19	0.027074999999968\\
21	0.0330750000000535\\
23	0.0396749999999884\\
25	0.046875\\
27	0.0546749999999747\\
29	0.0630750000000262\\
31	0.0720750000000407\\
33	0.0816750000000184\\
35	0.0918749999999591\\
37	0.102674999999977\\
39	0.114074999999957\\
41	0.126075000000014\\
43	0.138675000000035\\
45	0.151875000000018\\
47	0.165674999999965\\
49	0.180074999999988\\
51	0.195074999999974\\
53	0.210675000000037\\
55	0.22687499999995\\
57	0.243675000000053\\
59	0.261075000000005\\
61	0.279075000000034\\
63	0.297675000000027\\
65	0.316874999999982\\
67	0.336675000000014\\
69	0.357075000000009\\
71	0.378074999999967\\
73	0.399675000000002\\
75	0.421875\\
77	0.444674999999961\\
79	0.468074999999999\\
81	0.492075\\
83	0.516674999999964\\
85	0.541875000000005\\
87	0.567675000000008\\
89	0.594074999999975\\
91	0.621075000000019\\
93	0.648675000000026\\
95	0.676874999999995\\
97	0.705675000000042\\
99	0.735075000000052\\
100	0.75\\
101	0.769599999999969\\
102	0.788400000000024\\
103	0.806400000000053\\
104	0.823600000000056\\
105	0.840000000000032\\
106	0.855599999999981\\
107	0.870400000000018\\
108	0.884400000000028\\
109	0.897600000000011\\
110	0.909999999999968\\
111	0.921600000000012\\
112	0.93240000000003\\
113	0.942400000000021\\
114	0.951599999999985\\
115	0.960000000000036\\
116	0.967599999999948\\
117	0.974399999999946\\
118	0.980400000000031\\
119	0.985599999999977\\
120	0.990000000000009\\
121	0.993600000000015\\
122	0.996399999999994\\
123	0.998399999999947\\
124	0.999599999999987\\
125	1\\
126	0.999599999999987\\
127	0.998399999999947\\
128	0.996399999999994\\
129	0.993600000000015\\
130	0.990000000000009\\
131	0.985599999999977\\
132	0.980400000000031\\
133	0.974399999999946\\
134	0.967599999999948\\
135	0.960000000000036\\
136	0.951599999999985\\
137	0.942400000000021\\
138	0.93240000000003\\
139	0.921600000000012\\
140	0.909999999999968\\
141	0.897600000000011\\
142	0.884400000000028\\
143	0.870400000000018\\
144	0.855599999999981\\
145	0.840000000000032\\
146	0.823600000000056\\
147	0.806400000000053\\
148	0.788400000000024\\
149	0.769599999999969\\
150	0.75\\
151	0.720592079364224\\
152	0.692337259789952\\
153	0.665190327537857\\
154	0.639107841724694\\
155	0.614048064808458\\
156	0.589970895799866\\
157	0.566837806091826\\
158	0.544611777805244\\
159	0.523257244553292\\
160	0.50274003452671\\
161	0.483027315812365\\
162	0.464087543854589\\
163	0.445890410977654\\
164	0.428406797886623\\
165	0.411608727070529\\
166	0.395469318032269\\
167	0.379962744274167\\
168	0.365064191969964\\
169	0.35074982025742\\
170	0.336996723087964\\
171	0.323782892571785\\
172	0.311087183761174\\
173	0.298889280813341\\
174	0.287169664481326\\
175	0.275909580878533\\
176	0.265091011469053\\
177	0.254696644233718\\
179	0.235114635661944\\
181	0.217038163454276\\
183	0.200351476474339\\
185	0.184947722956167\\
187	0.170728266287824\\
189	0.157602053400524\\
191	0.145485031718181\\
193	0.134299610933567\\
195	0.123974166166136\\
197	0.114442579317711\\
199	0.105643815690769\\
201	0.0975215331587833\\
203	0.0900237213835453\\
205	0.0831023687717334\\
208	0.0737051892032241\\
211	0.0653706385965052\\
214	0.057978555332511\\
217	0.0514223656157355\\
220	0.0456075469688813\\
223	0.040450265475215\\
227	0.0344694424867384\\
231	0.0293729213242386\\
235	0.0250299524702768\\
240	0.0204927918355224\\
245	0.0167780788921164\\
251	0.0131981043117548\\
258	0.00997491265684403\\
266	0.007243273220638\\
275	0.00505346024931441\\
286	0.00325461245301994\\
300	0.00185906413253178\\
319	0.000869421880452137\\
348	0.000272551744842531\\
408	2.4725336515985e-05\\
846	0\\
1000	0\\
};
\addlegendentry{$\bm{x}^{\star}$}

\addplot [color=black, dashed, line width=1.5pt]
  table[row sep=crcr]{%
1	6.99877940633087e-05\\
3	0.000674526075272297\\
5	0.00187992243888857\\
7	0.00367354600518865\\
9	0.00607490940967637\\
11	0.00907180137971864\\
13	0.0126708338276558\\
15	0.0168809464182686\\
17	0.0216767219899339\\
19	0.0270730531170784\\
21	0.0330776833775417\\
23	0.0396730966378982\\
25	0.0468751303300223\\
27	0.0546684084997651\\
29	0.0630752918949611\\
31	0.072071075680924\\
33	0.0816754390240249\\
35	0.0918825543320736\\
37	0.102677682591548\\
39	0.114081538760388\\
41	0.126072450853599\\
43	0.138660848886161\\
45	0.151861374121154\\
47	0.165669283128182\\
49	0.180087507337248\\
51	0.195096911653309\\
53	0.210704001082263\\
55	0.226889140448293\\
57	0.243668081573105\\
59	0.261028564794174\\
61	0.279007884426051\\
63	0.297639117008089\\
65	0.316900425376843\\
67	0.336755076581312\\
69	0.357198123406874\\
72	0.388847582181825\\
74	0.410611637438933\\
76	0.43298404616155\\
78	0.456052382077019\\
80	0.479867884029545\\
82	0.504433933757809\\
85	0.542391982603704\\
89	0.594142796622918\\
91	0.620490653254251\\
93	0.647544226453874\\
95	0.675822170552124\\
96	0.690563919893179\\
97	0.705792221983188\\
98	0.721483441916689\\
100	0.754139372149325\\
105	0.838197396695136\\
106	0.854191879861105\\
107	0.869550086931213\\
108	0.88413399184526\\
109	0.897843436422249\\
110	0.910632733212651\\
111	0.922462083441928\\
112	0.933308746253033\\
113	0.943216757282357\\
114	0.952205261333688\\
115	0.960329280981682\\
116	0.967635045021893\\
117	0.974146303245448\\
118	0.979934291352379\\
119	0.985011962706494\\
120	0.989383356922076\\
121	0.99304987685673\\
122	0.99601518237148\\
123	0.998235790464378\\
124	0.999686816680878\\
125	1.00033275057785\\
126	1.00014098249108\\
127	0.99906306938226\\
128	0.997080290699273\\
129	0.994193333038197\\
130	0.990394124327167\\
131	0.985720990826849\\
132	0.980173846130469\\
133	0.973822993435874\\
134	0.966719462248989\\
135	0.958901692753557\\
136	0.950458308550651\\
137	0.941395702342334\\
138	0.931715686109555\\
139	0.921432584624085\\
140	0.910481825570628\\
141	0.898792892725055\\
142	0.886245177155956\\
143	0.872705685677261\\
144	0.858045093429155\\
145	0.842126447039846\\
146	0.82484920363413\\
147	0.806149461027758\\
148	0.786016022175659\\
149	0.76450481262475\\
150	0.741750933706953\\
151	0.717937880980571\\
152	0.693304957159285\\
156	0.592705846093168\\
157	0.568524715292824\\
158	0.545198700051628\\
159	0.522878692444806\\
160	0.5016540088684\\
161	0.481535316128088\\
162	0.462486853585347\\
163	0.444434943831652\\
164	0.427277162803762\\
165	0.410916058322186\\
166	0.395246328460985\\
167	0.38014738748393\\
168	0.365559267161075\\
169	0.351451988414397\\
170	0.337759960707785\\
171	0.324500934199818\\
172	0.311660532198175\\
173	0.299262196835002\\
174	0.287332430206902\\
175	0.275869521468508\\
176	0.264872479898258\\
177	0.254380405798202\\
178	0.244341555183155\\
179	0.234755424589139\\
181	0.216824827862297\\
183	0.200342319006495\\
185	0.185087310647759\\
187	0.170911776309822\\
189	0.157719320964929\\
191	0.145510390680897\\
193	0.134241313706525\\
195	0.12388268254449\\
197	0.114380519980045\\
199	0.105623693919142\\
201	0.0975455734709385\\
203	0.0900610666192279\\
206	0.0798746590003248\\
209	0.0708063629516573\\
212	0.0627874957363019\\
215	0.0556971521775722\\
218	0.0494132613540614\\
221	0.0438307244337466\\
224	0.0388688998411908\\
228	0.0331104640615649\\
232	0.0282156575940462\\
237	0.0231071822936428\\
242	0.0189186452068952\\
248	0.0148842962432809\\
254	0.0116996392093824\\
261	0.00884990880319947\\
269	0.00642030635276569\\
279	0.00430728725712015\\
290	0.00277836728821512\\
303	0.00164976510586712\\
324	0.000705997043610296\\
388	5.45831080671633e-05\\
636	-1.01023852039361e-06\\
940	0.000117836974709462\\
945	0.000325655615256437\\
954	-0.000533497610376799\\
958	0.00111193881582494\\
960	0.00166389520165922\\
962	0.00132511753827202\\
964	-7.27326776086556e-05\\
968	-0.00367687856146404\\
969	-0.00393110819936737\\
970	-0.0036457778762724\\
971	-0.00274057877300038\\
972	-0.00119529739890822\\
973	0.000943835951261462\\
975	0.00651382015064428\\
978	0.0153793903388078\\
979	0.0177765332025501\\
980	0.0197195057556883\\
981	0.0211891893275151\\
983	0.0228796536860045\\
985	0.0234496495297662\\
990	0.0234550857599061\\
1000	0.0231307389640278\\
};
\addlegendentry{$\bm{x}^{[1]}_{\gamma}$}

\end{axis}
\end{tikzpicture}%
        \input{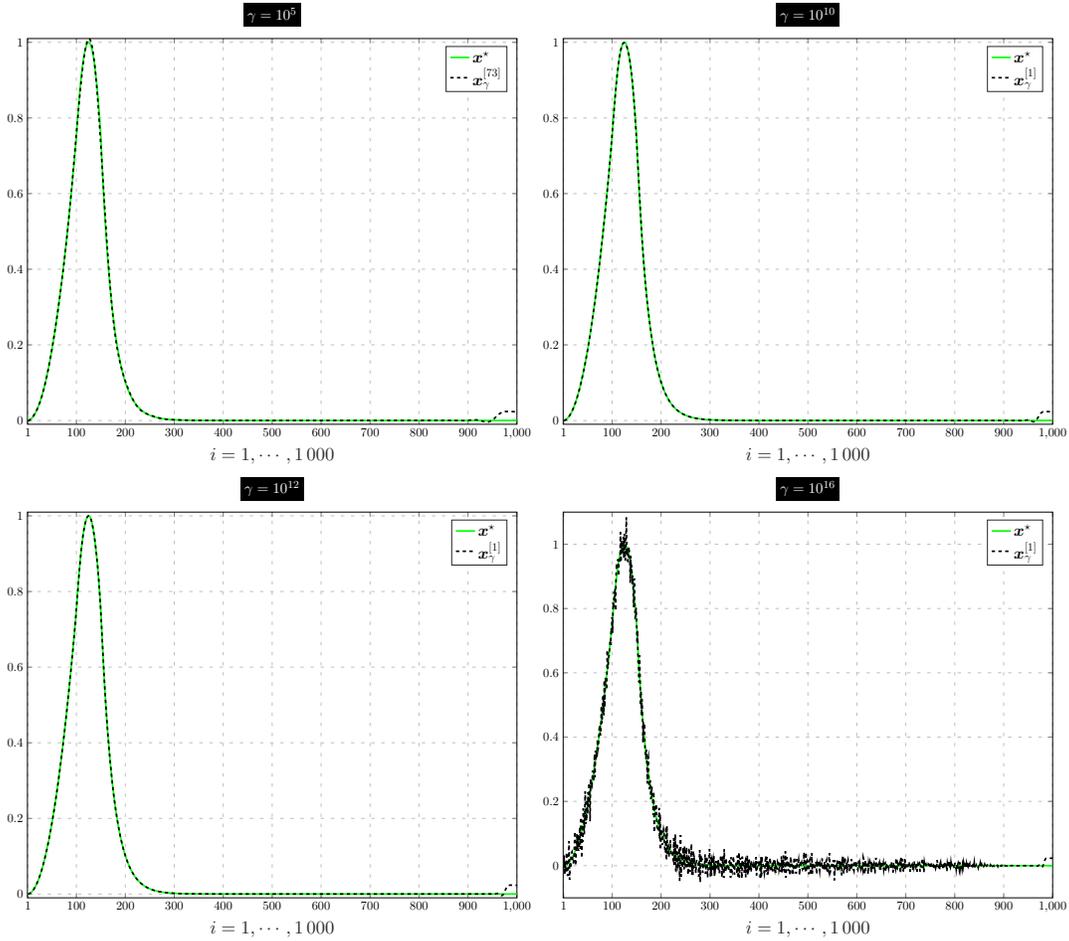}} 
    \caption{\footnotesize The exact solution $\bm{x}^{\star}$ (in green solid line) computed from \eqref{Heat2}, and the stabilized solution $\bm{x}^{[k]}_{\gamma}$ (in black dashed line) obtained from \eqref{EqP1}, are presented with respect to the index $i = 1, \cdots, 1\, 000$, for different values of $\gamma = 10^{5}, 10^{10}, 10^{12}, 10^{16}$. Data $\mathbf{A}$ and $\bm{b}$ are computed from the Matlab code \textsf{heat($n$,$1$)} in \cite{HansenBook07}, the dimension $n = 1\,000$, stepsize $\alpha = 1$, $k_{\text{max}} = n$ and tolerance $\varepsilon = 10^{-5}$, with $\textbf{cond}(\mathbf{A}) = 2.0332e^{+232}$, $\textbf{rank}(\mathbf{A}) = 588$.}
    \label{Figheat1}
\end{figure}

We also presented in Table \ref{Tabheat1}, results computed from the stabilized gradient method \eqref{EqP1} in the solving of this problem formulated under the linear system \eqref{Fredholm7} and the dimension $n = 1\,000$. We again present in Figure \ref{Figheat1}, both the exact solution $\bm{x}^{\star}$ (in green solid line) computed from \eqref{Heat2} and the stabilized solution $\bm{x}^{[k]}_{\gamma}$ (presented in Table \ref{Tabheat1}), with respect to the mesh index $i = 1, \cdots, 1\, 000$, and for the different values of the parameter $\gamma = 10^{5}, 10^{10}, 10^{12}, 10^{16}$. We can make the same observation than the previous example, despite the special ill-conditioning $\textbf{cond}(\mathbf{A}) = 2.0332e^{+232}$ faced and the rank-deficient $\textbf{rank}(\mathbf{A}) = 588$. However, we note a kind of boss at the end of the domain in the Figure \ref{Figheat1}.
\begin{table}[h]
\begin{center}
\begin{tikzpicture}
\draw[gray] [white, fill=gray!15.5]
    (0.0,6).. controls (15.5,6) ..(15.5,6)
    .. controls (15.5, 7.25) ..(15.5, 7.25)
    .. controls (0.0, 7.25) ..(0.0, 7.25)
    .. controls (0.0,6) ..(0.0,6);

\draw[gray] [white, fill=gray!10]
    (0.0,4.0).. controls (15.5,4.0) ..(15.5,4.0)
    .. controls (15.5, 5) .. (15.5, 5)
    .. controls (0.0, 5)  .. (0.0, 5)
    .. controls (0.0,4.0) .. (0.0,4.0);
\draw[gray] [white, fill=gray!10]
    (0.0,2.0).. controls (15.5,2.0) ..(15.5,2.0)
    .. controls (15.5, 3) .. (15.5, 3)
    .. controls (0.0, 3)  .. (0.0, 3)
    .. controls (0.0,2.0) .. (0.0,2.0);
\draw[gray] [white, fill=gray!10]
    (0.0,0.0).. controls (15.5,0.0) ..(15.5,0.0)
    .. controls (15.5, 1) .. (15.5, 1)
    .. controls (0.0, 1)  .. (0.0, 1)
    .. controls (0.0,0.0) .. (0.0,0.0);

\draw[, very thin] (1.5,-1) -- (1.5,7.25);
\draw[, very thin] (4.5,-1) -- (4.5,7.25);
\draw[, very thin] (7.5,-1) -- (7.5,7.25);
\draw[, very thin] (10.5,-1) -- (10.5,7.25);
\draw[, very thin] (13.5,-1) -- (13.5,7.25);

\draw[, very thin] (0,-1) -- (15.5,-1);
\draw[, very thin] (0,0) -- (15.5,0);
\draw[, very thin] (0,1) -- (15.5,1);
\draw[, very thin] (0,2)   -- (15.5,2);
\draw[, very thin] (0,3)   -- (15.5,3);
\draw[, very thin] (0,4)   -- (15.5,4);
\draw[, very thin] (0,5) -- (15.5,5);

\draw[black](0.56, 7.)  node [rotate=40, scale=0.9] {\text{Stab.}};
\draw[black](0.7, 6.65)  node [rotate=40, scale=0.9] {\text{Parameter}};
\draw[black](1., 6.3)  node [rotate=40, scale=1] {$(\gamma)$};
\draw[black](0.75, 5.5)  node [rotate=0, scale=1] {$10^{3}$};
\draw[black](0.75, 4.5)  node [rotate=0, scale=1] {$10^{4}$};
\draw[black](0.75, 3.5)  node [rotate=0, scale=1] {$10^{5}$};
\draw[black](0.75, 2.5)  node [rotate=0, scale=1] {$10^{6}$};
\draw[black](0.75, 1.5)  node [rotate=0, scale=1] {$10^{10}$};
\draw[black](0.75, 0.5)  node [rotate=0, scale=1] {$10^{12}$};
\draw[black](0.75, -0.5)  node [rotate=0, scale=1] {$10^{16}$};

\draw[black](3, 6.7)   node [rotate=0, scale=1] {$\|\mathbf{A}\bm{x}^{[k]}_{\gamma} - \bm{b}\|_{n}$};
\draw[black](3, 5.5)  node [rotate=0, scale=1] {$5.009976e^{-05}$};
\draw[black](3, 4.5)  node [rotate=0, scale=1] {$1.476816e^{-05}$};
\draw[black](3, 3.5)  node [rotate=0, scale=1] {$1.475814e^{-05}$};
\draw[black](3, 2.5)  node [rotate=0, scale=1] {$1.465276e^{-05}$};
\draw[black](3, 1.5)  node [rotate=0, scale=1] {$3.223123e^{-07}$};
\draw[black](3, 0.5)  node [rotate=0, scale=1] {$2.337872e^{-08}$};
\draw[black](3, -0.5)  node [rotate=0, scale=1] {$8.569897e^{-10}$};

\draw[black](6, 6.6)   node [rotate=0, scale=1] {$\frac{\|\mathbf{A}\bm{x}^{[k]}_{\gamma} - \bm{b}\|_{n}}{\|\mathbf{A}\bm{x}^{[0]}_{\gamma} - \bm{b}\|_{n}}$};
\draw[black](6, 5.5)  node [rotate=0, scale=1] {$3.390948e^{-05}$};
\draw[black](6, 4.5)  node [rotate=0, scale=1] {$9.995670e^{-06}$};
\draw[black](6, 3.5)  node [rotate=0, scale=1] {$9.988888e^{-06}$};
\draw[black](6, 2.5)  node [rotate=0, scale=1] {$9.917561e^{-06}$};
\draw[black](6, 1.5)  node [rotate=0, scale=1] {$2.181536e^{-07}$};
\draw[black](6, 0.5)  node [rotate=0, scale=1] {$1.582363e^{-08}$};
\draw[black](6, -0.5)  node [rotate=0, scale=1] {$5.800442e^{-10}$};

\draw[black](9., 6.7)   node [rotate=0, scale=1] {$\|\bm{x}^{\star} - \bm{x}^{[k]}_{\gamma}\|_{n}$};
\draw[black](9, 5.5)  node [rotate=0, scale=1] {$2.030031e^{-01}$};
\draw[black](9, 4.5)  node [rotate=0, scale=1] {$1.745498e^{-01}$};
\draw[black](9, 3.5)  node [rotate=0, scale=1] {$1.718651e^{-01}$};
\draw[black](9, 2.5)  node [rotate=0, scale=1] {$1.711090e^{-01}$};
\draw[black](9, 1.5)  node [rotate=0, scale=1] {$1.236628e^{-01}$};
\draw[black](9, 0.5)  node [rotate=0, scale=1] {$1.109604e^{-01}$};
\draw[black](9, -0.5)  node [rotate=0, scale=1] {$6.524621e^{-01}$};

\draw[black](12., 6.6)   node [rotate=0, scale=1] {$\frac{\|\bm{x}^{\star} - \bm{x}^{[k]}_{\gamma}\|_{n}}{\|\bm{x}^{\star}\|_{n}}$};
\draw[black](12, 5.5)  node [rotate=0, scale=1] {$2.608322e^{-02}$};
\draw[black](12, 4.5)  node [rotate=0, scale=1] {$2.242734e^{-02}$};
\draw[black](12, 3.5)  node [rotate=0, scale=1] {$2.208240e^{-02}$};
\draw[black](12, 2.5)  node [rotate=0, scale=1] {$2.198525e^{-02}$};
\draw[black](12, 1.5)  node [rotate=0, scale=1] {$1.588904e^{-02}$};
\draw[black](12, 0.5)  node [rotate=0, scale=1] {$1.425694e^{-02}$};
\draw[black](12, -0.5)  node [rotate=0, scale=1] {$8.383277e^{-02}$};

\draw[black](14.5, 6.85)   node [rotate=0, scale=1] {\text{$\#$ Iterations}};
\draw[black](14.5, 6.3)   node [rotate=0, scale=1] {$[k]$};
\draw[black](14.5, 5.5)  node [rotate=0, scale=1] {$1\,000$};
\draw[black](14.5, 4.5)  node [rotate=0, scale=1] {$730$};
\draw[black](14.5, 3.5)  node [rotate=0, scale=1] {$73$};
\draw[black](14.5, 2.5)  node [rotate=0, scale=1] {$8$};
\draw[black](14.5, 1.5)  node [rotate=0, scale=1] {$1$};
\draw[black](14.5, 0.5)  node [rotate=0, scale=1] {$1$};
\draw[black](14.5, -0.5)  node [rotate=0, scale=1] {$1$};

\end{tikzpicture} 
\caption{Numerical results from the stabilized method \eqref{EqP1} of the {\bf heat} problem for different values of the stabilization parameter $\gamma$, under the fixed dimension $n=1\,000$, stepsize $\alpha = 1$, maximum iterations number $k_{\text{max}} = n$ and stopping tolerance $\varepsilon = 10^{-5}$ in \eqref{CTR1}. The data $\mathbf{A}$, $\bm{b}$ and exact solution $\bm{x}^{\star}$ are computed from the Matlab code \textsf{heat($n$,$1$)} provided in {\it Regularization Tools} \cite{HansenBook07}, where $\textbf{cond}(\mathbf{A}) = 2.0332e^{+232}$, $\textbf{rank}(\mathbf{A}) = 588$.}
\label{Tabheat1}
\end{center}
\end{table}

$\bullet$ \,{\bf The gravity surveying test problem}

The gravity problem finally investigated here is also defined under the integral \eqref{Fredholm1}, with the bounds given by  $a = 0$, $b = 1$, $c = 0$ and $d = 1$, an analytical solution, and kernel given by
\begin{align}
    x^{\star}(t) & = sin(\pi t) + \frac{1}{2}sin(2\pi t), \; 0 \leq t \leq 1, \label{GravSol}\\
    K(s, t) & = \frac{1}{\left(1 + \left(s - t\right)^{2}\right)^{\frac{3}{2}}}, \;\, 0 \leq s \leq 1, \; 0 \leq t \leq 1. \nonumber
\end{align} 
The problem is discretized by means of the midpoint quadrature rule with $n$ points, where the integration in the interval $[0,1]$ leads to a symmetric Toeplitz matrix $\mathbf{A}$. The right-hand side is computed as follows $\bm{b} = \mathbf{A}\bm{x}^{\star}$, where the $n$-vector exact solution $\bm{x}^{\star} = (x^{\star}_{j})_{1 \leq j \leq n}$ is derived from the analytical formula in \eqref{GravSol} by $x^{\star}_{j} = x^{\star}(x^{\star}_j), \;\, i = 1, \cdots, n$.  For the computing of the data $\mathbf{A}$ and $\bm{x}^{\star}$, the Matlab code \textsf{gravity(n,1,0,1,1)} implemented in {\it Regularization Tools} \cite{HansenBook07} is used with the dimension $n = 1\, 000$.

The same results than the previous two test problems are computed and presented in Figure \ref{RayonShawHeatGravity}-{\small\colorbox{black}{\color{white}3}}, in Table \ref{Tabgravity1} and in Figure \ref{Figgravity1}. Despite the behavior of the spectral radius $\varrho(\mathbf{M}^{-1}_{\gamma}\mathbf{N}_{k})$, the ill-conditioning $\textbf{cond}(\mathbf{A}) = 1.5823e^{+20}$, and the very severe rank-deficient $\textbf{rank}(\mathbf{A}) = 15$, we observe similar results than the previous {\bf shaw} test problem.
\begin{table}[h]
\begin{center}
\begin{tikzpicture}
\draw[gray] [white, fill=gray!15.5]
    (0.0,6).. controls (15.5,6) ..(15.5,6)
    .. controls (15.5, 7.25) ..(15.5, 7.25)
    .. controls (0.0, 7.25) ..(0.0, 7.25)
    .. controls (0.0,6) ..(0.0,6);

\draw[gray] [white, fill=gray!10]
    (0.0,4.0).. controls (15.5,4.0) ..(15.5,4.0)
    .. controls (15.5, 5) .. (15.5, 5)
    .. controls (0.0, 5)  .. (0.0, 5)
    .. controls (0.0,4.0) .. (0.0,4.0);
\draw[gray] [white, fill=gray!10]
    (0.0,2.0).. controls (15.5,2.0) ..(15.5,2.0)
    .. controls (15.5, 3) .. (15.5, 3)
    .. controls (0.0, 3)  .. (0.0, 3)
    .. controls (0.0,2.0) .. (0.0,2.0);
\draw[gray] [white, fill=gray!10]
    (0.0,0.0).. controls (15.5,0.0) ..(15.5,0.0)
    .. controls (15.5, 1) .. (15.5, 1)
    .. controls (0.0, 1)  .. (0.0, 1)
    .. controls (0.0,0.0) .. (0.0,0.0);

\draw[, very thin] (1.5,-1) -- (1.5,7.25);
\draw[, very thin] (4.5,-1) -- (4.5,7.25);
\draw[, very thin] (7.5,-1) -- (7.5,7.25);
\draw[, very thin] (10.5,-1) -- (10.5,7.25);
\draw[, very thin] (13.5,-1) -- (13.5,7.25);

\draw[, very thin] (0,-1) -- (15.5,-1);
\draw[, very thin] (0,0) -- (15.5,0);
\draw[, very thin] (0,1) -- (15.5,1);
\draw[, very thin] (0,2)   -- (15.5,2);
\draw[, very thin] (0,3)   -- (15.5,3);
\draw[, very thin] (0,4)   -- (15.5,4);
\draw[, very thin] (0,5) -- (15.5,5);

\draw[black](0.46, 6.9)  node [rotate=40, scale=0.9] {\text{Stab.}};
\draw[black](0.7, 6.65)  node [rotate=40, scale=0.9] {\text{Parameter}};
\draw[black](1., 6.3)  node [rotate=40, scale=1] {$(\gamma)$};
\draw[black](0.75, 5.5)  node [rotate=0, scale=1] {$10^{3}$};
\draw[black](0.75, 4.5)  node [rotate=0, scale=1] {$10^{4}$};
\draw[black](0.75, 3.5)  node [rotate=0, scale=1] {$10^{5}$};
\draw[black](0.75, 2.5)  node [rotate=0, scale=1] {$10^{6}$};
\draw[black](0.75, 1.5)  node [rotate=0, scale=1] {$10^{10}$};
\draw[black](0.75, 0.5)  node [rotate=0, scale=1] {$10^{12}$};
\draw[black](0.75, -0.5)  node [rotate=0, scale=1] {$10^{15}$};

\draw[black](3, 6.7)   node [rotate=0, scale=1] {$\|\mathbf{A}\bm{x}^{[k]}_{\gamma} - \bm{b}\|_{n}$};
\draw[black](3, 5.5)  node [rotate=0, scale=1] {$1.732721e^{-04}$};
\draw[black](3, 4.5)  node [rotate=0, scale=1] {$1.732523e^{-04}$};
\draw[black](3, 3.5)  node [rotate=0, scale=1] {$1.727488e^{-04}$};
\draw[black](3, 2.5)  node [rotate=0, scale=1] {$1.612163e^{-04}$};
\draw[black](3, 1.5)  node [rotate=0, scale=1] {$3.793328e^{-07}$};
\draw[black](3, 0.5)  node [rotate=0, scale=1] {$1.489148e^{-08}$};
\draw[black](3, -0.5)  node [rotate=0, scale=1] {$2.309456e^{-09}$};

\draw[black](6, 6.6)   node [rotate=0, scale=1] {$\frac{\|\mathbf{A}\bm{x}^{[k]}_{\gamma} - \bm{b}\|_{n}}{\|\mathbf{A}\bm{x}^{[0]}_{\gamma} - \bm{b}\|_{n}}$};
\draw[black](6, 5.5)  node [rotate=0, scale=1] {$9.992355e^{-06}$};
\draw[black](6, 4.5)  node [rotate=0, scale=1] {$9.991215e^{-06}$};
\draw[black](6, 3.5)  node [rotate=0, scale=1] {$9.962179e^{-06}$};
\draw[black](6, 2.5)  node [rotate=0, scale=1] {$9.297114e^{-06}$};
\draw[black](6, 1.5)  node [rotate=0, scale=1] {$2.187558e^{-08}$};
\draw[black](6, 0.5)  node [rotate=0, scale=1] {$8.587706e^{-10}$};
\draw[black](6, -0.5)  node [rotate=0, scale=1] {$1.331830e^{-10}$};

\draw[black](9., 6.7)   node [rotate=0, scale=1] {$\|\bm{x}^{\star} - \bm{x}^{[k]}_{\gamma}\|_{n}$};
\draw[black](9, 5.5)  node [rotate=0, scale=1] {$4.920645e^{-01}$};
\draw[black](9, 4.5)  node [rotate=0, scale=1] {$5.038643e^{-01}$};
\draw[black](9, 3.5)  node [rotate=0, scale=1] {$5.024578e^{-01}$};
\draw[black](9, 2.5)  node [rotate=0, scale=1] {$4.631275e^{-01}$};
\draw[black](9, 1.5)  node [rotate=0, scale=1] {$4.844738e^{-02}$};
\draw[black](9, 0.5)  node [rotate=0, scale=1] {$1.817421e^{-02}$};
\draw[black](9, -0.5)  node [rotate=0, scale=1] {$1.933920e^{+00}$};

\draw[black](12., 6.6)   node [rotate=0, scale=1] {$\frac{\|\bm{x}^{\star} - \bm{x}^{[k]}_{\gamma}\|_{n}}{\|\bm{x}^{\star}\|_{n}}$};
\draw[black](12, 5.5)  node [rotate=0, scale=1] {$1.968258e^{-02}$};
\draw[black](12, 4.5)  node [rotate=0, scale=1] {$2.015457e^{-02}$};
\draw[black](12, 3.5)  node [rotate=0, scale=1] {$2.009831e^{-02}$};
\draw[black](12, 2.5)  node [rotate=0, scale=1] {$1.852510e^{-02}$};
\draw[black](12, 1.5)  node [rotate=0, scale=1] {$1.937895e^{-03}$};
\draw[black](12, 0.5)  node [rotate=0, scale=1] {$7.269684e^{-04}$};
\draw[black](12, -0.5)  node [rotate=0, scale=1] {$7.735678e^{-02}$};

\draw[black](14.5, 6.85)   node [rotate=0, scale=1] {\text{$\#$ Iterations}};
\draw[black](14.5, 6.3)   node [rotate=0, scale=1] {$[k]$};
\draw[black](14.5, 5.5)  node [rotate=0, scale=1] {$600$};
\draw[black](14.5, 4.5)  node [rotate=0, scale=1] {$165$};
\draw[black](14.5, 3.5)  node [rotate=0, scale=1] {$21$};
\draw[black](14.5, 2.5)  node [rotate=0, scale=1] {$3$};
\draw[black](14.5, 1.5)  node [rotate=0, scale=1] {$1$};
\draw[black](14.5, 0.5)  node [rotate=0, scale=1] {$1$};
\draw[black](14.5, -0.5)  node [rotate=0, scale=1] {$1$};

\end{tikzpicture} 
\caption{Numerical results from the stabilized method \eqref{EqP1} of the {\bf gravity} problem for different values of the stabilization parameter $\gamma$, under the fixed dimension $n=1\,000$, stepsize $\alpha = 1$, maximum iterations number $k_{\text{max}} = n$ and stopping tolerance $\varepsilon = 10^{-5}$ in \eqref{CTR1}. The data $\mathbf{A}$, $\bm{b}$ and exact solution $\bm{x}^{\star}$ are computed from the Matlab code \textsf{gravity(n,1,0,1,1)} provided in {\it Regularization Tools} \cite{HansenBook07}, where $\textbf{cond}(\mathbf{A}) = 1.5823e^{+20}$, $\textbf{rank}(\mathbf{A}) = 15$.}
\label{Tabgravity1}
\end{center}
\end{table}
\begin{figure}[h]  
    \centering  
    \scalebox{0.425}{
        \begin{tikzpicture}

\begin{axis}[%
width=6.028in,
height=4.754in,
at={(1.011in,0.642in)},
scale only axis,
xmin=0,
xmax=1000,
xlabel style={font=\color{white!15!black}, scale=1.5},
xlabel={$i = 1, \cdots, 1\,000$},
ymin=0,
ymax=1.4,
axis background/.style={fill=white},
title style={font=\bfseries, scale=1.2},
title={\normalsize\colorbox{black}{\color{white}$\gamma = 10^{5}$}},
xmajorgrids,
ymajorgrids,
grid style={dashed},
legend style={legend cell align=left, align=left, draw=white!15!black, font=\Large},
xtick={1, 100, 200, 300, 400, 500, 600, 700, 800, 900, 1000},
ytick={0, 0.2, 0.4, 0.6, 0.8, 1, 1.2, 1.4},
grid style = loosely dashed,
]
\addplot [color=green, line width=1.5pt]
  table[row sep=crcr]{%
1	0.00314158942376253\\
25	0.153558437359379\\
38	0.234260079385763\\
49	0.301798326798576\\
59	0.362423029556794\\
68	0.416228939123812\\
77	0.469208815704746\\
85	0.515525005474842\\
93	0.561035961726816\\
101	0.605672503286769\\
108	0.643959216630378\\
115	0.681481047420789\\
122	0.718195402020797\\
129	0.754061001637297\\
136	0.789037947782163\\
142	0.81828170161009\\
148	0.846820949914672\\
154	0.874633191014823\\
160	0.901696787094238\\
166	0.927990988594615\\
172	0.953495957411064\\
178	0.978192788854471\\
184	1.0020635323566\\
190	1.02509121088644\\
196	1.04725983905598\\
202	1.06855443988854\\
208	1.08896106023087\\
214	1.1084667847872\\
220	1.12705974875792\\
226	1.14472914906753\\
232	1.16146525416696\\
238	1.17725941239951\\
244	1.19210405891863\\
250	1.2059927211518\\
255	1.21683242882057\\
260	1.22700197709435\\
265	1.23649936794868\\
270	1.24532313376994\\
275	1.25347233685284\\
280	1.26094656840962\\
285	1.267745947089\\
290	1.27387111701046\\
295	1.27932324531366\\
300	1.2841040192277\\
305	1.2882156426615\\
311	1.29227017311769\\
317	1.29537072580047\\
323	1.29752400286179\\
329	1.29873773512202\\
335	1.29902066724947\\
341	1.29838254161473\\
347	1.29683408084145\\
353	1.29438696908153\\
359	1.29105383204183\\
365	1.28684821579191\\
371	1.28178456438457\\
377	1.27587819632129\\
383	1.269145279898\\
389	1.2616028074666\\
395	1.2532685686499\\
401	1.24416112254926\\
407	1.23429976898501\\
413	1.22370451881159\\
419	1.21239606335007\\
426	1.19832995797481\\
433	1.18335711706732\\
440	1.16751355623956\\
447	1.1508362894341\\
454	1.1333632542545\\
461	1.1151332358404\\
468	1.09618578944003\\
476	1.07370473891308\\
484	1.05040027353118\\
492	1.02633432660662\\
501	0.998427972556783\\
510	0.969727256463784\\
520	0.937016228269385\\
531	0.900179745830314\\
543	0.859163222383245\\
557	0.810493025921346\\
577	0.740075167088776\\
609	0.627399645429477\\
624	0.575467045537494\\
636	0.534663278727749\\
647	0.497996274541492\\
658	0.462159617606972\\
668	0.430397683301067\\
678	0.399496427101212\\
687	0.372482461248637\\
696	0.346275341952719\\
705	0.320920573485523\\
713	0.299131970681856\\
721	0.278075490277274\\
729	0.257774379035595\\
737	0.238249054807625\\
745	0.219517074448163\\
753	0.201593109529995\\
761	0.184488929902614\\
769	0.168213395124667\\
777	0.152772453778653\\
785	0.13816915065479\\
793	0.124403641773142\\
801	0.111473217190792\\
809	0.0993723315225452\\
817	0.088092642082529\\
825	0.0776230545349108\\
833	0.0679497759242622\\
841	0.0590563749337889\\
849	0.0509238492046507\\
857	0.0435306995299243\\
866	0.0360673231375586\\
875	0.0294718214415752\\
884	0.0237030687911783\\
893	0.0187165602947061\\
903	0.0140353339228341\\
913	0.0101911675070596\\
923	0.00711015370279711\\
934	0.00450962662182519\\
946	0.002491287999419\\
959	0.00110336092905072\\
975	0.000256651175618572\\
997	6.64676917949691e-07\\
1000	1.93790583580267e-09\\
};
\addlegendentry{$\bm{x}^{\star}$}

\addplot [color=black, dashed, line width=1.5pt]
  table[row sep=crcr]{%
1	0.0501778457777391\\
53	0.340757872035624\\
66	0.41220154842722\\
77	0.471863263061323\\
87	0.52532181373806\\
96	0.572691658585768\\
105	0.619262365602935\\
113	0.659913793715191\\
121	0.699798045410375\\
129	0.738852313292227\\
136	0.772295285046539\\
143	0.805014654025172\\
150	0.836970404459635\\
157	0.868123355658213\\
164	0.898435252083459\\
171	0.927868852033384\\
178	0.95638801629309\\
184	0.980078670783996\\
190	1.00305052749991\\
196	1.02528304303257\\
202	1.04675643838243\\
208	1.06745173179286\\
214	1.08735077685003\\
220	1.10643629560775\\
226	1.12469191047683\\
232	1.14210217649452\\
238	1.15865260976591\\
244	1.17432971782466\\
250	1.18912102380409\\
256	1.20301509324327\\
262	1.21600155569968\\
268	1.22807112844043\\
274	1.23921563283193\\
280	1.24942801442864\\
286	1.25870235656816\\
292	1.26703389555837\\
298	1.27441903052068\\
304	1.2808553340567\\
310	1.28634155846657\\
316	1.29087764111966\\
322	1.29446470717517\\
328	1.29710507035315\\
334	1.29880222960878\\
340	1.29956086903076\\
346	1.2993868482572\\
352	1.29828719507077\\
358	1.29627009553258\\
364	1.29334488114387\\
370	1.28952201447476\\
376	1.2848130702971\\
382	1.27923071882276\\
388	1.27278870333305\\
394	1.26550181753487\\
400	1.25738588070953\\
406	1.24845770899526\\
412	1.23873508975214\\
418	1.22823674812776\\
424	1.21698231666039\\
430	1.20499229949371\\
436	1.19228803791873\\
443	1.17659325742704\\
450	1.15999294617291\\
457	1.14252489454907\\
464	1.12422812944624\\
471	1.10514282210238\\
478	1.08531019391432\\
486	1.06178338772065\\
494	1.03739991463647\\
502	1.01222535044315\\
511	0.983041376325446\\
520	0.953036723939135\\
530	0.918853949931645\\
541	0.880381563967489\\
553	0.837575028015522\\
568	0.783183471778557\\
592	0.695168597039356\\
614	0.614741362810264\\
628	0.564304879245356\\
640	0.521823154034109\\
651	0.483659120737002\\
661	0.449728562387918\\
671	0.41662556737856\\
680	0.38761549456342\\
689	0.359411087914395\\
698	0.332069993456912\\
706	0.308534785580832\\
714	0.285758486794634\\
722	0.263772889981851\\
730	0.242606981659492\\
738	0.222286923814636\\
746	0.202836044721835\\
754	0.184274840924786\\
762	0.16662098475922\\
769	0.151929927732681\\
776	0.137953203915913\\
783	0.124697423725934\\
790	0.112167613104589\\
797	0.100367240067499\\
804	0.08929824356062\\
811	0.0789610669528429\\
818	0.0693546943709862\\
825	0.0604766886415291\\
832	0.052323235701806\\
839	0.0448891857283797\\
846	0.0381681022452085\\
853	0.0321523099521528\\
861	0.0261293176648678\\
869	0.0210008671854212\\
877	0.016750086692582\\
885	0.0133583818169427\\
893	0.0108055455505109\\
901	0.0090698714619748\\
909	0.0081282676118235\\
917	0.0079563710369257\\
925	0.00852866286834342\\
934	0.0100289230978206\\
943	0.0123981996925977\\
952	0.0155960075643407\\
961	0.0195807835502819\\
970	0.0243100833071139\\
980	0.0303855912577546\\
990	0.0372671313965611\\
1000	0.0448943768890331\\
};
\addlegendentry{$\bm{x}^{[21]}_{\gamma}$}

\end{axis}
\end{tikzpicture}%
        \begin{tikzpicture}

\begin{axis}[%
width=6.028in,
height=4.754in,
at={(1.011in,0.642in)},
scale only axis,
xmin=0,
xmax=1000,
xlabel style={font=\color{white!15!black}, scale=1.5},
xlabel={$i = 1, \cdots, 1\,000$},
ymin=0,
ymax=1.4,
axis background/.style={fill=white},
title style={font=\bfseries, scale=1.2},
title={\normalsize\colorbox{black}{\color{white}$\gamma = 10^{10}$}},
xmajorgrids,
ymajorgrids,
grid style={dashed},
legend style={legend cell align=left, align=left, draw=white!15!black, font=\Large},
xtick={1, 100, 200, 300, 400, 500, 600, 700, 800, 900, 1000},
ytick={0, 0.2, 0.4, 0.6, 0.8, 1, 1.2, 1.4},
grid style = loosely dashed,
]
\addplot [color=green, line width=1.5pt]
  table[row sep=crcr]{%
1	0.00314158942376253\\
25	0.153558437359379\\
38	0.234260079385763\\
49	0.301798326798576\\
59	0.362423029556794\\
68	0.416228939123812\\
77	0.469208815704746\\
85	0.515525005474842\\
93	0.561035961726816\\
101	0.605672503286769\\
108	0.643959216630378\\
115	0.681481047420789\\
122	0.718195402020797\\
129	0.754061001637297\\
136	0.789037947782163\\
142	0.81828170161009\\
148	0.846820949914672\\
154	0.874633191014823\\
160	0.901696787094238\\
166	0.927990988594615\\
172	0.953495957411064\\
178	0.978192788854471\\
184	1.0020635323566\\
190	1.02509121088644\\
196	1.04725983905598\\
202	1.06855443988854\\
208	1.08896106023087\\
214	1.1084667847872\\
220	1.12705974875792\\
226	1.14472914906753\\
232	1.16146525416696\\
238	1.17725941239951\\
244	1.19210405891863\\
250	1.2059927211518\\
255	1.21683242882057\\
260	1.22700197709435\\
265	1.23649936794868\\
270	1.24532313376994\\
275	1.25347233685284\\
280	1.26094656840962\\
285	1.267745947089\\
290	1.27387111701046\\
295	1.27932324531366\\
300	1.2841040192277\\
305	1.2882156426615\\
311	1.29227017311769\\
317	1.29537072580047\\
323	1.29752400286179\\
329	1.29873773512202\\
335	1.29902066724947\\
341	1.29838254161473\\
347	1.29683408084145\\
353	1.29438696908153\\
359	1.29105383204183\\
365	1.28684821579191\\
371	1.28178456438457\\
377	1.27587819632129\\
383	1.269145279898\\
389	1.2616028074666\\
395	1.2532685686499\\
401	1.24416112254926\\
407	1.23429976898501\\
413	1.22370451881159\\
419	1.21239606335007\\
426	1.19832995797481\\
433	1.18335711706732\\
440	1.16751355623956\\
447	1.1508362894341\\
454	1.1333632542545\\
461	1.1151332358404\\
468	1.09618578944003\\
476	1.07370473891308\\
484	1.05040027353118\\
492	1.02633432660662\\
501	0.998427972556783\\
510	0.969727256463784\\
520	0.937016228269385\\
531	0.900179745830314\\
543	0.859163222383245\\
557	0.810493025921346\\
577	0.740075167088776\\
609	0.627399645429477\\
624	0.575467045537494\\
636	0.534663278727749\\
647	0.497996274541492\\
658	0.462159617606972\\
668	0.430397683301067\\
678	0.399496427101212\\
687	0.372482461248637\\
696	0.346275341952719\\
705	0.320920573485523\\
713	0.299131970681856\\
721	0.278075490277274\\
729	0.257774379035595\\
737	0.238249054807625\\
745	0.219517074448163\\
753	0.201593109529995\\
761	0.184488929902614\\
769	0.168213395124667\\
777	0.152772453778653\\
785	0.13816915065479\\
793	0.124403641773142\\
801	0.111473217190792\\
809	0.0993723315225452\\
817	0.088092642082529\\
825	0.0776230545349108\\
833	0.0679497759242622\\
841	0.0590563749337889\\
849	0.0509238492046507\\
857	0.0435306995299243\\
866	0.0360673231375586\\
875	0.0294718214415752\\
884	0.0237030687911783\\
893	0.0187165602947061\\
903	0.0140353339228341\\
913	0.0101911675070596\\
923	0.00711015370279711\\
934	0.00450962662182519\\
946	0.002491287999419\\
959	0.00110336092905072\\
975	0.000256651175618572\\
997	6.64676917949691e-07\\
1000	1.93790583580267e-09\\
};
\addlegendentry{$\bm{x}^{\star}$}

\addplot [color=black, dashed, line width=1.5pt]
  table[row sep=crcr]{%
1	0.00493903909568871\\
21	0.129926611414362\\
35	0.216620765939297\\
47	0.290085457308123\\
57	0.350531742547105\\
66	0.404213749522341\\
75	0.457111872411929\\
83	0.503395063114567\\
91	0.548915180183144\\
99	0.593604608019746\\
106	0.631975892158948\\
113	0.669616890666362\\
120	0.706487146085806\\
127	0.742544599027156\\
134	0.777747678996434\\
141	0.812057858749995\\
147	0.840730381072376\\
153	0.868693064973741\\
159	0.895927628120489\\
165	0.92241031862568\\
171	0.948119128302665\\
177	0.973038325986067\\
183	0.99714608780971\\
189	1.02042270688901\\
195	1.04285150400494\\
201	1.06441610008289\\
207	1.08510342800139\\
213	1.10489641977995\\
219	1.12378103413437\\
225	1.14174569713714\\
231	1.15878122068136\\
237	1.17487397763728\\
243	1.19001720443418\\
248	1.20190471534818\\
253	1.21312388413969\\
258	1.22367098583106\\
263	1.23354154132289\\
268	1.24273377992665\\
273	1.25124777481506\\
278	1.25907913960384\\
283	1.26622932107\\
288	1.2726991201049\\
293	1.27848787259211\\
299	1.28453603726007\\
304	1.28883251642867\\
310	1.29309616680644\\
316	1.29639141287294\\
322	1.29872794368987\\
328	1.30010571640776\\
333	1.30053506462048\\
339	1.30019052661953\\
345	1.29892017341672\\
351	1.29673462695018\\
357	1.29364434109868\\
363	1.28966733714401\\
369	1.28481566402047\\
375	1.27910368937694\\
381	1.27255128243223\\
387	1.26517459678132\\
393	1.25698944211604\\
399	1.24802092269385\\
405	1.2382827048599\\
411	1.22780058121805\\
417	1.21659093730182\\
423	1.20467985227663\\
430	1.18992513517173\\
437	1.17428410457262\\
444	1.15779283437291\\
451	1.14049033315575\\
458	1.12241964348016\\
465	1.10361716154887\\
473	1.08129340276946\\
481	1.05813363308562\\
489	1.03420254648245\\
498	1.00643945504225\\
507	0.977877571364502\\
517	0.945312786993895\\
528	0.908633774112786\\
540	0.867790715375463\\
554	0.819324108879641\\
574	0.749202078179451\\
605	0.640475345470691\\
620	0.588687940915179\\
633	0.544607636495016\\
645	0.504765521410832\\
656	0.469089880925026\\
666	0.437449803578147\\
676	0.406647324017968\\
685	0.379692990641615\\
694	0.353516425774615\\
703	0.328158326377775\\
711	0.30633117879438\\
720	0.282616095100252\\
728	0.262300909022883\\
736	0.242727900199384\\
744	0.223910810179632\\
752	0.205864480282798\\
760	0.188597120648524\\
768	0.172120662735438\\
776	0.156437243199093\\
784	0.14155486322511\\
792	0.127469704219948\\
800	0.114189332789692\\
808	0.101702547154218\\
816	0.0900111777042412\\
824	0.0791060844954927\\
832	0.0689773290355333\\
840	0.0596179699973618\\
848	0.0510160677712292\\
856	0.0431561609603932\\
864	0.0360257257872263\\
873	0.0288548234439077\\
882	0.0225639628534964\\
891	0.0171247804836412\\
900	0.0125093987610398\\
909	0.00868747589788654\\
918	0.00563174103695019\\
928	0.00309460201413003\\
938	0.00141807931299809\\
948	0.000555572129769644\\
958	0.000462668734940053\\
969	0.00119097138804136\\
980	0.00273096081184576\\
992	0.00526098137686404\\
1000	0.00740183781806536\\
};
\addlegendentry{$\bm{x}^{[1]}_{\gamma}$}

\end{axis}
\end{tikzpicture}%
} 
        \\
        \scalebox{0.425}{
        \begin{tikzpicture}

\begin{axis}[%
width=6.028in,
height=4.754in,
at={(1.011in,0.642in)},
scale only axis,
xmin=0,
xmax=1000,
xlabel style={font=\color{white!15!black}, scale=1.5},
xlabel={$i = 1, \cdots, 1\,000$},
ymin=0,
ymax=1.4,
axis background/.style={fill=white},
title style={font=\bfseries, scale=1.2},
title={\normalsize\colorbox{black}{\color{white}$\gamma = 10^{12}$}},
xmajorgrids,
ymajorgrids,
grid style={dashed},
legend style={legend cell align=left, align=left, draw=white!15!black, font=\Large},
xtick={1, 100, 200, 300, 400, 500, 600, 700, 800, 900, 1000},
ytick={0, 0.2, 0.4, 0.6, 0.8, 1, 1.2, 1.4},
grid style = loosely dashed,
]
\addplot [color=green, line width=1.5pt]
  table[row sep=crcr]{%
1	0.00314158942376253\\
25	0.153558437359379\\
38	0.234260079385763\\
49	0.301798326798576\\
59	0.362423029556794\\
68	0.416228939123812\\
77	0.469208815704746\\
85	0.515525005474842\\
93	0.561035961726816\\
101	0.605672503286769\\
108	0.643959216630378\\
115	0.681481047420789\\
122	0.718195402020797\\
129	0.754061001637297\\
136	0.789037947782163\\
142	0.81828170161009\\
148	0.846820949914672\\
154	0.874633191014823\\
160	0.901696787094238\\
166	0.927990988594615\\
172	0.953495957411064\\
178	0.978192788854471\\
184	1.0020635323566\\
190	1.02509121088644\\
196	1.04725983905598\\
202	1.06855443988854\\
208	1.08896106023087\\
214	1.1084667847872\\
220	1.12705974875792\\
226	1.14472914906753\\
232	1.16146525416696\\
238	1.17725941239951\\
244	1.19210405891863\\
250	1.2059927211518\\
255	1.21683242882057\\
260	1.22700197709435\\
265	1.23649936794868\\
270	1.24532313376994\\
275	1.25347233685284\\
280	1.26094656840962\\
285	1.267745947089\\
290	1.27387111701046\\
295	1.27932324531366\\
300	1.2841040192277\\
305	1.2882156426615\\
311	1.29227017311769\\
317	1.29537072580047\\
323	1.29752400286179\\
329	1.29873773512202\\
335	1.29902066724947\\
341	1.29838254161473\\
347	1.29683408084145\\
353	1.29438696908153\\
359	1.29105383204183\\
365	1.28684821579191\\
371	1.28178456438457\\
377	1.27587819632129\\
383	1.269145279898\\
389	1.2616028074666\\
395	1.2532685686499\\
401	1.24416112254926\\
407	1.23429976898501\\
413	1.22370451881159\\
419	1.21239606335007\\
426	1.19832995797481\\
433	1.18335711706732\\
440	1.16751355623956\\
447	1.1508362894341\\
454	1.1333632542545\\
461	1.1151332358404\\
468	1.09618578944003\\
476	1.07370473891308\\
484	1.05040027353118\\
492	1.02633432660662\\
501	0.998427972556783\\
510	0.969727256463784\\
520	0.937016228269385\\
531	0.900179745830314\\
543	0.859163222383245\\
557	0.810493025921346\\
577	0.740075167088776\\
609	0.627399645429477\\
624	0.575467045537494\\
636	0.534663278727749\\
647	0.497996274541492\\
658	0.462159617606972\\
668	0.430397683301067\\
678	0.399496427101212\\
687	0.372482461248637\\
696	0.346275341952719\\
705	0.320920573485523\\
713	0.299131970681856\\
721	0.278075490277274\\
729	0.257774379035595\\
737	0.238249054807625\\
745	0.219517074448163\\
753	0.201593109529995\\
761	0.184488929902614\\
769	0.168213395124667\\
777	0.152772453778653\\
785	0.13816915065479\\
793	0.124403641773142\\
801	0.111473217190792\\
809	0.0993723315225452\\
817	0.088092642082529\\
825	0.0776230545349108\\
833	0.0679497759242622\\
841	0.0590563749337889\\
849	0.0509238492046507\\
857	0.0435306995299243\\
866	0.0360673231375586\\
875	0.0294718214415752\\
884	0.0237030687911783\\
893	0.0187165602947061\\
903	0.0140353339228341\\
913	0.0101911675070596\\
923	0.00711015370279711\\
934	0.00450962662182519\\
946	0.002491287999419\\
959	0.00110336092905072\\
975	0.000256651175618572\\
997	6.64676917949691e-07\\
1000	1.93790583580267e-09\\
};
\addlegendentry{$\bm{x}^{\star}$}

\addplot [color=black, dashed, line width=1.5pt]
  table[row sep=crcr]{%
1	0.00147185921548498\\
5	0.0269875302219589\\
7	0.0396347478086909\\
13	0.0776422569548458\\
40	0.24690467747655\\
51	0.314515248777525\\
56	0.344912116888622\\
63	0.387191616579003\\
67	0.411007942470405\\
72	0.440516588882815\\
79	0.481474297767249\\
84	0.510427787545723\\
97	0.584123848963145\\
101	0.606102703503211\\
104	0.62267549154592\\
114	0.676624173859636\\
116	0.687113213610928\\
121	0.713264662985239\\
130	0.759328558259767\\
140	0.808670260139706\\
153	0.870071356806989\\
156	0.883543968451022\\
159	0.897195650653885\\
162	0.910384344641784\\
167	0.932095419075267\\
177	0.973902088498676\\
181	0.98986213355613\\
187	1.01337968854534\\
192	1.03210361660479\\
194	1.03959106785555\\
204	1.0750079968991\\
209	1.09184511758917\\
211	1.09842952294514\\
217	1.11744797543952\\
222	1.13253587660313\\
226	1.14429032566647\\
235	1.16903375551612\\
243	1.18948724973689\\
246	1.1963742867249\\
249	1.20330515482112\\
258	1.22261472826881\\
265	1.23608624058988\\
271	1.24674014057393\\
273	1.25003498460683\\
276	1.2547438563405\\
285	1.26754026458627\\
296	1.2801885959301\\
302	1.28576504325781\\
311	1.29219136046572\\
317	1.2953378025245\\
328	1.29872292161724\\
335	1.29914700565541\\
341	1.2985039036439\\
349	1.29647246127706\\
353	1.29466323668601\\
358	1.29200290436734\\
367	1.28550648945679\\
372	1.28125696894949\\
386	1.26599722211813\\
389	1.26209092603779\\
393	1.25669582892056\\
397	1.25070876194377\\
399	1.24770228074044\\
412	1.22609954072232\\
418	1.21476297833851\\
420	1.21097494775904\\
425	1.20080660177723\\
430	1.19036120521162\\
436	1.17714544559158\\
441	1.16567591752664\\
456	1.12864836242102\\
464	1.10747740301326\\
473	1.08249521032326\\
485	1.04764734958565\\
494	1.02044473911985\\
511	0.966508968958806\\
519	0.940256813096653\\
532	0.896618992357617\\
536	0.882940164245497\\
541	0.865935816058141\\
545	0.851837888510886\\
549	0.838163380820106\\
552	0.827547113034939\\
554	0.820617668758587\\
558	0.80656424342294\\
566	0.778441174508771\\
585	0.711355204359734\\
588	0.700625629412343\\
599	0.661849526072388\\
606	0.637327337400052\\
609	0.626731572423751\\
611	0.619907324887549\\
615	0.605881673741123\\
619	0.592161843520557\\
623	0.578346101992338\\
632	0.547646682504023\\
654	0.474588111574008\\
670	0.423873992920676\\
674	0.41156431412594\\
680	0.393139943893061\\
684	0.381383921040651\\
688	0.369376226184841\\
692	0.35778204224755\\
698	0.340637256801188\\
705	0.320969523748886\\
707	0.315537473433324\\
713	0.299335836044861\\
723	0.273238299748414\\
732	0.250766014299302\\
740	0.231713900948535\\
745	0.220056079876031\\
760	0.187114150919683\\
763	0.181061960018724\\
770	0.166840218423886\\
773	0.161287717944219\\
776	0.155472520806143\\
785	0.138901202634429\\
794	0.1234348870953\\
798	0.117006921713596\\
801	0.11228009886679\\
812	0.0956543106065055\\
814	0.0930443370968987\\
817	0.088705324591956\\
819	0.0860630223678527\\
822	0.0821074553947483\\
834	0.067222683833279\\
842	0.0582285898314012\\
844	0.0562713571639506\\
848	0.052036926714095\\
856	0.0444253723263728\\
866	0.0359002456876851\\
876	0.0284231618012427\\
883	0.0238285881009688\\
894	0.0174448236184617\\
913	0.00915458451152062\\
925	0.00538864203429057\\
945	0.00157795071902456\\
949	0.00115711555781672\\
956	0.000544593110930691\\
981	0.000895267478313144\\
989	0.00161135406574431\\
992	0.00187547266318688\\
997	0.00260328371518881\\
1000	0.00300492560745624\\
};
\addlegendentry{$\bm{x}^{[1]}_{\gamma}$}

\end{axis}
\end{tikzpicture}%
        \input{GravitySol5.tex}} 
    \caption{\footnotesize The exact solution $\bm{x}^{\star}$ (in green solid line) computed from \eqref{GravSol}, and the stabilized solution $\bm{x}^{[k]}_{\gamma}$ (in black dashed line) obtained from \eqref{EqP1}, are presented with respect to the index $i = 1, \cdots, 1\, 000$, for different values of $\gamma = 10^{5}, 10^{10}, 10^{12}, 10^{15}$. Data $\mathbf{A}$ and $\bm{b}$ are computed from the Matlab code \textsf{gravity(n,1,0,1,1)} in \cite{HansenBook07}, the dimension $n = 1\,000$, stepsize $\alpha = 1$, $k_{\text{max}} = n$ and tolerance $\varepsilon = 10^{-5}$, with $\textbf{cond}(\mathbf{A}) = 1.5823e^{+20}$, $\textbf{rank}(\mathbf{A}) = 15$.}
    \label{Figgravity1}
\end{figure}

\clearpage
\subsubsection{\bf Large-scale finite element problems}\label{LSP1}
We now consider the following reaction-diffusion partial differential equation (see \cite{leary09}, Chap.-28)  
\begin{align}
    -\Delta\bm{u} + \kappa^{2}\bm{u} & = \bm{f}, \text{ in } \Omega,               \label{Equat1}\\
                              \bm{u} & = \bm{0}, \text{ on } \partial\Omega,       \label{Equat2}
\end{align}
posed in the square domain $\Omega = [-1, \; 1]\times[-1, \; 1] \subset \mathbb{R}^{2}$ with boundary $\partial\Omega$, where $\kappa$ is the reaction parameter which is chosen equal to $1$. The right-hand side $\bm{f}$ is computed from the equation \eqref{Equat1} using the following exact solution  
\begin{align}\label{Appro3}
    \bm{u}^{\star}(x,y) = \cos\left(\frac{\pi}{2}x\right)\sin\left(4\pi y\right) - (x^2 - 1)(y^2 - 1),
\end{align}
which satisfies the Dirichlet boundary conditions \eqref{Equat2}.

We apply $P_{1}$ finite elements (i.e. continuous piecewise linear polynomials) in the approximation of the exact solution \eqref{Appro3} through the weak formulation of the boundary value problem \eqref{Equat1} - \eqref{Equat2}. Denoting the mesh size by $h$, the made error is of order $\mathcal{O}(h)$. Therefore, it is sufficient to reach accuracy of the same order in the computing of the approximate solution $\bm{x}^{[k]}_{\gamma,h}$. We use the stabilized gradient method \eqref{EqP1} to solve the finite element equation starting with $\bm{x}^{[0]}_{\gamma} = (0, 0)$, and stopping when the relative error \eqref{CTR2} falls below the length of the longest edge ($\varepsilon = h_{i}$) among the triangles in each \textsf{mesh$_{i}$}, $i = 1, 2, 3, \cdots$. For the discretization, we start with a triangular grid (or \textsf{mesh$_{1}$}) of size $h_{1} = \sqrt{2}$, and we apply regular refinements to this starting triangulation. The computations performed here are with respect to the different meshes \textsf{mesh$_{i}$}, $i = 4, 5, 6, 7, 8$, and whose sizes are $h_{4} = \frac{\sqrt{2}}{2^{3}}$, $h_{5} = \frac{\sqrt{2}}{2^{4}}$, $h_{6} = \frac{\sqrt{2}}{2^{5}}$, $h_{7} = \frac{\sqrt{2}}{2^{6}} \approx 0.0220$ and $h_{8} = \frac{\sqrt{2}}{2^{7}} \approx 0.0110$, respectively. Note that for the \textsf{mesh$_{7}$}, we have approximately $16\,129$ unknowns (corresponding to $n = 128$). 
\begin{figure}[h]
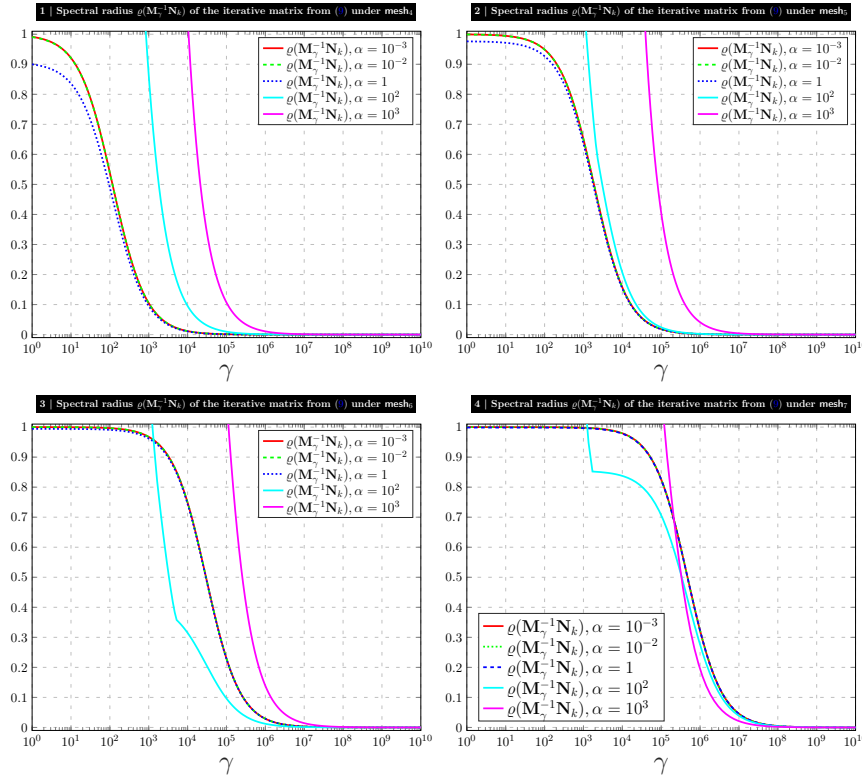
  
    \centering  
    \scalebox{0.45}{
        \input{RayonProb1mesh4.tex}     %
        \input{RayonProb1mesh5.tex}}    
        \scalebox{0.45}{
        \input{RayonProb1mesh6.tex}     
        \input{RayonProb1mesh7.tex}}
    \caption{\footnotesize The spectral radius $\varrho(\mathbf{M}^{-1}_{\gamma}\mathbf{N}_{k})$ of the stabilized method \eqref{EqP1} (or \eqref{EqP1Gen}) applied to the linear system derived from the finite element approximation of the problem \eqref{Equat1}-\eqref{Equat2}, represented with respect to the stabilization parameter $\gamma$ for the different values of the fixed stepsize $\alpha_{k} = \alpha = 10^{i}, i = -3, -2, 0, 2, 3$, and for the different meshes \textsf{mesh$_{4}$} in {\tiny\colorbox{black}{\color{white}1}}, \textsf{mesh$_{5}$} in {\tiny\colorbox{black}{\color{white}2}}, \textsf{mesh$_{6}$} in {\tiny\colorbox{black}{\color{white}3}} and  \textsf{mesh$_{7}$} in {\tiny\colorbox{black}{\color{white}4}}.}
    \label{RayonEF1}
\end{figure}

We present in Figure \ref{RayonEF1}, the spectral radius $\varrho(\mathbf{M}^{-1}_{\gamma}\mathbf{N}_{k})$ of the stabilized gradient method \eqref{EqP1} as a function of the parameter $\gamma \in [10^{0}, 10^{10}]$, for each of the stepsizes $\alpha = 10^{i}, i = -3, -2, 0, 2, 3$ and the meshes \textsf{mesh$_{i}$}, $i = 4, 5, 6, 7$. For $\alpha = 10^{i}, i = -3, -2, 0$, we note that this method converges from any initial value $\bm{x}^{[0]}_{\gamma}$. Whereas for $\alpha = 10^{i}, i = 2, 3$, the convergence is reached when the stabilization parameter $\gamma$ is greater from a certain threshold. 

In Table \ref{TabEF0}, we present numerical results computed from the stabilized gradient method \eqref{EqP1} using the finite element matrix $\mathbf{A}_{6}$ and right-hand side $\mathbf{b}_{6}$ obtained under the mesh \textsf{mesh$_{6}$}. We obtain accurate and very stable results with respect to the fixed mesh size $h_{6}$, and we note that the stabilized gradient method behaves like a direct method as soon as $\gamma > \frac{1}{\varepsilon}$. 
\begin{table}[h]
\begin{center} 
\begin{tikzpicture}
\draw[gray] [white, fill=gray!15.5]
    (0.0,6).. controls (15.5,6) ..(15.5,6)
    .. controls (15.5, 7.25) ..(15.5, 7.25)
    .. controls (0.0, 7.25) ..(0.0, 7.25)
    .. controls (0.0,6) ..(0.0,6);

\draw[gray] [white, fill=gray!10]
    (0.0,4.0).. controls (15.5,4.0) ..(15.5,4.0)
    .. controls (15.5, 5) .. (15.5, 5)
    .. controls (0.0, 5)  .. (0.0, 5)
    .. controls (0.0,4.0) .. (0.0,4.0);
\draw[gray] [white, fill=gray!10]
    (0.0,2.0).. controls (15.5,2.0) ..(15.5,2.0)
    .. controls (15.5, 3) .. (15.5, 3)
    .. controls (0.0, 3)  .. (0.0, 3)
    .. controls (0.0,2.0) .. (0.0,2.0);
\draw[gray] [white, fill=gray!10]
    (0.0,0.0).. controls (15.5,0.0) ..(15.5,0.0)
    .. controls (15.5, 1) .. (15.5, 1)
    .. controls (0.0, 1)  .. (0.0, 1)
    .. controls (0.0,0.0) .. (0.0,0.0);

\draw[, very thin] (1.,0) -- (1.,7.25);
\draw[, very thin] (3.5,0) -- (3.5,7.25);
\draw[, very thin] (6.,0) -- (6.,7.25);
\draw[, very thin] (8.5,0) -- (8.5,7.25);
\draw[, very thin] (11.,0) -- (11.,7.25);
\draw[, very thin] (13.5,0) -- (13.5,7.25);

\draw[, very thin] (0,0) -- (15.5,0);
\draw[, very thin] (0,1) -- (15.5,1);
\draw[, very thin] (0,2)   -- (15.5,2);
\draw[, very thin] (0,3)   -- (15.5,3);
\draw[, very thin] (0,4)   -- (15.5,4);
\draw[, very thin] (0,5) -- (15.5,5);
\draw[, very thin] (0,6) -- (15.5,6);

\draw[black](0.4, 6.9)  node [rotate=55, scale=0.9] {\text{Stab.}};
\draw[black](0.56, 6.67)  node [rotate=55, scale=0.8] {\text{Parameter}};
\draw[black](0.8, 6.5)  node [rotate=55, scale=0.9] {$(\gamma)$};
\draw[black](0.5, 5.5)  node [rotate=0, scale=1] {$10^{4}$};
\draw[black](0.5, 4.5)  node [rotate=0, scale=1] {$10^{6}$};
\draw[black](0.5, 3.5)  node [rotate=0, scale=1] {$10^{8}$};
\draw[black](0.5, 2.5)  node [rotate=0, scale=1] {$10^{10}$};
\draw[black](0.5, 1.5)  node [rotate=0, scale=1] {$10^{15}$};
\draw[black](0.5, 0.5)  node [rotate=0, scale=1] {$10^{20}$};

\draw[black](2.25, 6.7)   node [rotate=0, scale=1] {$\|\mathbf{A}\bm{x}^{[k]}_{\gamma,h_{6}} - \bm{b}\|_{n}$};
\draw[black](2.25, 5.5)  node [rotate=0, scale=1] {$1.030300e^{-02}$};
\draw[black](2.25, 4.5)  node [rotate=0, scale=1] {$5.681840e^{-03}$};
\draw[black](2.25, 3.5)  node [rotate=0, scale=1] {$5.849127e^{-05}$};
\draw[black](2.25, 2.5)  node [rotate=0, scale=1] {$5.851086e^{-07}$};
\draw[black](2.25, 1.5)  node [rotate=0, scale=1] {$2.233742e^{-11}$};
\draw[black](2.25, 0.5)  node [rotate=0, scale=1] {$2.740189e^{-11}$};

\draw[black](4.75, 6.6)   node [rotate=0, scale=1] {$\frac{\|\mathbf{A}\bm{x}^{[k]}_{\gamma,h_{6}} - \bm{b}\|_{n}}{\|\mathbf{A}\bm{x}^{[0]}_{\gamma} - \bm{b}\|_{n}}$};
\draw[black](4.75, 5.5)  node [rotate=0, scale=1] {$2.041260e^{-03}$};
\draw[black](4.75, 4.5)  node [rotate=0, scale=1] {$1.125703e^{-03}$};
\draw[black](4.75, 3.5)  node [rotate=0, scale=1] {$1.158846e^{-05}$};
\draw[black](4.75, 2.5)  node [rotate=0, scale=1] {$1.159234e^{-07}$};
\draw[black](4.75, 1.5)  node [rotate=0, scale=1] {$4.425556e^{-12}$};
\draw[black](4.75, 0.5)  node [rotate=0, scale=1] {$5.428944e^{-12}$};

\draw[black](7.25, 6.7)   node [rotate=0, scale=1] {$\|\bm{u}^{\star} - \bm{x}^{[k]}_{\gamma,h_{6}}\|_{n}$};
\draw[black](7.25, 5.5)  node [rotate=0, scale=1] {$1.823648e^{+00}$};
\draw[black](7.25, 4.5)  node [rotate=0, scale=1] {$1.060131e^{+00}$};
\draw[black](7.25, 3.5)  node [rotate=0, scale=1] {$4.056218e^{-01}$};
\draw[black](7.25, 2.5)  node [rotate=0, scale=1] {$4.055075e^{-01}$};
\draw[black](7.25, 1.5)  node [rotate=0, scale=1] {$4.055076e^{-01}$};
\draw[black](7.25, 0.5)  node [rotate=0, scale=1] {$4.055076e^{-01}$};

\draw[black](9.75, 6.6)   node [rotate=0, scale=1] {$\frac{\|\bm{u}^{\star} - \bm{x}^{[k]}_{\gamma,h_{6}}\|_{n}}{\|\bm{u}^{\star}\|_{n}}$};
\draw[black](9.75, 5.5)  node [rotate=0, scale=1] {$3.897713e^{-02}$};
\draw[black](9.75, 4.5)  node [rotate=0, scale=1] {$2.265835e^{-02}$};
\draw[black](9.75, 3.5)  node [rotate=0, scale=1] {$8.669420e^{-03}$};
\draw[black](9.75, 2.5)  node [rotate=0, scale=1] {$8.666977e^{-03}$};
\draw[black](9.75, 1.5)  node [rotate=0, scale=1] {$8.666979e^{-03}$};
\draw[black](9.75, 0.5)  node [rotate=0, scale=1] {$8.666979e^{-03}$};

\draw[black](12.25, 6.85)   node [rotate=0, scale=1] {\text{$\#$ Iterations}};
\draw[black](12.25, 6.3)   node [rotate=0, scale=1] {$[k]$};
\draw[black](12.25, 5.5)  node [rotate=0, scale=1] {$10$};
\draw[black](12.25, 4.5)  node [rotate=0, scale=1] {$1$};
\draw[black](12.25, 3.5)  node [rotate=0, scale=1] {$1$};
\draw[black](12.25, 2.5)  node [rotate=0, scale=1] {$1$};
\draw[black](12.25, 1.5)  node [rotate=0, scale=1] {$1$};
\draw[black](12.25, 0.5)  node [rotate=0, scale=1] {$1$};

\draw[black](14.5, 6.85)   node [rotate=0, scale=1] {\text{CPU times}};
\draw[black](14.5, 6.35)   node [rotate=0, scale=1] {\text{($s$)}};
\draw[black](14.5, 5.5)  node [rotate=0, scale=1] {$3.1616$};
\draw[black](14.5, 4.5)  node [rotate=0, scale=1] {$0.3995$};
\draw[black](14.5, 3.5)  node [rotate=0, scale=1] {$0.3996$};
\draw[black](14.5, 2.5)  node [rotate=0, scale=1] {$0.4300$};
\draw[black](14.5, 1.5)  node [rotate=0, scale=1] {$0.4053$};
\draw[black](14.5, 0.5)  node [rotate=0, scale=1] {$0.4067$};

\end{tikzpicture} 
\caption{Numerical results of the method \eqref{EqP1} for different values of the stabilization parameter $\gamma$, under the fixed \textsf{mesh$_{6}$} with $n = 3\,969$, the stepsize $\alpha = 1$, the maximum iterations number $k_{\text{max}} = 100$ and the tolerance $\varepsilon = h_{6} \backsimeq 0.044$ in \eqref{CTR2}. The data associated to the \textsf{mesh$_{6}$}, namely the matrix $\mathbf{A}_{6}$ and the right hand side $\bm{b}_{6}$, are computed from the finite element approximation of the reaction-diffusion problem \eqref{Equat1}-\eqref{Equat2}, where $\textbf{cond}(\mathbf{A}_{6}) = 1.9728e^{+03}$.}
\label{TabEF0}
\end{center}
\end{table}
\begin{table}[h]
\begin{center} 
\begin{tikzpicture}
\draw[gray] [white, fill=gray!15.5]
    (0.0,6).. controls (16.5,6) ..(16.5,6)
    .. controls (16.5, 7.25) ..(16.5, 7.25)
    .. controls (0.0, 7.25) ..(0.0, 7.25)
    .. controls (0.0,6) ..(0.0,6);

\draw[gray] [white, fill=gray!10]
    (0.0,4.0).. controls (16.5,4.0) ..(16.5,4.0)
    .. controls (16.5, 5) .. (16.5, 5)
    .. controls (0.0, 5)  .. (0.0, 5)
    .. controls (0.0,4.0) .. (0.0,4.0);
\draw[gray] [white, fill=gray!10]
    (0.0,2.0).. controls (16.5,2.0) ..(16.5,2.0)
    .. controls (16.5, 3) .. (16.5, 3)
    .. controls (0.0, 3)  .. (0.0, 3)
    .. controls (0.0,2.0) .. (0.0,2.0);

\draw[, very thin] (1.,1) -- (1.,7.25);
\draw[, very thin] (3.5,1) -- (3.5,7.25);
\draw[, very thin] (6.,1) -- (6.,7.25);
\draw[, very thin] (8.5,1) -- (8.5,7.25);
\draw[, very thin] (11.,1) -- (11.,7.25);
\draw[, very thin] (13.5,1) -- (13.5,7.25);
\draw[, very thin] (14.85,1) -- (14.85,7.25);

\draw[, very thin] (0,1) -- (16.5,1);
\draw[, very thin] (0,2)   -- (16.5,2);
\draw[, very thin] (0,3)   -- (16.5,3);
\draw[, very thin] (0,4)   -- (16.5,4);
\draw[, very thin] (0,5) -- (16.5,5);
\draw[, very thin] (0,6) -- (16.5,6);

\draw[black](0.36, 6.78)  node [rotate=55, scale=0.8] {\text{Meshes}};
\draw[black](0.53, 6.6)  node [rotate=55, scale=0.8] {\text{size}};
\draw[black](0.8, 6.4)  node [rotate=60, scale=0.8] {$(h_{i})$};
\draw[black](0.5, 5.5)  node [rotate=0, scale=1] {$h_{4}$};
\draw[black](0.5, 4.5)  node [rotate=0, scale=1] {$h_{5}$};
\draw[black](0.5, 3.5)  node [rotate=0, scale=1] {$h_{6}$};
\draw[black](0.5, 2.5)  node [rotate=0, scale=1] {$h_{7}$};
\draw[black](0.5, 1.5)  node [rotate=0, scale=1] {$h_{8}$};

\draw[black](2.25, 6.7)   node [rotate=0, scale=1] {$\|\mathbf{A}\bm{x}^{[k]}_{\gamma,h_{i}} - \bm{b}\|_{n}$};
\draw[black](2.25, 5.5)  node [rotate=0, scale=1] {$1.682288e^{-13}$};
\draw[black](2.25, 4.5)  node [rotate=0, scale=1] {$1.692714e^{-12}$};
\draw[black](2.25, 3.5)  node [rotate=0, scale=1] {$2.233742e^{-11}$};
\draw[black](2.25, 2.5)  node [rotate=0, scale=1] {$2.056490e^{-10}$};
\draw[black](2.25, 1.5)  node [rotate=0, scale=1] {$3.363334e^{-10}$};

\draw[black](4.75, 6.6)   node [rotate=0, scale=1] {$\frac{\|\mathbf{A}\bm{x}^{[k]}_{\gamma,h_{i}} - \bm{b}\|_{n}}{\|\mathbf{A}\bm{x}^{[0]}_{\gamma} - \bm{b}\|_{n}}$};
\draw[black](4.75, 5.5)  node [rotate=0, scale=1] {$8.332645e^{-15}$};
\draw[black](4.75, 4.5)  node [rotate=0, scale=1] {$1.676837e^{-13}$};
\draw[black](4.75, 3.5)  node [rotate=0, scale=1] {$4.425556e^{-12}$};
\draw[black](4.75, 2.5)  node [rotate=0, scale=1] {$8.148734e^{-11}$};
\draw[black](4.75, 1.5)  node [rotate=0, scale=1] {$2.665403e^{-10}$};

\draw[black](7.25, 6.7)   node [rotate=0, scale=1] {$\|\bm{u}^{\star} - \bm{x}^{[k]}_{\gamma,h_{i}}\|_{n}$};
\draw[black](7.25, 5.5)  node [rotate=0, scale=1] {$1.820886e^{+00}$};
\draw[black](7.25, 4.5)  node [rotate=0, scale=1] {$8.295110e^{-01}$};
\draw[black](7.25, 3.5)  node [rotate=0, scale=1] {$4.055076e^{-01}$};
\draw[black](7.25, 2.5)  node [rotate=0, scale=1] {$2.016228e^{-01}$};
\draw[black](7.25, 1.5)  node [rotate=0, scale=1] {$1.006708e^{-01}.$};

\draw[black](9.75, 6.6)   node [rotate=0, scale=1] {$\frac{\|\bm{u}^{\star} - \bm{x}^{[k]}_{\gamma,h_{i-1}}\|_{n}}{\|\bm{u}^{\star} - \bm{x}^{[k]}_{\gamma,h_{i}}\|_{n}}$};
\draw[black](9.75, 5.5)  node [rotate=0, scale=1] {$-$};
\draw[black](9.75, 4.5)  node [rotate=0, scale=1] {$2.1951$};
\draw[black](9.75, 3.5)  node [rotate=0, scale=1] {$2.0456$};
\draw[black](9.75, 2.5)  node [rotate=0, scale=1] {$2.0112$};
\draw[black](9.75, 1.5)  node [rotate=0, scale=1] {$2.0028$};

\draw[black](12.25, 6.7)   node [rotate=0, scale=1] {$\frac{\|\bm{u}^{\star} - \bm{x}^{[k]}_{\gamma,h_{i}}\|_{n}}{\|\bm{u}^{\star}\|_{n}}$};
\draw[black](12.25, 5.5)  node [rotate=0, scale=1] {$1.556736e^{-01}$};
\draw[black](12.25, 4.5)  node [rotate=0, scale=1] {$3.545856e^{-02}$};
\draw[black](12.25, 3.5)  node [rotate=0, scale=1] {$8.666979e^{-03}$};
\draw[black](12.25, 2.5)  node [rotate=0, scale=1] {$2.154659e^{-03}$};
\draw[black](12.25, 1.5)  node [rotate=0, scale=1] {$5.379134e^{-04}$};

\draw[black](14.2, 7.0)   node [rotate=0, scale=1] {\text{CPU}};
\draw[black](14.2, 6.6)   node [rotate=0, scale=1] {\text{times}};
\draw[black](14.2, 6.20)   node [rotate=0, scale=1] {\text{($s$)}};
\draw[black](14.2, 5.5)  node [rotate=0, scale=1] {$0.0064$};
\draw[black](14.2, 4.5)  node [rotate=0, scale=1] {$0.0208$};
\draw[black](14.2, 3.5)  node [rotate=0, scale=1] {$0.3932$};
\draw[black](14.2, 2.5)  node [rotate=0, scale=1] {$4.6163$};
\draw[black](14.2, 1.5)  node [rotate=0, scale=1] {$171.70$};

\draw[black](15.7, 6.6)   node [rotate=0, scale=1] {\text{Cond($\mathbf{A}$)}};
\draw[black](15.7, 5.5)  node [rotate=0, scale=1] {$1.231e^{+02}$};
\draw[black](15.7, 4.5)  node [rotate=0, scale=1] {$4.930e^{+02}$};
\draw[black](15.7, 3.5)  node [rotate=0, scale=1] {$1.972e^{+03}$};
\draw[black](15.7, 2.5)  node [rotate=0, scale=1] {$7.891e^{+03}$};
\draw[black](15.7, 1.5)  node [rotate=0, scale=1] {$3.157e^{+04}$};

\end{tikzpicture} 
\caption{Numerical results of the method \eqref{EqP1} for different meshes \textsf{mesh$_{i}$}, $i = 4, 5, 6, 7, 8, 9$, under the fixed stabilization parameter $\gamma = 10^{15}$, the stepsize $\alpha = 1$, the maximum iterations number $k_{\text{max}} = 100$ and the tolerance $\varepsilon_{i} = h_{i}$, $i = 4, 5, 6, 7, 8, 9$, in \eqref{CTR2}. The data associated to the mesehes \textsf{mesh$_{i}$}, namely the matrix $\mathbf{A}$ and the right hand side $\bm{b}$, are computed from the finite element approximation of the reaction-diffusion problem \eqref{Equat1}-\eqref{Equat2}.}
\label{TabEF1}
\end{center}
\end{table}

We also present in Table \ref{TabEF1} the same results, but with respect to the meshes size $h_{i}, i = 4, 5, 6, 7, 8$ and the fixed stabilization parameter $\gamma = 10^{15}$. We also note a very stable solution $\bm{x}^{[k]}_{\gamma}$ despite the large values of the parameter $\gamma$. This is explained by the moderate conditioning of the matrices we face here (see the last column of the Table \ref{TabEF1}). Through these computations, the number of iterations made is always one for each mesh and is no longer proportional to the condition number of the matrix $\mathbf{A}_{i}$.
\section{General Conclusion}
We proposed through this work a stabilization of the gradient method by replacing the gradient algorithm \eqref{Eq2} by the iterative method \eqref{EqP1}. This stabilization method is based on the idea to build an iterative scheme from the gradient method, by gradually controlling the residual as the iterations go along. We established theoretical convergence results in Corollary \ref{CorMSN12} and made numerical computations in the Section \ref{sec5} that illustrated our analysis. These established results are made under the nonsingular matrix property alone, and the performed numerical tests take into account the both symmetric and non-symmetric matrices.

We first analyse the filtered SVD expansion of the stabilized iterative solution $\bm{x}^{[k]}_{\gamma}$ in \eqref{SVDSol}, and note that as the stabilization parameter $\gamma$ increases, the filter factors \eqref{FilFac} effectively include more and more SVD components in the iterations $\bm{x}^{[k]}_{\gamma}$ by shifting toward the unit value. The limit expansion obtained in \eqref{SVD1} (namely the naïve solution reached when $\gamma \to +\infty$), shows that under ill-conditioning problem we should not take the stabilization parameter as large as we want since it is known that this solution is very inaccurate. This explained the oscillating solution we observe in sub-subsection \ref{NEx3SSP} under severely ill-conditioned systems (see Table \ref{Tabshaw1}, Table \ref{Tabheat1}, Table \ref{Tabgravity1}, Figure \ref{Figshaw1}, Figure \ref{Figheat1}, and Figure \ref{Figgravity1}, when the stabilization parameter is fixed to $\gamma = 10^{14}, 10^{15}, 10^{16}$). This hence proves that only controlling the residual (as we did in the stabilized gradient method \eqref{EqP1}), may not provide accurate approximations under severely ill-conditioning context. However, with large stabilization parameter (namely when $\gamma$ is around $10^{10}$), we obtain through severely ill-conditioned systems analysed in sub-subsection \ref{NEx3SSP}, efficient and accurate approximations from the stabilized solution $\bm{x}^{[k]}_{\gamma}$. Through these examples, we only faced errors from the quadrature in the evaluation of the coefficients of the matrix $\mathbf{A}$ and the right-hand side $\mathbf{b}$, and the rounding errors. 

The stabilized gradient method has spared us the intrinsic difficulty associated with determining the stepsize $\alpha$ in the gradient method. We first observed in \eqref{Reslt1Fin2} and \eqref{Reslt1FinBis1} that the convergence is speeded up by taking the stabilization parameter $\gamma$ as large as possible, and is no longer based on the choix of the stepsize $\alpha$. For well-conditioned linear systems in sub-sections \ref{NEx1} and \ref{NEx2} whose the matrix is symmetric (or not), the numerical tests have confirmed this result. And for ill-conditioned linear systems, we observe the same behaviour as long as the ill-conditioning is not severe as faced under the {\bf shaw}, {\bf heat} and {\bf gravity} test problems (in sub-subsection \ref{NEx3SSP}). The convergence of the stabilized gradient method \eqref{EqP1} is as faster than the stabilization parameter $\gamma$ is large (it sometimes behaves as a direct method). Moreover, the convergence of the stabilized method \eqref{EqP1} no longer depends on the eigenvalues of the matrix $\mathbf{A}$ (the structure of the matrix $\mathbf{A}$). It is only guided by the choice of the stabilization parameter $\gamma$, independently of the linear system data $\mathbf{A}$ and $\mathbf{b}$.

As we have seen through the theoretical results \eqref{Reslt1Fin2}-\eqref{Reslt1FinBis1} and the previous numerical examples, the value of the stepsize $\alpha$ has a limit influence on the convergence of the stabilized gradient method \eqref{EqP1} (meaning that we may fix the stepsize to neutral value $\alpha = 1$). What really matter is the value of the stabilization parameter $\gamma$, which should be fixed as large as possible in order to quickly reduce the residual and hence to make less iterations. 

In forthcoming works, it will be interesting to remedy to the oscillating solutions observed under severely ill-conditioned systems of the sub-subsection \ref{NEx3SSP}. The idea which may be tested is, in addition to the stabilization made here, to add some regularization of the solution as in the direct solving method proposed in \cite{Dione24}.


\begin{thebibliography}{9}
\bibitem{Akaike59}
Akaike H., \emph{On a successive transformation of probability distribution and its application to the analysis of the optimum gradient method}, Annals of the Institute of Statistical Mathematics, Tokyo, 11, 1-16 (1959).

\bibitem{AllaireKaber08}
Allaire G. and Kaber S. M., \emph{Numerical Linear Algebra}, Springer (2008).

\bibitem{AllaireKaber07}
Allaire G. and Kaber S. M., \emph{Numerical Linear Algebra. Solutions Manual for Instructors}, Springer (2007).


\bibitem{Armijo66}
Armijo L., \emph{Minimization of functions having Lipschitz continuous first partial derivatives}, Pacific Journal of mathematics, 16, 1-3 (1966).


\bibitem{BarzilaiBorwein88}
Barzilai J. and Borwein J. M., \emph{Two point stepsize gradient methods}, IMA Journal of Numerical Analysis, 8, 141-148 (1988).

\bibitem{BazaraaSheraliShetty06}
Bazaraa M. S., Sherali H. D. and Shetty C. M., \emph{Nonlinear programming : theory and algorithms}, Third Edition, John Wiley \& Sons, Inc. (2006).

\bibitem{Beck14}
Beck A., \emph{Introduction to nonlinear optimization : theory, algorithms, and applications with MATLAB}, MOS-SIAM series on optimization (2014).

\bibitem{BirginMartinezRaydan00}
Birgin E. G., Martinez J. M. and Raydan M., \emph{Nonmonotone spectral projected gradient methods on convex sets}, SIAM Journal on Optimization, 10, 1196-1211 (2000).


\bibitem{BjÄorck96}
Björck Å., \emph{Numerical methods for least squares problems}, Society for Industrial and Applied Mathematics (1996).






\bibitem{Cauchy47}
Cauchy A. L., \emph{Méthode générale pour la résolution des systemes d'équations simultanées}, Comptes Rendus Sciences Paris, 25, 536-538 (1847).


\bibitem{DaiYuanYuan02}
Dai Y. H., Yuan J. Y. and Yuan Y., \emph{Modified two-point stepsize gradient methods for unconstrained optimization}, Computational Optimization and Applications, 22, 103-109 (2002).

\bibitem{DaiYuan03}
Dai Y. H. and Yuan Y., \emph{Alternate minimization gradient method}, IMA Journal of Numerical Analysis, 23, 377-393 (2003).


\bibitem{DennisSchnabel83}
Dennis J. E. and Schnabel R. B., \emph{Numerical Methods for Unconstrained Optimization and Nonlinear Equations}, Prentice-Hall, Englewoods Cliffs, New Jersey (1983).

\bibitem{Dione24}
Dione, I., \emph{Regularization of Discrete Ill-Conditioned Problems Done Right - I}, Preprint submitted to Applied Mathematics and Computation (2024).

\bibitem{ELDÉN95}
Eldén L., \emph{Numerical solution of the sideways heat equation}, Inverse Problem 11, 913-923 (1995).



\bibitem{Fletcher05}
Fletcher R., \emph{On the Barzilar-Borwein method}, The International Workshop on Optimization and Control with Applications, Springer, 235-256 (2005).

\bibitem{Fletcher90}
Fletcher R., \emph{Low storage methods for unconstrained optimization}, Lectures in Applied Mathematics, 26, 165-179 (1990).

\bibitem{Fletcher87}
Fletcher R., \emph{Practical Methods of Optimization}, Wiley, New York (1987).


\bibitem{ForsytheMotzkin51}
Forsythe G. E., Motzkin T. S., \emph{Asymptotic properties of the optimum gradient method}, Bull.
Am. Soc. 57, 183 (1951).

\bibitem{FriedlanderMartinezRaydan95}
Friedlander A., Martinez J. M. and Raydan M., \emph{A new method for large-scale box constrained convex quadratic minimization problems}, Optimization Methods and Software, 5, 57-74 (1995).



\bibitem{Goldstein65}
Goldstein, A. A., \emph{On steepest descent}, SIAM Journal Control, 3, 147-151 (1965).

\bibitem{GolubLoan13}
Golub G. H. and Van Loan C. F., \emph{Matrix Computations}, Fourth Edition, The Johns Hopkins University Press
Baltimore (2013).


\bibitem{GrivaNashSofer09}
Griva I., Nash S. G. and Sofer A., \emph{Linear and Nonlinear Optimization}, $2^{nd}$ edition, Society for Industrial and Applied Mathematics (2009).


\bibitem{Groetsch84}
Groetsch C. W., \emph{The theory of Tikhonov regularization for Fredholm equations of the first kind}, Longman Higher Education (1984).


\bibitem{Hackbusch10}
Hackbusch W., \emph{Iterative Solution of Large Sparse Systems of Equations}, Springer, Second Edition (2010).

\bibitem{Hanke95}
Hanke M., \emph{Conjugate gradient type methods for ill-posed problems}, Longman Scientific and Technical, Wiley (1995).

\bibitem{HansenBook10}
Hansen P. C., \emph{Discrete inverse problems: insight and algorithms}, Society for Industrial and Applied Mathematics (2010).

\bibitem{HansenBook07}
Hansen P. C., \emph{REGULARIZATION TOOLS: A Matlab package for analysis and solution of discrete ill-posed problems}, Numerical Algorithms, 46, 189-194 (2007). 


\bibitem{HansenBook98}
Hansen P. C., \emph{Rank-deficient and discrete ill-posed problems: numerical aspects of linear inversion}, SIAM monographs on mathematical modeling and computation (1998).


\bibitem{HestenesStiefel52}
Hestenes M. R. and Stiefel E. L., \emph{Methods of conjugate gradients for solving linear systems}, Journal of Research of the National Bureau of Standards, 49, 409-436 (1952).

\bibitem{Humphrey66}
Humphrey W. E., \emph{A general minimising routine - minfun. In: Lavi, A., Vogl, T.P. (eds.) Recent
Advances in Optimisation Techniques.}, Wiley, New York (1966).





\bibitem{Lemaréchal81}
Lemaréchal C., \emph{Aview of line search}, In: Auslander, A., Oettli, W., Stoer, J. (eds.) Optimization
and Optimal Control, pp. 59–78. Springer, Berlin (1981).


\bibitem{MoréThuente90}
Moré J. J. and Thuente D. J., \emph{On line search algorithms with guaranteed sufficient decrease}, Mathematics
and Computer Science Division Preprint MCS-P153-0590, Argonne National Laboratory, Argonne (1990).




\bibitem{leary09}
O'leary D. P., \emph{Scientific computing with case studies}, Society for Industrial and Applied Mathematics (2009).



\bibitem{Powell76}
Powell M. J. D., \emph{Some global convergence properties of a variable-metric algorithm for minimization without exact line searches}, SIAM, 9, 53-72 (1976).

\bibitem{PotraShi95}
Potra F. A. and Shi Y., \emph{Efficient line search algorithm for unconstrained optimization}, Journal of Optimization Theory and Applications, 85, 677-704 (1995).


\bibitem{QuarteroniSaccoSaleri07}
Quarteroni A., Sacco R. and Saleri F., \emph{Numerical Mathematics}, Springer, Second Edition (2007).

\bibitem{Raydan93}
Raydan M., \emph{On the Barzilai and Borwein choice of steplength for the gradient method}, IMA Journal of Numerical Analysis, 13, 321-326 (1993).

\bibitem{Raydan97}
Raydan M., \emph{The Barzilai and Borwein gradient method for the large scale unconstrained minimization problem}, SIAM Journal on Optimization 7, 26-33 (1997).


\bibitem{Schinzinger66}
Schinzinger, R., \emph{Optimization in electromagnetic system design. In: Lavi, A., Vogl, T.P. (eds.) Recent Advances in Optimisation Techniques}, Wiley, New York (1966).


\bibitem{SibonyMardon88}
Sibony M., Mardon J. C., \emph{Analyse numérique I. Systèmes linéaires et non linéaires}, Hermann (1988).





\bibitem{Wolfe68}
Wolfe P., \emph{Convergence conditions for ascent methods}, SIAM Review, 11, 226-235 (1968).

\bibitem{Yuan06}
Yuan Ya-xiang, \emph{A new stepsize for the steepest descent method}, Journal of Computational Mathematics, 24, 149-156 (2006).

\end{thebibliography}
\end{document}